\documentclass[final]{siamltex}

\usepackage{amssymb,latexsym,amsmath}
\usepackage{epsfig}
\usepackage{caption}
 
 \usepackage{epsfig}

\usepackage{makeidx}

\allowdisplaybreaks

 \usepackage{graphicx,fancybox,latexsym,epsfig}
\usepackage{fancyhdr,amsmath,times,amsxtra,amssymb}
\usepackage{color}

\usepackage{fancybox}
\usepackage{tikz}

\usepackage{amsmath}
\usepackage{amsfonts}
\usepackage{array}
\usepackage{dsfont}
\usepackage{hyperref}
\usepackage{amssymb}
\usepackage{bbold}
\usepackage{standalone}
\usepackage{accents}
\usepackage{enumitem}
\usepackage{tikz}
\usepackage{pgfplots}
\pgfplotsset{compat=1.6}
\usepgfplotslibrary{groupplots}
\usepackage{cancel}
 \usepackage{color}

\allowdisplaybreaks

\newtheorem{assumption}{Assumption}[section]

  \newtheorem{remark}{Remark}
  \newtheorem{example}{Example}

\usepackage{amssymb,latexsym,amsmath}
\usepackage{mathrsfs}
\usepackage{epsfig}
\usepackage{caption}

\newcommand\scalemath[2]{\scalebox{#1}{\mbox{\ensuremath{\displaystyle #2}}}}

\pgfplotsset{soldot/.style={color=blue,only marks,mark=*}}
\pgfplotsset{holdot/.style={color=blue,fill=white,only marks,mark=*}}
\fussy

\allowdisplaybreaks

\begin{document}
\sloppy
\title{Optimality of Independently Randomized Symmetric Policies for Exchangeable Stochastic Teams with Infinitely Many Decision Makers
\thanks{This research was supported by the Natural Sciences and Engineering Research Council (NSERC) of Canada.  Sina Sanjari and Serdar Y\"uksel are with the Department of Mathematics and Statistics,
     Queen's University, Kingston, ON, Canada,
     Email: \{16sss3,yuksel@queensu.ca\}}.
     Naci Saldi is with the Department of Natural and Mathematical Sciences, {\"{O}}zye\u{g}in University, Cekmekoy, Turkey, Email: \{naci.saldi@ozyegin.edu.tr\}}
\author{Sina Sanjari, Naci Saldi and Serdar Y\"uksel}
\maketitle
\tableofcontents

\begin{abstract}{
We study stochastic team (known also as decentralized stochastic control or identical interest stochastic dynamic game) problems with large or countably infinite number of decision makers, and characterize existence and structural properties for (globally) optimal policies. We consider both static and dynamic non-convex team problems where the cost function and dynamics satisfy an exchangeability condition. {To arrive at existence and structural results on optimal policies, we first introduce a topology on control policies, which involves various relaxations given the decentralized information structure. This is then utilized to arrive at a de Finetti type representation theorem for exchangeable policies.} This leads to a representation theorem for policies which admit an infinite exchangeability condition. For a general setup of stochastic team problems with $N$ decision makers, under exchangeability of observations of decision makers and the cost function, we show that without loss of global optimality, the search for optimal policies can be restricted to those that are $N$-exchangeable. Then, by extending $N$-exchangeable policies to infinitely-exchangeable ones, establishing a convergence argument for the induced costs, and using the presented de Finetti type theorem, we establish the existence of an optimal decentralized policy for static and dynamic teams with countably infinite number of decision makers, which turns out to be symmetric (i.e., identical) and randomized. In particular, unlike prior work, convexity of the cost in policies is not assumed. Finally, we show near optimality of symmetric independently randomized policies for finite $N$-decision maker team problems and thus establish approximation results for $N$-decision maker weakly coupled stochastic teams.}
\end{abstract}




\section{Introduction}

\label{sec:intro}

Stochastic team problems consist of a collection of decision makers or agents acting together to optimize a common cost function, but not necessarily sharing all the available information. %
At each time stage, each decision maker only has partial access to the global information which is defined by the \textit{information structure} (IS) of the problem \cite{wit75}. When there is a pre-defined order according to which the decision makers act then the team is called a \textit{sequential team}. For sequential teams, if each agent's information depends only on primitive random variables, the team is \textit{static}. If at least one agent's information is affected by an action of another agent, the team is said to be \textit{dynamic}. 

In this paper, we study stochastic team problems with a large but finite, and countably infinite number of decision makers. We characterize existence and structural properties of (globally) optimal policies in such problems. While teams can be at first sight viewed as a narrow class of (identical interest) stochastic dynamic games, when viewed as a generalization of classical single decision maker (DM) stochastic control, they are quite general with increasingly common applications involving many areas of applied mathematics such as decentralized stochastic control \cite{sandell1978survey, ho1980team, CDCTutorial, arrow1979allocation}, networked control \cite{ho1980team, hespanha2007survey}, communication networks \cite{hespanha2007survey}, cooperative systems \cite{mar55, Radner, mcguire1961some, beckmann1958decision}, large sensor networks \cite{tsitsiklis1988decentralized}, and energy, or more generally, smart grid design \cite{sandell1978survey, davison1973optimal}.

{\bf Connections to convex stochastic teams.} For teams with finitely many decision makers, Marschak \cite{mar55} studied static teams and Radner \cite{Radner} established connections between person-by-person optimality, stationarity, and team-optimality. Radner's results were generalized in \cite{KraMar82} by relaxing optimality conditions. A summary of these results is that in the context of static team problems, the convexity of the cost function, subject to minor regularity conditions, suffices for the global optimality of person-by-person-optimal solutions. In the particular case for LQG (Linear Quadratic Gaussian) static teams, this result leads to the optimality of linear policies \cite{Radner}, which also applies to dynamic LQG problems under partially nested information structures \cite{HoChu}. These results are applicable to static teams with finitely many decision makers. 

In our paper, the main focus is on teams with infinitely many decision makers. In this direction, we note that in our prior works \cite{sanjari2018optimal, sanjari2019optimal}, we studied static and dynamic teams where under convexity and symmetry conditions, global optimality of the limit points of  the sequence of $N$ decision maker optimal policies was established. These works also provided existence and structural results for convex static and dynamic teams with infinitely many decision makers. We also note \cite{mahajan2013static} where LQG static teams with countably infinite number of decision makers have been studied and sufficient conditions for global optimality have been established. In our paper here, convexity is not imposed.

{\bf Connections with the literature on mean-field games/teams.} Team problems can be considered as games with identical interests. For the case with infinitely many decision makers, a related set of results involves mean-field games: mean-field games (see e.g., \cite{CainesMeanField2,CainesMeanField3,LyonsMeanField}) can be viewed as limit models of symmetric non-zero-sum non-cooperative finite player games with a mean-field interaction. We note that in team problems, person-by-person optimality (Nash equilibrium when viewed as games) does not in general imply global optimality both for $N$-decision maker teams and teams with countably infinite number of DMs. As we have mentioned, for static teams, a sufficient condition is the convexity of the cost function, subject to minor regularity conditions \cite{KraMar82}. However, mean-field teams under decentralized information structures generally correspond to dynamic team problems with non-classical information structures (an observation of a decision maker $i$ is affected by the action of a decision maker $j$ where decision maker $i$ does not have access to the observation of decision maker $j$), hence, mean-field team problems may be non-convex even under the convexity of the cost function due to non-classical information structures (see \cite[Section 3.3]{YukselSaldiSICON17} and the celebrated counterexample of Witsenhausen \cite{wit68}). Hence, person-by-person optimality  is generally inconclusive for global optimality.

The existence of equilibria has been established for mean-field games in \cite{LyonsMeanField, bardi2019non, carmona2016mean, light2018mean, lacker2015mean}. Furthermore, person-by-person optimal solutions may perform arbitrarily poorly. There have also been several studies for mean-field games where the limits of sequences of Nash equilibria have been investigated as the number of decision makers tends to infinity (see e.g., \cite{fischer2017connection, lacker2018convergence, bardi2014linear, LyonsMeanField, arapostathis2017solutions}). 

Social optima for mean-field linear quadratic Gaussian control problems under both centralized and restricted decentralized information structure have been considered in \cite{huang2012social,  huang2016linear, wang2017social, arabneydi2015team}. We also note a result in \cite{yu2021teamwise} where two large teams compete in a mean-field competition game. We refer readers to \cite{carmona2018probabilistic, caines2018peter} for a literature review and a detailed summary of some recent results on mean-field games and social optima problems. 

Some relevant studies on the existence and convergence of equilibria from the mean-field games literature are the following: In \cite{cardaliaguet2010notes}, for one-shot mean-field games, under regularity assumptions on the cost function, it has been shown that mixed Nash strategies of $N$-player symmetric games  converge through a subsequence to a limit (which is a weak-solution of the mean-field limit). In \cite{fischer2017connection},  through a concentration of measures argument, it has been shown that a subsequence of  symmetric local approximate Nash equilibria for $N$ player games converges to a solution for the mean-field game under the assumption that the normalized occupational measures converges weakly to a deterministic measure. Furthermore, using a similar method in \cite{lacker2016general}, assumptions on equilibrium policies of large population mean-field symmetric stochastic differential games have been presented to allow for convergence of asymmetric approximate Nash equilibria to a weak solution of the mean-field game \cite[Theorem 2.6]{lacker2016general} in the presence of common randomness. Using martingale methods and relaxed controls (see also \cite{fischer2017connection, lacker2016general, lacker2015mean, carmona2016mean}), an existence result and a limit theory have been established for controlled McKean-Vlasov dynamics \cite{lacker2017limit}. We note that in \cite{lacker2016general, lacker2017limit, lacker2015mean, carmona2016mean}, it has been assumed that each player has full access to the information available to all players, i.e., the controls are functions of all initial states, Wiener processes of all players, and common randomness. 

We further note that the existence results for equilibria have been established in \cite{lacker2016general, carmona2016mean, carmona2018probabilistic, fischer2017connection} where strategies of each player are assumed to be progressively measurable to the filtration generated by initial states and Wiener processes (also called \textit{open-loop} controllers in the mean-field games' literature \cite{lacker2016general, carmona2016mean, carmona2018probabilistic, fischer2017connection}). We note that in our setup under these strategies, the information structure corresponds to static information structure. The equilibria with respect to \textit{closed-loop} (in the team problem setup, with respect to dynamic information structure) is completely different since the deviating player can still influence the information  of other players and hence it can influence the average of states or actions substantially. 


{In \cite{lacker2018convergence}, under a convexity condition (which has been introduced in \cite{filippov1962certain} and also considered in \cite{lacker2017limit,  lacker2015mean}), and under the classical information structure (or full information, i.e., what would be a centralized problem in the team theoretic setup), convergence of Nash equilibria induced by (path-dependent and feedback Markovian) closed-loop controllers to a {\it weak (semi-Markov) mean-field equilibrium} has been established. We also note a result in \cite{cardaliaguet2019master} for the convergence of Markov feedback equilibria, where an infinite-dimensional partial differential equations  referred to as \textit{master equation} (obtained as a limit of Hamilton-Jacobi-Bellman systems) has been considered and its unique smooth solution has been used to show the convergence of empirical measures to the unique mean-field game equilibrium. We note that, the approach in \cite{cardaliaguet2019master} requires uniqueness of the mean-field equilibrium but the one in \cite{lacker2018convergence} applies even if mean-field equilibria are non-unique. In addition, the notion of a weak (semi-Markov) solution considered in \cite{lacker2018convergence} allows for an additional randomization in stochastic flows of measures, but under uniqueness, the limit solution becomes the unique (weak) mean-field equilibrium, and hence recovers the related convergence results in \cite{cardaliaguet2019master}. We also note that the convergence problem of Markov feedback equilibria for a finite state model with multiple mean-field equilibria has been studied in \cite{hajek2019non, bayraktar2020non, cecchin2019convergence}.} Recently, in \cite{campi2020correlated}, both a convergence  result for all correlated equilibrium solutions of discrete finite state mean-field games as limits of exchangeable correlated equilibria restricted to Markov open-loop strategies and an approximation result for $N$-player correlated equilibria have been established. For infinite horizon problems, in \cite{cardaliaguet2019example}, an example of ergodic differential games with mean-field coupling has been constructed such that limits of sequences of expected costs induced by symmetric Nash-equilibria of $N$-player games capture expected costs induced by many more Nash-equilibria policies including a mean-field equilibrium and social optimum. In \cite{lacker2018convergence}, the classical information structure (a centralized problem) has been considered, where in \cite{cardaliaguet2019example} it has been assumed that players have access to all the history of states of all players but not controls (we note that in the team problem setup with the classical information structure through using a classical result of Blackwell \cite{Blackwell2} in the case where each decision maker knows all the history of states of all decision makers, optimal policies can be realized as one in the centralized problem where just the global state is a sufficient statistic for optimality). As we see, information structure aspects lead to subtle differences in analysis and conclusions.

Furthermore, in the context of stochastic teams with countably infinite number of decision makers, the gap between person by person optimality (Nash equilibrium in the game-theoretic context) and global team optimality is significant since a perturbation of finitely many policies fails to deviate the value of the expected cost, thus person by person optimality is a weak condition for such a setup. Hence, without establishing the uniqueness of the mean-field solution (which may hold under strong monotonicity assumptions \cite{LyonsMeanField}), the results presented in the aforementioned papers may be inconclusive regarding global optimality of the limit equilibrium. {For example, we refer the reader to \cite{bardi2019non, delarue2020selection, cardaliaguet2019example} for non-uniqueness results and to \cite{hajek2019non, bayraktar2020non, cecchin2019convergence, lacker2018convergence} for connections between limit theories and non-uniqueness of mean-field equilibria.}
For teams and social optima control problems, the analysis has primarily focused on the LQG model where the centralized performance has been shown to be achieved asymptotically by decentralized controllers (see e.g., \cite{huang2012social, arabneydi2015team}).


In this paper, we will adopt a different and novel approach. First, under symmetry of information structures and cost functions, we  show that optimal policies are of an exchangeable type for both teams with finite and countably infinite number of decision makers. {Then, in view of our topology on policies, we develop a de Finetti type representation theorem that characterizes the set of optimal policies as the extreme points of a convex set.  }

{\bf Connections with existence results on decentralized stochastic control.} We also note that compared to the results on the existence of a globally optimal policy in team problems where (finite) $N$-decision maker team problems has been considered \cite{yuksel2018general, gupta2014existence, YukselSaldiSICON17, saldi2019topology}, we study stochastic team problems with countably infinite number of decision makers. 

In our approach, we use randomized policies for our analysis and we define a topology on control policies for decentralized stochastic control. A consequence of our analysis is that, in the limit of countably infinitely many decision makers, one can characterize the set of optimal policies as the extreme points of a convex set of policies, which is, in turn, a subset of decentralized, independently randomized and identical policies. Such a result is not applicable to teams with finitely many decision makers. This geometric representation of the set of policies is related to the celebrated de Finetti's theorem. De Finetti's theorem implies that infinitely-exchangeable joint probability measures can be represented as mixtures (convex combination) of identical and independent probability measures \cite{aldous2006ecole,hewitt1955symmetric, kingman1978uses}. 


There has been related work in the quantum information/mechanics literature. Let us first note, however, that in \cite{diaconis1980finite}, it has been shown that finite number of exchangeable probability measures can be approximated by a mixture of identical and independent probability measures, and this approximation asymptotically becomes more accurate when the number of exchangeable random variables increases. The de Finetti representation type results have been extended for quantum systems where conditional probability measures have been considered \cite{brandao2017quantum, renner2007symmetry, christandl2009finite, banica2012finetti, caves2002unknown}. In fact, for permutation-symmetric conditional probability measures, approximation results have been obtained, provided that the non-signaling property holds (a conditional independence property between local actions and other measurements given local measurement) \cite{brandao2017quantum, renner2007symmetry, christandl2009finite, banica2012finetti, caves2002unknown}. We refer readers to \cite{brunner2014bell, popescu2014nonlocality}, for a review on the connection between the non-signaling conditional probability measures and the  conditional probability measures with private and common randomness. 

We note that de Finetti type results developed for conditional probability measures in quantum information literature give us a geometric interpretation we require for strategic measures (a geometric connection between non-signaling infinitely-exchangeable conditional probability measures and conditional probability measures induced by common and private randomness). However, in the team problem setup, in addition to show this geometric connection, one is required to show that the common randomness is independent of the observations. {We address this issue by introducing an appropriate topology on policies and establishing a de Finetti type representation theorem on space of policies, properly defined and metrized.}


{
{\bf Contributions.} 
In view of the above, this paper makes the following contributions.
\begin{itemize}[wide]
\item [(i)] Under symmetry of information structures and exchangeability of the cost function, we first consider teams with $N$ DMs ($N$-DM teams) and establish the optimality of $N$-exchangeable randomized policies.  
\item[(ii)] {We introduce a suitable topology on control policies which facilitates our analysis using a de Finetti type representation theorem for decentralized relaxed policies, that is, for the probability measures induced on actions and measurements under decentralized information structures. This leads to a representation theorem for decentralized relaxed policies which admit an infinite exchangeability condition.}
\item [(iii)] By extending $N$-exchangeable policies to infinitely-exchangeable ones, establishing a convergence argument for the induced costs, and using the presented de Finetti theorem for decentralized relaxed policies, we establish the structure, and also the existence of optimal decentralized policies for static and dynamic teams with countably infinite number of decision makers, which turns out to be symmetric (i.e., identical) and randomized. Compared to our previous results for static and dynamic mean-field teams in \cite[Theorem 12 or Proposition 1]{sanjari2018optimal} and \cite[Theorem 3.4]{sanjari2019optimal}: i) the cost function is not necessarily convex in actions, ii) action spaces are not necessarily convex, and iii) the mean-field coupling is considered in dynamics, which leads to a non-classical information structure (a consequence being that the problem is in general non-convex in policies).
\item[(iv)] For $N$-decision maker symmetric teams with a symmetric information structure, we show that symmetric (identical) randomized policies of mean-field teams are nearly optimal. 
\end{itemize}
}


\section{Preliminaries and Statement of Main Results}
We begin by Witsenhausen's intrinsic model for team problems, and then, we provide a description for main problems studied in this paper.
\subsection{Preliminaries}\label{sec:pre}
In this section, we introduce Witsenhausen's \textit{Intrinsic Model} for sequential teams \cite{wit75}. 

\begin{itemize}[wide]
\item There exists a collection of \textit{measurable spaces} $\{(\Omega, {\mathcal F}), \allowbreak(\mathbb{U}^i,{\mathcal U}^i), (\mathbb{Y}^i,{\mathcal Y}^i), i \in {\mathcal{N}}\}$, specifying the system's distinguishable events, and control and measurement spaces. The set $\mathcal{N}$ denotes the collection of decision makers. The set $\mathcal{N}$ can be a finite set $\{1,2,\dots, N\}$ or a countable set $\mathbb{N}$. The pair $(\Omega, {\mathcal F})$ is a
measurable space (on which an underlying probability may be defined). The pair $(\mathbb{U}^i, {\mathcal U}^i)$
denotes the standard Borel space from which the action $u^i$ of DM$^i$ is selected. The pair $(\mathbb{Y}^i,{\mathcal Y}^i)$ denotes the standard Borel observation/measurement space for each decision maker $i$ (DM$^{i}$).

\item There is a \textit{measurement constraint} to establish the connection between the observation variables and the system's distinguishable events. The $\mathbb{Y}^i$-valued observation variables are given by $y^i=h^i(\omega,{\underline u}^{[1,i-1]})$, where ${\underline u}^{[1,i-1]}:=(u^1,\dots, u^{i-1})$ and $h^i$s are measurable functions. 
\item The set of admissible control laws $\underline{\gamma}:= (\gamma^i)_{i \in \mathcal{N}}$, also called
{\textit{designs}} or {\textit{policies}}, are measurable control functions, so that $u^i = \gamma^i(y^i)$. Let $\Gamma^i$ denote the set of all admissible policies for DM$^i$ and let ${\Gamma} = \prod_{i \in \mathcal{N}} \Gamma^i$. These policies will later be allowed to be randomized and accordingly the image will be ${\mathcal P}(\mathbb{U}^i)$, where ${\mathcal P}(\cdot)$ denotes the space of probability measures.
\end{itemize}
\begin{itemize}
\item There is a {\textit{probability measure}} $\mathbb{P}$ on $(\Omega, {\mathcal F})$ describing the probability space on which the system is defined.
\end{itemize}

Under this intrinsic model, a sequential team problem is {\textit{dynamic}} if the
information available to at least one DM is affected by the action of at least one other DM. A team problem is {\it static}, if for every DM the information available is only affected by exogenous disturbances; that is no other DM can affect the information at any given DM. Information structures can also be categorized as {\it classical}, {\it quasi-classical} or {\it non-classical}. An Information Structure (IS) $\{y^i, i \in \mathcal{N} \}$ is {\it classical} if $y^i$ contains all of the information available to DM$^k$ for $k < i$. An IS is {\it quasi-classical} or {\it partially nested}, if whenever $u^k$, for some $k < i$, affects $y^i$ through the measurement function $h^i$, $y^i$ contains $y^k$ (that is $\sigma(y^k) \subset \sigma(y^i)$). An IS which is not partially nested is {\it non-classical}.\par

{In the paper, we will also allow for randomized policies, where in addition to $y^i$, each decision maker DM$^i$ has access to common and private randomization. This will be made precise later in Section \ref{sec:strategic}. }
 
\subsection{Problem statement}\label{sec:2.2}
We consider stochastic team problems with finite but large as well as team problems with countably infinite number of DMs. We address three main problems: (i) existence and structural results for static teams with countably infinite number of DMs (Section \ref{sec:section4}) (ii) existence and structural results for dynamic teams with countably infinite number of DMs (Section \ref{sec:4}) (iii) approximation results for $N$-DM static and dynamic teams (Section \ref{sec:aproxi}).

\subsubsection{Static Teams}
As we consider exchangeable team problems, we let action and observation spaces be identical through DMs $\mathbb{U}^{i}=\mathbb{U}\subseteq\mathbb{R}^{n}$ and $\mathbb{Y}^{i}=\mathbb{Y}\subseteq\mathbb{R}^{m}$ for all $i \in \mathcal{N}$, where $n$ and $m$ are positive integers. 

\begin{itemize}[wide]
\item[\bf\text{Problem} \bf{($\mathcal{P}_{N}$)}:]
Let ${\mathcal N}=\{1,\dots, N\}$. Let $\underline{\gamma}_{N} := (\gamma^1, \cdots, \gamma^N)$ and ${\Gamma}_{N} := \prod_{i=1}^{N} \Gamma^i$. Let an expected cost function of $\underline{\gamma}_{N}$ be given as
\begin{equation}\label{eq:1.1}
J_{N}(\underline{\gamma}_{N}) = \mathbb{E}^{\underline{\gamma}_{N}}[c(\omega_{0},\underline{u}_{N})]:= \mathbb{E}[c(\omega_{0},\gamma^1(y^1),\cdots,\gamma^N(y^N))],
\end{equation}
for some Borel measurable cost function $c: \Omega_{0} \times \prod_{k=1}^{N} \mathbb{U} \to \mathbb{R}_{+}$. We define $\omega_{0}$ as the $\Omega_{0}$-valued, cost function relevant, exogenous random variable as $\omega_{0}:(\Omega,\mathcal{F}, \mathbb{P}) \to (\Omega_{0},\mathcal{F}_{0})$, where $\Omega_{0}$ is a Borel space with its Borel $\sigma$-field $\mathcal{F}_{0}$. Here, we have the notation $\underline {u}_{N}:=(u^1,\dots,u^{N})$.
\end{itemize}
\begin{definition}\label{eq:gof}
For a given stochastic team problem  ($\mathcal{P}_{N}$) with a given information
structure, a policy (strategy)\/ ${\underline \gamma}^{*}_{N}:=({\gamma^1}^*,\ldots, {\gamma^N}^*)\in { \Gamma}_{N}$\@ is
 \textit{(globally) optimal} for ($\mathcal{P}_{N}$) if
\begin{equation*}
J_{N}({\underline \gamma}^*_{N})=\inf_{{{\underline \gamma}_{N}}\in {{\Gamma}_{N}}}
J_{N}({{\underline \gamma}}_{N}).
\end{equation*}  
\end{definition}

{Our focus in this paper is on a class of exchangeable team problems satisfying an exchangeability assumption on the cost function. 
\begin{assumption}\label{assump:exccost}
The cost function is exchangeable with respect to actions for all $\omega_{0}$, i.e., for any permutation $\sigma$ of $\{1,\dots,N\}$, $c(\omega_{0}, u^{1},\dots, u^{N})=c(\omega_{0}, u^{\sigma(1)},\dots, u^{\sigma(N)})$ for all $\omega_{0}$.
\end{assumption}

In particular, for our main results, we focus on team problems with the following expected cost function instead of \eqref{eq:1.1}:
\begin{equation}\label{eq:exceee}
 \mathbb{E}^{\underline{\gamma}_{N}}\bigg[\frac{1}{N}\sum_{i=1}^{N}c\bigg(\omega_{0},u^{i},\frac{1}{N}\sum_{p=1}^{{N}}u^{p}\bigg)\bigg].
\end{equation}

Clearly, $\frac{1}{N}\sum_{i=1}^{N}c(\omega_{0},u^{i},\frac{1}{N}\sum_{p=1}^{{N}}u^{p})$ satisfies Assumption \ref{assump:exccost}.} 
Now, we introduce a stochastic team problem with countably infinite number of decision makers.
\begin{itemize}[wide]
\item[\bf\text{Problem} \bf{($\mathcal{P}_{\infty}$)}:]
Consider a stochastic team with countably infinite number of decision makers, that is, $\mathcal{N}=\mathbb{N}$. Let ${\Gamma}:=\prod_{i \in \mathbb{N}} \Gamma^{i}$ and $\underline{\gamma}:=(\gamma^{1},\gamma^{2},\dots)$. Let an expected cost of $\underline{\gamma}$ be given as
\begin{equation}\label{eq:2.5.5}
J(\underline{\gamma})=\limsup\limits_{N\rightarrow \infty}  \mathbb{E}^{\underline{\gamma}}\bigg[\frac{1}{N}\sum_{i=1}^{N}c\bigg(\omega_{0},u^{i},\frac{1}{N}\sum_{p=1}^{{N}}u^{p}\bigg)\bigg],
\end{equation}
for some Borel measurable cost function $c:\Omega_{0} \times \mathbb{U}\times \mathbb{U} \rightarrow \mathbb{R}_{+}$.
\end{itemize}
\begin{definition}\label{eq:goi}
For a given stochastic team problem ($\mathcal{P}_{\infty}$) with a given information
structure, a policy\/ ${\underline \gamma}^{*}:=({\gamma^{1*}},\gamma^{2*},\ldots)\in {\Gamma}$\@ is
 \textit{optimal} for ($\mathcal{P}_{\infty}$) if
\begin{equation*}
J({\underline \gamma}^*)=\inf_{{{\underline \gamma}}\in {{\Gamma}}}
J({{\underline \gamma}}).
\end{equation*}   
\end{definition}

Later on, we allow DMs to apply randomized policies and provide a description of the problems within randomized policies; see \eqref{eq:pf} and \eqref{eq:pinf}. Our first goal here is to establish the existence of a symmetric (identical) randomized globally optimal policy for static mean-field team problems ${(\mathcal{P}_{\infty})}$.  To this end, we first establish $N$-exchangeability of randomized optimal policies for ${(\mathcal{P}_{N})}$ and symmetry for optimal randomized policies of ${(\mathcal{P}_{\infty})}$. Then in our Theorem \ref{the:2}, using symmetry, we establish an existence result for ${(\mathcal{P}_{\infty})}$. {Our theorems require the following absolute continuity condition under which we can equivalently view the observations of each DM as independent and also independent of $\omega_{0}$ via change of measure argument (due to Witsenhausen \cite{wit88}).

\begin{assumption}\label{assump:ind}
Assume that for every $N\in \mathbb{N}\cup \{\infty\}$, there exists a probability measure $Q^{i}$ on $\mathbb{Y}$ and a function $f^{i}$ for all $i\in \mathcal{N}$ such that for all Borel set $B^{i}$ in $\mathbb{Y}$ (with $B := B^{1}\times \cdots \times B^{N}$)
\begin{flalign}
&\tilde{\mu}^{N}(B \big|\omega_{0})=\prod_{i=1}^{N}\int_{B^i} f^{i}(y^{i}, \omega_{0},y^{1},\dots, y^{i-1})Q^{i}(dy^{i})\label{eq:abscon},
\end{flalign}
 where $\tilde{\mu}^{N}$ is the conditional distribution of observations $(y^{1},\dots,y^{N})$ given $\omega_{0}$.
 \end{assumption}

{\begin{remark}   
In particular, if $y^{i}$ takes values from a countable set, Assumption \ref{assump:ind} always holds e.g., with the reference measure taken as $Q^{i}(r)=\sum_{p\geq1}2^{-p}1_{\{r=m_{p}\}}$ where $\mathbb{Y}=\{m_{p}~|~ p\in \mathbb{N}\}$ (see \cite{wit88}).
\end{remark}}

The above allows us to introduce a suitable topology under which the space of randomized policies is Borel (see Section \ref{sec:strategic}). In addition, our main Theorem \ref{the:2} imposes the following assumptions on the observations and action space. 

\begin{assumption}\label{assump:oc}
\hfill 
\begin{itemize}
\item [(i)] Observations $(y^{i})_{i\in \mathcal{N}}$ are i.i.d. conditioned on $\omega_{0}$; 
\item [(ii)] $\mathbb{U}$ is compact.
\end{itemize}
\end{assumption}

We note that under Assumption \ref{assump:ind} and Assumption \ref{assump:oc}(i), there exists an identical reference probability measure $Q$ and function $f$ such that the absolute continuity condition \eqref{eq:abscon} holds; that is, for any Borel set $B^{i}$ in $\mathbb{Y}$ (with $B := B^{1}\times \cdots \times B^{N}$)
\begin{flalign*}
\tilde{\mu}^{N}(B \big|\omega_{0})&=\prod_{i=1}^{N}\hat{\mu}(B^{i}|\omega_{0})\\
&=\prod_{i=1}^{N}\int_{B^{i}}f(y^{i},\omega_{0})Q(dy^{i}),
\end{flalign*} 
where $\hat{\mu}$ is the conditional distribution of each observation $y^{i}$ given $\omega_{0}$. We note that the function $f$ and the measure $Q$ are identical through DMs since observations are identically distributed conditioned on $\omega_{0}$. Furthermore, our main Theorem \ref{the:2} imposes the following continuity assumption on the cost function.

\begin{assumption}\label{assump:cont}
The cost function in \eqref{eq:exceee}, $c:\Omega_{0} \times \mathbb{U}\times \mathbb{U} \rightarrow \mathbb{R}_{+}$, is continuous in its second and third arguments for all $\omega_{0}$.
\end{assumption}

For our results in Section \ref{sec:section4}, we impose Assumption \ref{assump:exccost} and Assumption \ref{assump:ind}, but we only impose Assumption \ref{assump:oc} and Assumption \ref{assump:cont} when they are needed.}

\subsubsection{Dynamic Teams}\label{sec:dn}

Our second goal here is to establish the existence of a symmetric (identical) randomized globally optimal policy for mean-field dynamic team problems where DMs are weakly coupled through the average of states and actions in dynamics and/or the cost function. {Again, we consider exchangeable teams, and hence, we let action, observation, and state spaces, respectively, be identical through DMs $i\in \mathcal{N}$, and for simplicity, also through time $t=0,\dots,T-1$, $\mathbb{U}^{i}_{t}=\mathbb{U}\subseteq\mathbb{R}^{n}$, $\mathbb{Y}^{i}_{t}=\mathbb{Y}\subseteq\mathbb{R}^{n^{\prime}}$, $\mathbb{X}^{i}_{t}=\mathbb{X}\subseteq\mathbb{R}^{n^{\prime\prime}}$ for all $i \in \mathcal{N}$ and $t=0,\dots,T-1$, where $n$, $n^{\prime}$ and $n^{\prime\prime}$  are positive integers. Define state dynamics and observation dynamics of DMs as follows:
\begin{flalign}
x_{t+1}^{i}&=f_{t}\bigg(x_{t}^{i},u_{t}^{i},\frac{1}{N}\sum_{p=1}^{N}x_{t}^{p}, \frac{1}{N}\sum_{p=1}^{N}u_{t}^{p}, w_{t}^{i}\bigg),\label{eq:mfdynamics2}\\
y_{t}^{i}&=h_{t}\bigg(x_{0:t}^{i},u_{0:t-1}^{i},v_{0:t}^{i}\bigg)\label{eq:mfobs2},
\end{flalign}
where functions $f_{t}$ and $h_{t}$ are measurable functions and $v_{t}^{i}$ and $w_{t}^{i}$ are random vectors representing uncertainties in state dynamics and observations. We denote $x_{0:t}^{i}:=(x_{0}^{i},\dots, x_{t}^{i})$, $u_{0:t-1}^{i}:=(u_{0}^{i},\dots, u_{t-1}^{i})$, and $v_{0:t}^{i}:=(v_{0}^{i},\dots, v_{t}^{i})$. Let the admissible policies $(\gamma^{i}_{0:T-1})_{i\in \mathcal{N}}$ (with $\gamma^{i}_{0:T-1}:=(\gamma^{i}_{0},\dots, \gamma^{i}_{T-1})$) be measurable control functions so that $u^i_{t} = \gamma^i_{t}(y^i_{t})$ for all $i\in \mathcal{N}$ and $t=0,\dots, T-1$. }
 \begin{itemize}[wide]
  \item[ \textbf{Problem} \bf{($\mathcal{P}_{T}^{N}$)}:] Consider $N$-DM mean-field dynamic teams with the expected cost function of $\underline{\gamma}^{1:N}$ as 
  \begin{flalign}
&\scalemath{0.95}{J_{T}^{N}(\underline{\gamma}^{1:N})=\mathbb{E}^{\underline{\gamma}^{1:N}}\bigg[\frac{1}{N}\sum_{t=0}^{T-1}\sum_{i=1}^{N}c\bigg(\omega_{0}, x_{t}^{i},u_{t}^{i},\frac{1}{N}\sum_{p=1}^{N}u_{t}^{p},\frac{1}{N}\sum_{p=1}^{N}x_{t}^{p}\bigg)\bigg]}\label{eq:mfcost},
\end{flalign}
where $\underline{\gamma}^{1:N}:=(\gamma^{1}_{0:T-1}, \dots, \gamma^{N}_{0:T-1})$ and $\gamma^{i}_{0:T-1}:=(\gamma^{i}_{0},\dots, \gamma^{i}_{T-1})$. Again, $\omega_{0}:(\Omega,\mathcal{F}, \mathbb{P}) \to (\Omega_{0},\mathcal{F}_{0})$ is a cost-related random variable, where $\Omega_{0}$ is a Borel space with its Borel $\sigma$-field $\mathcal{F}_{0}$.
\end{itemize}
\begin{itemize}[wide]
\item[\textbf{Problem} \bf{($\mathcal{P}_{T}^{\infty}$)}:] Consider mean-field dynamic teams with the expected cost function of $\underline{\gamma}$ as
 \begin{flalign}
&\scalemath{0.95}{J_{T}^{\infty}(\underline{\gamma})=\limsup\limits_{N\rightarrow \infty}J_{T}^{N}(\underline{\gamma}^{1:N})}\label{eq:dmfcost},
\end{flalign}
where $\underline{\gamma}:=(\gamma^{1}_{0:T-1},\gamma^{2}_{0:T-1},\dots)$ and $\underline{\gamma}^{1:N}:=({\gamma}^{1}_{0:T-1}, \dots, {\gamma}^{N}_{0:T-1})$.
\end{itemize}

Analogous to Definition \ref{eq:gof} and Definition \ref{eq:goi}, we can define globally optimal policies for ${(\mathcal{P}^{N}_{T})}$ and  ${(\mathcal{P}^{\infty}_{T})}$. Again, we allow DMs to apply randomized policies and provide a description of the problems within randomized policies; see \eqref{eq:finiterandom} and \eqref{eq:infiniterandom}. In Section \ref{sec:4}, we establish the existence of a symmetric (identical through DMs) randomized globally optimal policy for ${(\mathcal{P}^{\infty}_{T})}$. Similar to the static case, we first establish $N$-exchangeablity of randomized optimal policies for ${(\mathcal{P}^{N}_{T})}$ and symmetry for optimal randomized policies of ${(\mathcal{P}^{\infty}_{T})}$. Then using symmetry, we establish an existence result for ${(\mathcal{P}^{\infty}_{T})}$. 

{Our solution technique for dynamic problems is similar to the static one, which requires more technical arguments and additional assumptions. Our theorems for dynamic case impose an absolute continuity condition (see Assumption \ref{assump:ind1}) to allow us to introduce a suitable topology on control policies and to facilitate our analysis (our main Theorem \ref{the:exmftdy2}  requires an additional technical Assumption \ref{assump:ind2}). Furthermore, our main Theorem \ref{the:exmftdy2} imposes the following:

\begin{assumption}\label{assump:c}
\hfill 
\begin{itemize}
\item[(i)] For $t=0,\dots,T-1$, functions $f_{t}$ and $h_{t}$ in \eqref{eq:mfdynamics2} and \eqref{eq:mfobs2} are continuous in the states and actions and $f_{t}$s are bounded;
\item [(ii)]  The cost function in \eqref{eq:mfcost}, $c:\Omega_{0}\times \mathbb{X}\times \mathbb{U} \times \mathbb{U} \times \mathbb{X} \to \mathbb{R}_{+}$, is continuous in the second, third, fourth, and fifth arguments.
\end{itemize}
 \end{assumption}
\begin{assumption}\label{assump:2}
\hfill 
\begin{itemize}
\item[(i)] $(x_{0}^{i})_{i\in \mathcal{N}}$ are i.i.d. random vectors conditioned on $\omega_{0}$;  
\item[(ii)] For $t=0,\dots,T-1$, $(w^{i}_{t})_{i\in \mathcal{N}}$ are i.i.d. random vectors, and for $i\in \mathcal{N}$, $(w^{i}_{t})_{t=0}^{T-1}$ are mutually independent, and independent of $\omega_{0}$ and $(x_{0}^{i})_{i\in \mathcal{N}}$. For $t=0,\dots,T-1$, $(v^{i}_{t})_{i\in \mathcal{N}}$ are i.i.d. random vectors, and for $i\in \mathcal{N}$, $(v^{i}_{t})_{t=0}^{T-1}$ are mutually independent, and independent of $\omega_{0}$, $(x_{0}^{i})_{i\in \mathcal{N}}$, and $w^{i}_{t}$s for $i\in \mathcal{N}$ and $t=0,\dots,T-1$.
\item [(iii)] $\mathbb{U}$ is compact.
\end{itemize}
\end{assumption}

In view of Assumption \ref{assump:2}(i), we note that $\omega_{0}$ also introduces a correlation between initial states. For our results in Section \ref{sec:4}, we impose Assumption \ref{assump:ind1}, but we impose Assumption \ref{assump:ind2}, Assumption \ref{assump:c}, and Assumption \ref{assump:2} only when they are needed.}

\subsubsection{Approximations}

Finally, we address the following problem in Section \ref{sec:aproxi}. If $P_{\pi}^{*}$ is a (randomized) symmetric optimal policy for ($\mathcal{P}_{\infty}$) (($\mathcal{P}_{T}^{\infty}$)) then there exist $\epsilon_{N}\geq 0$, with $\epsilon_{N} \to 0$ as $N \to\infty$, such that $P_{\pi}^{*}|_{N}$ is $\epsilon_{N}$-optimal for ($\mathcal{P}_{N}$) (($\mathcal{P}_{T}^{N}$)) where $P_{\pi}^{*}|_{N}$ is the restriction of $P_{\pi}^{*}$ to the first $N$ decision makers. We use our symmetry results and analysis for ($\mathcal{P}_{\infty}$) (($\mathcal{P}_{T}^{\infty}$)).

\subsection{Discussion of main results}
 In mean-field team problems, one may be interested in the existence and structure of globally optimal policies. In particular, one can ask if there is a globally optimal policy and whether this optimal policy is symmetric for these type of problems (by a symmetric policy we mean that a policy is identical through DMs). One may be also interested in the connection between optimal policies for mean-field teams and approximation of optimal policies for the pre-limit $N$-DM teams when $N$ is large. The purpose of this paper is to address these questions for mean-field team problems where the problem can be non-convex. The non-convexity of the problem can arise as a result of non-convexity of the action space and/or non-convexity of the cost function in actions. Also, even if the action space is convex and the cost function is convex in actions, the information structure of the problem may lead to non-convexity of the problem in policies (see for example \cite[Section 3.3]{YukselSaldiSICON17}).  A celebrated example is the counterexample of Witsenhausen \cite{wit68}.


One of the main difficulties in studying non-convex mean-field team problems is to show that globally optimal policies for mean-field team problems are symmetric (identical for each DM). This difficulty stems from the observation that, in general, globally optimal policies are not symmetric for non-convex pre-limit $N$-DM team problems (which can be seen in Example \ref{ex:1}). This is in contrast to the convex mean-field teams where symmetry can be established for both pre-limit $N$-DM and mean-field team problems  \cite{sanjari2018optimal, sanjari2019optimal}. In our approach: 
\begin{itemize}[wide]
\item [(i)] {We introduce a topology on control polices which is used to establish a de Finetti representation result for probability measures on policies identified as randomized policies.} In Theorem \ref{the:defin}, we show that any infinitely-exchangeable randomized  policies can be represented by elements of the set of randomized policies with common and private independent randomness where conditioned on common randomness, randomization of the policies are independent and identical through DMs.   
 \item [(ii)] In Section \ref{sec:section4} for static and Section \ref{sec:4} for dynamic $N$-DM stochastic teams (see Lemma \ref{lem:exc} and Lemma \ref{lem:3.1dy}), we show that by exchangeability of the cost function and considering symmetric information structures (under a causality condition for the dynamic case), one can establish $N$-exchangeability of randomized optimal policies.
 \item [(iii)] In Section \ref{sec:section4} for static and Section \ref{sec:4} for dynamic mean-field teams (see Lemma \ref{lem:findef} and Lemma \ref{lem:findefdyn2}) under regularity conditions on the cost function and dynamics, by constructing infinitely-exchangeable randomized policies by relabeling $N$-exchangeable randomized optimal policies, as $N$ goes to infinity, we show the asymptotic optimality of infinitely-exchangeable randomized optimal policies. Hence, this, following from our de Finetti type theorem (see Theorem \ref{the:defin}), establishes asymptotic global optimality of symmetric and conditionally independent policies. 
 \item [(iv)] Using extreme point and lower semi-continuity arguments, we establish the existence of a symmetric optimal policy (which is privately randomized) for static and dynamic mean-field teams (see Theorem \ref{the:2} and Theorem \ref{the:exmftdy2}).
 \item [(v)] In Section \ref{sec:aproxi}, using our analysis for mean-field problems, as $N$ goes to infinity, we show that symmetric optimal policies of mean-field teams are asymptotically optimal for $N$-DM weakly coupled teams, hence, it establishes approximation results for this class of problems.
 


 In the following, we first study static teams, then we study dynamic teams where the analysis is similar to the static case but is somewhat more technical.
\end{itemize}
\section{Topology on Control Policies and a de Finetti Representation Result}
\subsection{Topology on control policies.}\label{sec:strategic}

In this section, we introduce a topology using which, we can introduce Borel probability measures on policies. We first consider $N$-DM static team problems. Following from \cite{yuksel2018general, wit88}, Assumption \ref{assump:ind} allows us to reduce the problem as a static team problem where now the observation of each DM is independent of observations of other DMs and also independent of $\omega_{0}$ {(since under  the measure transformation \eqref{eq:abscon}, a probability measure on the observation of each DM is $Q^{i}$, which is independent of observations of other DMs and $\omega_{0}$). Hence, under Assumption \ref{assump:ind}, we can focus on each DM$^{i}$ separately. Let us define 
 \begin{flalign}
 &\Theta^{i}:=\bigg\{P \in \mathcal{P}(\mathbb{U}\times \mathbb{Y}) \bigg| P(B)=\int_{B}1_{\{g^{i}(y^{i})\in du^{i}\}}Q^{i}(dy^{i}),~g:\mathbb{Y}\to \mathbb{U}, ~ B\in \mathcal{B}(\mathbb{U}\times \mathbb{Y})\bigg\}\label{eq:topogamd},
 \end{flalign}
where $\mathcal{P}(\cdot)$ denotes the space of probability measures, and $1_{\{\cdot ~\in~A\}}$ denotes the indicator function of the set $A$. The above set is the set of extreme points of the set of probability measures 
on $(\mathbb{U}\times \mathbb{Y})$ with fixed marginals $Q^{i}$ on $\mathbb{Y}$, that is,
\begin{flalign}
\mathcal{R}^{i}:=\bigg\{P \in \mathcal{P}(\mathbb{U}\times \mathbb{Y}) \bigg| P(B)=\int_{B}\Pi^{i}(du^{i}|y^{i})Q^{i}(dy^{i}), ~ B\in \mathcal{B}(\mathbb{U}\times \mathbb{Y})\bigg\}\label{eq:topogam},
\end{flalign}
where $\Pi^{i}$ is a stochastic kernel from $\mathbb{Y}$ to $\mathbb{U}$. Hence, it inherits Borel measurability and topological properties of that Borel measurable set \cite{BorkarRealization}. We note that this set corresponds to Young measures \cite{young1937generalized} and this representation result is due to Borkar \cite{BorkarRealization}. Now, we identify the set of relaxed policies $\Gamma^{i}$ by $\mathcal{R}^{i}$ and we define convergence on policies as $\gamma^{i}_{n} \to \gamma^{i}$ if and only if $\gamma^{i}_{n}(du^{i}|y^{i})Q^{i}(dy^{i}) \to \gamma^{i}(du^{i}|y^{i})Q^{i}(dy^{i})$ (in the weak convergence topology) as $n\to \infty$.  }


 In view of the above standard Borel space formulation for $\Gamma^{i}$ for each $i \in \mathcal{N}$, we can define the set of Borel probability measures on admissible policies $\Gamma_{N}$ (which is referred to as a set of randomized policies) as $L^{N}:=\mathcal{P}(\Gamma_{N})$, where Borel $\sigma$-field $\mathcal{B}(\Gamma^{i})$ is induced by the topology defined above. Define the set of randomized policies induced by a common and individual randomness as: 
\begin{flalign*}
L_{\text{CO}}^{N}:=\bigg\{&P_{\pi} \in L^{N}\bigg{|}\text{for all}~A_{i} \in \mathcal{B}(\Gamma^{i}): \nonumber\\
&P_{\pi}(\gamma^{1} \in A_{1},\dots,\gamma^{N}\in A_{N})=\int_{z\in [0,1]}\prod_{i=1}^{N}P_{\pi}^{i}(\gamma^{i}\in A_{i}|z)\eta(dz), ~~~~\eta \in \mathcal{P}([0, 1])\bigg\},
\end{flalign*} 
{where $\eta$ is the distribution of common, but independent (from intrinsic exogenous system variables), randomness}, and for every fixed $z$, ${P}_{\pi}^{i}\in \mathcal{P}(\Gamma^{i})$ indicates an identical independent randomized policy of each DM$^{i}$ ($i=1,\dots, N$). Note that conditioned on a $[0, 1]$-valued random variable $Z$, policies are independent. {It can be shown that $L_{\text{CO}}^{N}$ and $L^{N}$ are identical (see Theorem \ref{the:lcoln} in the Appendix), and hence,  the set of randomized policies $L^{N}$ corresponds to randomized policies induced by an individual and a common randomness. Since individual and a common randomness do not improve the optimal expected cost, the relaxation of the problem to sets of randomized policies  $L^{N}$ is a legitimate relaxation for the team problems with $N$-DMs.}

{Before, we introduce the set of exchangeable randomized policies, we recall the definition of {\it exchangeability} for random variables.
\begin{definition}
Random vectors $x^{1},x^{2},\dots,x^{N}$ defined on a common probability space are $N$-\it{exchangeable} if for any permutation $\sigma$ of the set $\{1,\dots,N\}$, 
\begin{flalign*}
&\mathcal{L}\bigg(x^{\sigma(1)},x^{\sigma(2)},\dots,x^{\sigma(N)}\bigg)=\mathcal{L}\bigg(x^{1},x^{2} ,\dots, x^{N}\bigg),
\end{flalign*}
where $\mathcal{L}$ denotes the joint distribution of random vectors. Random vectors $(x^{1},x^{2},\dots)$ is {\it infinitely-exchangeable} if finite distributions of $(x^{1},x^{2},\dots)$ and $(x^{\sigma(1)},x^{\sigma(2)},\dots)$ are identical for any finite permutation (affecting only finitely many elements) of $\mathbb{N}$. 
\end{definition}}
Now, we define the set of exchangeable randomized policies as:
\begin{flalign}
L_{\text{EX}}^{N}:=\bigg\{&P_{\pi} \in L^{N}\bigg{|}\text{for all}~A_{i} \in \mathcal{B}(\Gamma^{i})~\text{and for all}~\sigma \in S_{N}:\nonumber\\
&P_{\pi}(\gamma^{1} \in A_{1},\dots,\gamma^{N}\in A_{N})=P_{\pi}(\gamma^{\sigma(1)} \in A_{1},\dots,\gamma^{\sigma(N)}\in A_{N})\bigg\}\label{eq:LEXN},
\end{flalign}
where $S_{N}$ is the set of permutations of $\{1,\dots,N\}$. We note that $L_{\text{EX}}^{N}$ is a convex subset of $L^{N}$.  {We also define the set $L_{\text{CO,SYM}}^{N}$ as the set of identical randomized  policies induced by a common and individual randomness:
\begin{flalign*}
\scalemath{0.98}{L_{\text{CO,SYM}}^{N}:=\bigg\{}&\scalemath{0.98}{P_{\pi} \in L^{N}\bigg{|}\text{for all}~A_{i} \in \mathcal{B}(\Gamma^{i}):} \nonumber\\
&\scalemath{0.98}{P_{\pi}(\gamma^{1} \in A_{1},\dots,\gamma^{N}\in A_{N})=\int_{z\in [0,1]}\prod_{i=1}^{N}\tilde{P}_{\pi}(\gamma^{i} \in A_{i}|z)\eta(dz), ~~~~\eta \in \mathcal{P}([0, 1])\bigg\}},
\end{flalign*}
where for all $i\in \mathcal{N}$, and fixed $z$, $\tilde{P}_{\pi}\in \mathcal{P}(\Gamma^{i})$ indicates an identical independent randomized policy of each DM$^{i}$ ($i=1,\dots, N$). Also, define the set of randomized policies with only private independent randomness as:
\begin{flalign}
L_{\text{PR}}^{N}:=\bigg\{&P_{\pi} \in L^{N}\bigg{|}\text{for all}~A_{i} \in \mathcal{B}(\Gamma^{i}): \nonumber\\
&P_{\pi}(\gamma^{1} \in A_{1},\dots, \gamma^{N} \in A_{N})=\prod_{i=1}^{N}P_{\pi}^{i}(\gamma^{i}\in A_{i}),~\text{for}~P_{\pi}^{i}\in \mathcal{P}(\Gamma^{i})\bigg\}\nonumber.
\end{flalign} 
Finally, define the set of randomized policies with identical and independent randomness:
\begin{flalign}
L_{\text{PR,SYM}}^{N}:=\bigg\{&P_{\pi} \in L^{N}\bigg{|}\text{for all}~A_{i} \in \mathcal{B}(\Gamma^{i}):\nonumber\\
& P_{\pi}(\gamma^{1} \in A_{1},\dots, \gamma^{N} \in A_{N})=\prod_{i=1}^{N}\tilde{P}_{\pi}(\gamma^{i}\in A_{i}),~\text{for}~\tilde{P}_{\pi} \in \mathcal{P}(\Gamma^{i})\bigg\}\nonumber.
\end{flalign}

For a team with a countably infinite number of decision makers, we define sets of randomized policies $L, L_{\text{EX}}, L_{\text{CO}}, L_{\text{CO,SYM}}, L_{\text{PR}}, L_{\text{PR,SYM}}$ similarly using Ionescu Tulcea
extension theorem through the sequential formulation reviewed in Section \ref{sec:pre}, by iteratively adding new coordinates for our probability measure (see e.g., \cite{InfiniteDimensionalAnalysis, HernandezLermaMCP}). We define the set of randomized policies $L$ on the infinite product Borel spaces $\Gamma=\prod_{i\in \mathbb{N}}\Gamma^{i}$ as
$L:=\mathcal{P}(\Gamma)$.
Now, we define the set of infinitely-exchangeable randomized policies as:
\begin{flalign*}
L_{\text{EX}}:=\bigg\{&P_{\pi} \in L\bigg{|}\text{for all}~A_{i} \in \mathcal{B}(\Gamma^{i})~\text{and for all $N\in\mathbb{N}$, and for all $\sigma \in S_{N}$:} \nonumber\\
&P_{\pi}(\gamma^{1} \in A_{1},\dots,\gamma^{N}\in A_{N})=P_{\pi}(\gamma^{\sigma(1)} \in A_{1},\dots,\gamma^{\sigma(N)}\in A_{N})\bigg\},
\end{flalign*}
and we define
\begin{flalign*}
L_{\text{CO}}:=\bigg\{&P_{\pi} \in L\bigg{|}\text{for all}~A_{i} \in \mathcal{B}(\Gamma^{i}): \nonumber\\
&P_{\pi}(\gamma^{1} \in A_{1},\gamma^{2} \in A_{2},\dots)=\int_{z\in [0,1]}\prod_{i\in \mathbb{N}}P_{\pi}^{i}(\gamma^{i}\in A_{i}|z)\eta(dz), ~~~~\eta \in \mathcal{P}([0, 1])\bigg\}.
\end{flalign*} 
 Note that $L_{\text{CO}}$ is a convex subset of $L$ and its extreme points are in the set of randomized policies with private independent randomness:
\begin{flalign}
L_{\text{PR}}:=\bigg\{&P_{\pi} \in L\bigg{|}\text{for all}~A_{i} \in \mathcal{B}(\Gamma^{i}): \nonumber\\
&P_{\pi}(\gamma^{1} \in A_{1},\gamma^{2} \in A_{2},\dots)=\prod_{i\in \mathbb{N}}P_{\pi}^{i}(\gamma^{i}\in A_{i}),~\text{for}~P_{\pi}^{i}\in \mathcal{P}(\Gamma^{i})\bigg\}\nonumber.
\end{flalign} 
{Also, we define
\begin{flalign*}
\scalemath{0.98}{L_{\text{CO,SYM}}:=\bigg\{}&\scalemath{0.98}{P_{\pi} \in L\bigg{|}\text{for all}~A_{i} \in \mathcal{B}(\Gamma^{i}):} \nonumber\\
&\scalemath{0.98}{P_{\pi}(\gamma^{1} \in A_{1},\gamma^{2}\in A_{2}, \dots)=\int_{z\in [0,1]}\prod_{i\in \mathbb{N}}\tilde{P}_{\pi}(\gamma^{i} \in A_{i}|z)\eta(dz), ~~~~\eta \in \mathcal{P}([0, 1])\bigg\}},
\end{flalign*}
and 
\begin{flalign*}
L_{\text{PR,SYM}}:=\bigg\{&P_{\pi} \in L\bigg{|}\text{for all}~A_{i} \in \mathcal{B}(\Gamma^{i}):\\
& P_{\pi}(\gamma^{1} \in A_{1},\gamma^{2}\in A_{2}, \dots)=\prod_{i\in \mathbb{N}}\tilde{P}_{\pi}(\gamma^{i} \in A_{i}),~\text{for}~\tilde{P}_{\pi} \in \mathcal{P}(\Gamma^{i})\bigg\}.
\end{flalign*}}

\subsection{A de Finetti theorem for admissible team policies}

{{In view of the introduced topology and sets of Borel probability measures on policies (sets of randomized polices), we now establish a connection between $L_{\text{EX}}$ and $L_{\text{CO,SYM}}$ using the classical de Finetti's theorem; that is, infinitely-exchangeable randomized policies are a mixture of i.i.d. randomized policies. 

\begin{theorem}\label{the:defin}
Any infinitely-exchangeable randomized policy $P_{\pi}\in L_{\text{EX}}$ is in the set of randomized policies $L_{\text{CO,SYM}}$ ($P_{\pi}\in L_{\text{CO,SYM}}$), i.e., for any $P_{\pi}\in L_{\text{EX}}$, there exists a $[0, 1]$-valued random variable $Z$ such that for any $A_{i} \in \mathcal{B}(\Gamma^{i})$
\begin{flalign}
&P_{\pi}(\gamma^{1} \in A_{1},\gamma^{2} \in A_{2},\dots)=\int_{z\in [0, 1]}\prod_{i\in \mathbb{N}}\tilde{P}_{\pi}(\gamma^{i}\in A_{i}|z)\eta(dz), ~~~~\eta \in \mathcal{P}([0, 1]),\label{eq:definiti}
\end{flalign} 
where for every fixed $z$, $\tilde{P}_{\pi}\in \mathcal{P}(\Gamma^{i})$.
\end{theorem}

\begin{proof}
In view of the introduced weak convergence topology on $\Gamma^{i}$ (using Borel measurable sets \eqref{eq:topogam} and \eqref{eq:topogamd}), we have $\Gamma^{i}$ is a closed subset of the Borel space $\mathcal{P}(\mathbb{U}\times \mathbb{Y})$, and hence, $\Gamma^{i}$ is Borel for $i\in \mathbb{N}$. The proof follows from \cite[Theorem 1.1]{kallenberg2006probabilistic} since $\Gamma=\prod_{i=1}^{\infty}\Gamma^{i}$ is Borel. We note that the de Finetti representation in \cite[Theorem 1.1]{kallenberg2006probabilistic} is of the form $P_{\pi}(\gamma^{1} \in A_{1},\gamma^{2} \in A_{2},\dots)=\int_{\mathcal{P}(\Gamma^{i})}\prod_{i=1}^{\infty}m(A^{i})\hat{\eta}(dm)$ for $\hat{\eta}\in \mathcal{P}(\mathcal{P}(\Gamma^{i}))$ which can be written as in \eqref{eq:definiti}. That is because, $\mathcal{P}(\Gamma^{i})$ is an (uncountable) Borel space \cite[Corollary 7.25.1]{BertsekasShreve}, and hence, by Borel-isomorphism Theorem (see for example, \cite[Proposition 7.16]{BertsekasShreve}), it is Borel isomorphic to $[0,1]$. 
\end{proof}
}

\section{Existence and Structure of Optimal Policies for Symmetric Static Team Problems with Infinitely Many Decision Makers}\label{sec:section4}

In this section, we consider static stochastic team problems where we impose Assumption \ref{assump:exccost} and Assumption \ref{assump:ind}. We note that all the proofs regarding this section are presented in  Appendix \ref{appA}. {We again note that for our results in this section, we impose Assumption \ref{assump:exccost} and Assumption \ref{assump:ind}. Based on the definitions of randomized policies, we redefine the expected cost in $(\mathcal{P}_{N})$ of a randomized policy $P_{\pi}\in L^{N}$ as:
\begin{flalign}
\scalemath{0.95}{J_{N}^{\pi}(\underline{\gamma}_{N})}&\:\scalemath{0.95}{:=\int P_{\pi}(d\underline{\gamma})\mu^{N}(d\omega_{0},d\underline{y})c^{N}(\underline{\gamma},\underline{y}, \omega_{0})}\nonumber\\
&\scalemath{0.95}{:=\int \bigg( \int c(\omega_{0},u^{1},\dots,u^{N})\prod_{k=1}^{N}\gamma^{k}(du^k|y^k)\bigg) P_{\pi}(d\gamma^{1},\dots,d\gamma^{N})\mu^{N}(d\omega_{0}, dy^{1},\dots,dy^{N})}\label{eq:finitecost},
\end{flalign}
where $c^{N}(\underline{\gamma},\underline{y}, \omega_{0}):=\int c(\omega_{0},u^{1},\dots,u^{N})\prod_{k=1}^{N}\gamma^{k}(du^k|y^k)$, and $\mu^{N}$ is the joint probability measure on measurements $(y^{1},\dots,y^{N})$ and $\omega_{0}$. In the following, we characterize team problems in which the search for a randomized optimal policy can be restricted to policies in $L_{\text{EX}}^{N}$ without losing global optimality. 

\begin{assumption}\label{assump:exobs}
Let observations of DMs, $(y^{1},\cdots,y^{N})$, be exchangeable conditioned on $\omega_{0}$.
\end{assumption}

Note that Assumption \ref{assump:exobs} is weaker than Assumption \ref{assump:oc}(i).

\begin{lemma}\label{lem:exc}
For a fixed $N$, consider an $N$-DM static team. Assume $\bar{L}^{N}$ is an arbitrary convex subset of $L^{N}$. If Assumption \ref{assump:exobs} holds, then 
\begin{equation}\label{eq:3.5}
\scalemath{0.97}{ \inf\limits_{P_{\pi} \in \bar{L}^{N}}\int P_{\pi}(d\underline{\gamma})\mu^{N}(d\omega_{0},d\underline{y})c^{N}(\underline{\gamma},\underline{y}, \omega_{0})=\inf\limits_{P_{\pi} \in \bar{L}^{N} \cap L_{\text{EX}}^{N}}\int P_{\pi}(d\underline{\gamma})\mu^{N}(d\omega_{0},d\underline{y})c^{N}(\underline{\gamma},\underline{y}, \omega_{0})}.
\end{equation}
\end{lemma}
}

In the following, we present an existence result on globally optimal policies for static mean-field teams with infinitely many decision makers. First, we re-state the infinite decision maker mean-field team problem and its pre-limit.
\begin{itemize} [wide]
\item[\bf\text{Problem} \bf{($\mathcal{P}_{N}$)}:]
Consider an $N$-DM static team with the expected cost of a randomized policy $P_{\pi}^{N}\in L^{N}$ as:
\begin{flalign}
\scalemath{0.99}{\int P_{\pi}^{N}(d\underline{\gamma})\mu^{N}(d\omega_{0},d\underline{y})c^{N}(\underline{\gamma},\underline{y}, \omega_{0}):=\int}&\scalemath{0.99}{\bigg(\int\frac{1}{N}\sum_{i=1}^{N}c\big(\omega_{0}, u^{i},\frac{1}{N}\sum_{p=1}^{N}u^{p}\big)\prod_{k=1}^{N}\gamma^{k}(du^k|y^k)\bigg)}\nonumber\\
&\scalemath{0.99}{\times P_{\pi}^{N}(d\gamma^{1},\dots,d\gamma^{N})\mu^{N}(d\omega_{0}, dy^{1},\dots,dy^{N})}\label{eq:pf}.
\end{flalign}
\end{itemize}

The above problem is considered as a pre-limit problem for our infinite-decision maker team problem.  This problem is a special case of $(\mathcal{P}_{N})$ defined in the previous section since we have a special structure for the cost function $c^{N}$ which satisfies Assumption \ref{assump:exccost}. 
\begin{itemize} [wide]
\item[\bf\text{Problem} \bf{($\mathcal{P}_{\infty}$)}:]
Consider infinite-DM static team with the following expected cost of a randomized policy $P_{\pi}\in L$ as
\begin{flalign}
&\scalemath{0.97}{\limsup\limits_{N \to \infty}\int P_{\pi,N}(d\underline{\gamma})\mu^{N}(d\omega_{0},d\underline{y})c^{N}(\underline{\gamma},\underline{y}, \omega_{0})}\label{eq:pinf},
\end{flalign}
where $P_{\pi,N}$ is the marginal of the $P_{\pi}\in L$ to the first $N$ components and $\mu^{N}$ is the marginal of the  fixed probability measure on $(\omega_{0}, y^{1},y^{2},\dots)$ to the first $N+1$ components.
\end{itemize}


In the following, we present a key result required for our main theorem. Under mild conditions, we show that the optimal expected cost function induced by $L^{N}_{\text{EX}}$ and $L_{\text{EX}}$ are equal as $N$ goes to infinity. Hence, by Lemma \ref{lem:exc}, under symmetry, this allows us to show that  without loss of global optimality, optimal policies of static mean-field teams with countably infinite number of DMs can be considered to be an infinitely-exchangeable type. 

\begin{lemma}\label{lem:findef}
Suppose that Assumption \ref{assump:oc} and Assumption \ref{assump:cont}  hold. Assume further that the cost function is bounded. Then 
\begin{flalign}
&\limsup\limits_{N \to \infty}\inf\limits_{P_{\pi}^{N} \in L^{N}_{\text{EX}}}\int P_{\pi}^{N}(d\underline{\gamma})\mu^{N}(d\omega_{0},d\underline{y}) c^{N}(\underline{\gamma},\underline{y}, \omega_{0})\nonumber\\
&=\limsup\limits_{N \to \infty}\inf\limits_{P_{\pi} \in L_{\text{EX}}}\int P_{\pi, N}(d\underline{\gamma})\mu^{N}(d\omega_{0},d\underline{y})c^{N}(\underline{\gamma},\underline{y}, \omega_{0})\label{eq:ert},
\end{flalign}
where $P_{\pi,N}$ is the marginal of the $P_{\pi}\in L_{\text{EX}}$ to the first $N$ components. 
\end{lemma}

In the following, we establish an existence of a randomized optimal policy for ($\mathcal{P}_{\infty}$). 

\begin{theorem}\label{the:2}
Consider a static team problem ($\mathcal{P}_{\infty}$) where Assumption  \ref{assump:oc} and Assumption \ref{assump:cont} hold.  
Then, there exists a randomized optimal policy $P^{*}_{\pi}$ for {($\mathcal{P}_{\infty}$)} which is in $L_{\text{PR,SYM}}$:
\begin{flalign*}
\scalemath{0.95}{\inf\limits_{P_{\pi} \in L_{\text{PR,SYM}}}}&\scalemath{0.95}{ \limsup\limits_{N \to \infty} \int P_{\pi, N}(d\underline{\gamma})\mu
^{N}(d\omega_{0},d\underline{y})c^{N}(\underline{\gamma},\underline{y}, \omega_{0})}\nonumber\\
&\scalemath{0.95}{:=\limsup\limits_{N \to \infty} \int P_{\pi, N}^{*}(d\underline{\gamma})\mu
^{N}(d\omega_{0},d\underline{y})c^{N}(\underline{\gamma},\underline{y}, \omega_{0})}\\
&\scalemath{0.95}{=\inf\limits_{P_{\pi} \in L_{\text{PR}}} \limsup\limits_{N \to \infty} \int P_{\pi, N}(d\underline{\gamma})\mu
^{N}(d\omega_{0},d\underline{y})c^{N}(\underline{\gamma},\underline{y}, \omega_{0})}.
\end{flalign*}
\end{theorem}

Here, we present an example where Theorem \ref{the:2} can be applied but the existence result of \cite[Theorem 12]{sanjari2018optimal} cannot be applied because the assumption that $\mathbb{U}$ is convex in \cite[Theorem 12]{sanjari2018optimal} is violated.
\begin{example}\label{ex:1}
Consider a team problem with the following expected cost function
\begin{equation*}
J({\underline{\gamma}})=\limsup\limits_{N \to \infty}\mathbb{E}^{\underline{\gamma}}\bigg[\bigg((\frac{1}{N}\sum_{i=1}^{N}u^{i})-\frac{1}{2}\bigg)^{2}\bigg],
\end{equation*}
where $\sigma$-field $\sigma(y^{i})=\{\emptyset, \Omega\}$ (this corresponds to a team setup where DMs have no measurement, hence measurable functions (policies) are constant functions), and we consider $u^{i}\in \{0,1\}$ for each DM.
Clearly, an optimal policy that achieves zero is the one where a matching partition (such as even numbers vs. odd numbers) among DMs picking $u^{i}=0$ and $u^{i}=1$, that is because the cost function is non-negative. One can see that there is an optimal policy  in $L_{\text{PR},\text{SYM}}$ since each DM can choose independently an action zero or one with probability half and this achieves the expected cost of zero; however, there is no identically deterministic policy that achieves zero expected cost. We note also that the problem is not a convex problem, therefore the results in \cite[Theorem 12 or Proposition 1]{sanjari2018optimal} are not applicable to show the existence of a symmetric randomized optimal policy, in particular, the action sets are not convex.     
\end{example}

\section{Finite Horizon Dynamic Team Problems with a Symmetric Information Structure}\label{sec:4}
In this section, we study dynamic stochastic team problems. All the proofs regarding this section are presented in Appendix \ref{appB}. Similar to the static case, we first introduce the intrinsic model for general dynamic team problems, and then, we introduce a topology on control policies, and finally we establish our main results for dynamic problems. 

\subsection{A Revised intrinsic model for dynamic team problems}
Under the intrinsic model (see Section \ref{sec:pre}), every DM acts separately. However, depending on the information structure, it may be convenient to consider a collection of DMs as a single DM acting at different time instances. In fact, in the classical stochastic control, this is the standard approach.  In this subsection, we introduce the general (multi-stage) dynamic problems using the intrinsic model under deterministic policies. In the next subsections, we allow randomization equipped with a suitable topology.

According to the discussion above, by considering a collection of DMs as a single DM ($i=1,\dots,N$) acting at different time instances ($t=0,\dots,T-1$), we revise the intrinsic mode for (multi-stage) dynamic team problems with $(NT)$-DMs as a team with $N$-DMs (for $N\in \mathbb{N}\cup \{\infty\}$):

\begin{itemize}[wide]
\item  [{(i)}] Let the observation and action spaces be standard Borel spaces and be identical for each DM ($i=1,\dots,N$) with $\mathbb{Y}_{i}:={\bf{Y}}=\prod_{t=0}^{T-1}\mathbb{Y}^{t}$, $\mathbb{U}_{i}:={\bf{U}}=\prod_{t=0}^{T-1}\mathbb{U}^{t}$, respectively (later on, for simplicity of our notation and analysis, we assume that action and observation spaces are also identical through time). For each DM$^{i}$,  the set of all admissible policies are denoted by $\Gamma_{i}:=\prod_{t=0}^{T-1} \Gamma^{t}$. Later on, these policies will be allowed to be randomized and accordingly the image will be ${\mathcal P}({\bf{U}})$.
\item [(ii)] For $i=1,\dots,N$, $y^{i}_{t}:=h_{t}^{i}(x_{0}^{1:N}, \zeta^{1:N}_{0:t}, u_{0:t-1}^{1:N})$ represents the observation of DM$^i$ at time $t$ ($h_{t}^{i}$s are Borel measurable functions). Let ${\nu}_{t}^{N}$ be a stochastic kernel characterizing the joint distribution of observations  $y^{1:N}_{t}:=(y_{t}^{1},\dots,y_{t}^{N})$ at time $t$ induced by $h_{t}^{i}$s given the available information, and let $(\underline{\zeta}^{1:N}):=(\underline{\zeta}^{1},\dots,\underline{\zeta}^{N})$ where $\underline\zeta^{i}:=(x_{0}^{i}, \zeta_{0:T-1}^{i})$ denotes all the uncertainty associated with DM$^i$ including his/her initial states. We assume that $(\underline{\zeta}^{i})$ takes values in $\Omega_{\zeta}$ (where at each time instances $t$, it takes value in $\Omega_{\zeta_{t}}$). Let $\mu^{N}$ denote the law of $\underline{\zeta}^{1:N}$. {To be consistent with our notations in our analysis of the static case, we used the same notation $\mu^{N}$ as the fixed probability measures on observations and $\omega_0$ for the static case; however, we note that in the dynamic case, the probability measures on uncertainties $(\underline{\zeta}^{1:N})$ is fixed and not probability measures on observations.}
\end{itemize}

\subsection{Topology on dynamic control policies}
Similar to Section \ref{sec:strategic}, here, we allow randomization in policies, but first we introduce two reduction conditions (independent and nested reduction) that enable us to define sets of Borel probability measures on randomized policies for dynamic teams with different information structures by considering a policy of a single DM ($i=1,\dots,N$) acting at different time instances ($t=0,\dots,T-1$).

\begin{assumption}\label{assump:ind1}
One of the following conditions holds:
\begin{itemize}[wide]
\item [(i)] (Independent reduction): for every $N\in \mathbb{N}\cup\{\infty\}$ and for $i=1,\dots, N$ and $t=0,\dots,T-1$, there exists a probability measure $\tau_{t}^{i}$ on $\mathbb{Y}^{t}$ and a function $\psi_{t}^{i}:\mathbb{Y}^{t}\times \Omega_{0}\times \prod_{p=1}^{N}(\prod_{k=0}^{t-1}\Omega_{\zeta_{k}}\times\prod_{k=0}^{t-1}(\mathbb{U}^{k}\times \mathbb{Y}^{k}))\to \mathbb{R}_{+}$ such that for all Borel sets $A^{i}$ on $\mathbb{Y}^{t}$ (with $A=A^1 \times \dots \times A^N$)
\begin{flalign*}
&\scalemath{0.95}{\nu_{t}^{N}(A|\omega_{0}, x_{0}^{1:N}, \zeta_{0:t-1}^{1:N}, y^{1:N}_{0:t-1}, u_{0:t-1}^{1:N})}\scalemath{0.95}{=\prod_{i=1}^{N}\int_{A^i}\psi_{t}^{i}(y^{i}_{t}, \omega_{0}, x_{0}^{1:N},\zeta_{0:t-1}^{1:N},y_{0:t-1}^{1:N}, u_{0:t-1}^{1:N})\tau_{t}^{i}(dy_{t}^{i})}.
\end{flalign*}
\item [(ii)] (Nested reduction): for every $N\in \mathbb{N}\cup\{\infty\}$ and for $i=1,\dots, N$ and $t=0,\dots,T-1$, there exists a probability measure $\eta^{i}_{t}$ on $\mathbb{Y}^{t}$ and a function $\phi^{i}_{t}$ such that for all Borel sets $A^{i}$ on $\mathbb{Y}^{t}$ (with $A=A^1 \times \dots \times A^N$)
\begin{flalign*}
&\scalemath{0.95}{\nu_{t}^{N}(A|\omega_{0}, x_{0}^{1:N}, \zeta_{0:t-1}^{1:N}, y^{1:N}_{0:t-1}, u_{0:t-1}^{1:N})}\\
&\scalemath{0.95}{=\prod_{i=1}^{N}\int_{A^i}\phi_{t}^{i}(y^{i}_{t}, \omega_{0}, x_{0}^{-i},\zeta_{0:t-1}^{-i},y_{0:t-1}^{-i}, u_{0:t-1}^{-i})\eta_{t}^{i}(dy_{t}^{i}|x_{0}^{i}, \zeta_{0:t-1}^{i}, y^{i}_{0:t-1}, u_{0:t-1}^{i} )}\nonumber,
\end{flalign*}
and for each DM$^{i}$ through time ($t=0,\dots,T-1$), there exists a static reduction with the classical information structure (i.e., under the reduction, the information structure of each DM through time is expanding such that $\sigma(y_{t}^{i}) \subset \sigma(y_{t+1}^{i})$ for $t=0,\dots,T-1$).
\end{itemize}
\end{assumption}

We note that Assumption \ref{assump:ind1}(i) allows us to obtain an independent measurements reduction both through DMs and through time, $t=0,\dots,T-1$ (see Appendix \ref{sec:statreduc}).  Assumption \ref{assump:ind1}(ii) holds if an independent static reduction exists through DMs and there exists a nested static reduction for each DM through time, i.e., under the reduction, the information is expanding for each DM through time (see Appendix \ref{sec:statreduc}). In view of the above reduction conditions, we introduce a suitable topology for randomized policies. Similar to Section \ref{sec:strategic}, under Assumption \ref{assump:ind1}(i), we define convergence on policies as: 
\begin{flalign*}
\underline{\gamma}^{i}_{n}  \xrightarrow{n \to \infty} \underline{\gamma}^{i}\:\:\: \text{if and only if}\:\:\: {\gamma}^{i}_{t, n}(d{u}^{i}_{t}|{y}^{i}_{t})\tau^{i}_{t}(d{y}^{i}_{t}) \xrightarrow[\text{weakly}]{n \to \infty} {\gamma}^{i}_{t}(d{u}^{i}_{t}|{y}^{i}_{t})\tau^{i}_{t}(d{y}^{i}_{t})~~\forall ~~t=0,\dots,T-1.
\end{flalign*}
Under Assumption \ref{assump:ind1}(ii), we define convergence on policies as:
 \begin{flalign*}
\underline{\gamma}^{i}_{n}  \xrightarrow{n \to \infty} \underline{\gamma}^{i}\:\:\: \text{if and only if}\:\:\: {\gamma}^{i}_{t, n}(d{u}^{i}_{t}|{y}^{i}_{0:t})\eta^{i}_{t}(d{y}^{i}_{0:t}) \xrightarrow[\text{weakly}]{n \to \infty} {\gamma}^{i}_{t}(d{u}^{i}_{t}|{y}^{i}_{0:t})\eta^{i}_{t}(d{y}^{i}_{0:t})~~\forall ~~t=0,\dots,T-1.
\end{flalign*}
 Hence, under Assumption \ref{assump:ind1}, we define all the sets of randomized policies defined in Section \ref{sec:strategic} for the dynamic teams by considering $\underline{\gamma}^{i}$. 
 
{\begin{remark}
We note that our first reduction condition, independent reduction, is essentially a version of Girsanov's transformation \cite{girsanov1960transforming, benevs1971existence} which has been considered first in \cite[Eqn(4.2)]{wit88}, and later utilized in \cite[p. 114]{YukselBasarBook} and \cite[Section 2.2]{yuksel2018general} (for discrete-time partially observed stochastic control, similar arguments had been presented, e.g. by Borkar in \cite{Bor00}, \cite{Bor07}). We refer the reader to \cite{charalambous2016decentralized} for relations with the classical continuous-time stochastic control, where the relation with Girsanov's classical measure transformation \cite{girsanov1960transforming, benevs1971existence} is recognized. Our second reduction condition, nested reduction, holds when there exists a reduction for DMs through time under which each DM has a perfect recall of private history of information.
\end{remark} }

{Now, we provide examples under which either one of the conditions in Assumption  \ref{assump:ind1} holds.
 \begin{example}\label{lem:exmple}
For each $i=1,\dots, N$ and $t=0,\dots, T-1$, let
$x_{t+1}^{i}=f_{t}^{i}(x_{0:t}^{1:N},u_{0:t}^{1:N},w_{t}^{i})$ and $y_{t}^{i}={h}_{t}^{i}(\omega_0, x_{0:t}^{1:N},\zeta_{0:t-1}^{1:N}, u_{0:t-1}^{1:N})+v_{t}^{i}$, where $\zeta_{t}^{i}:=(w_{t}^{i},v_{t}^{i})$, and $v_{t}^{i}$ admits zero-mean Gaussian density function $\theta_{t}^{i}$ with positive-definite covariance. 
\begin{itemize}[wide]
\item [(i)] If the information structure for each DM at time $t$ is described as $I_{t}^{i}:=\{y_{t}^{i}\}$ for all $i=1,\dots, N$ and $t=0,\dots,T-1$, then Assumption \ref{assump:ind1}(i) holds.
\item[(ii)] If $I_{t}^{i}:=\{y_{0:t}^{i}, u_{0:t-1}^{i}\}$ for all $i=1,\dots, N$ and $t=0,\dots,T-1$ (or equivalently, $I_{t}^{i}:=\{\tilde{y}_{t}^{i}\}$ with $\tilde{y}_{t}^{i}:=\tilde{h}_{t}^{i}(\omega_0, x_{0:t}^{1:N},\zeta_{0:t-1}^{1:N}, u_{0:t-1}^{1:N}, v_{0:t}^{i})$ for some function $\tilde{h}_{t}^{i}$ and $\sigma(\tilde{y}_{t}^{i})\subset \sigma(\tilde{y}_{t+1}^{i})$  and $\sigma(u_{t}^{i})\subset \sigma(\tilde{y}_{t+1}^{i})$ for some function $\bar{h}_{t}$), then Assumption \ref{assump:ind1}(ii) holds.
\end{itemize}

 Part(i) is true since for all $t=0,\dots,T-1$ and $i=1,\dots, N$, we have
\begin{flalign*}
y_{t}^{i}&={h}_{t}^{i}(\omega_{0}, x_{0:t}^{1:N},\zeta_{0:t-1}^{1:N}, u_{0:t-1}^{1:N})+v_{t}^{i}=\kappa_{t}^{i}(\omega_{0}, x_{0}^{1:N},\zeta_{0:t-1}^{1:N}, u_{0:t-1}^{1:N})+v_{t}^{i},
\end{flalign*}
for some functions $\kappa_{t}^{i}$, and hence, we can define
\begin{flalign*}
&\psi_{t}^{i}(y^{i}_{t}, \omega_{0}, x_{0}^{1:N},\zeta_{0:t-1}^{1:N},y_{0:t-1}^{1:N}, u_{0:t-1}^{1:N}):=\frac{\theta_{t}^{i}(y_{t}^{i}-\kappa_{t}^{i}(\omega_{0}, x_{0}^{1:N},\zeta_{0:t-1}^{1:N}, u_{0:t-1}^{1:N}))}{\theta_{t}^{i}(y_{t}^{i})},\\
&\tau_{t}^{i}(dy_{t}^{i}):=\theta_{t}^{i}(y_{t}^{i})dy_{t}^{i}.
\end{flalign*}
Part(ii) can be shown similarly by first applying the independent reduction as above, and then, considering the nested information structure through time for each DM.
\end{example}

\begin{example}
Consider the following two information structures:
\begin{itemize}[wide]
\item [(i)] (Open-loop information structure): For each $i=1,\dots, N$ and $t=0,\dots, T-1$, let 
$x_{t+1}^{i}=f_{t}^{i}(x_{0:t}^{1:N},u_{0:t}^{1:N},w_{t}^{i})$ and $y_{t}^{i}={h}_{t}^{i}(\zeta_{0:t-1}^{i}, v_{t}^{i})$ such that $\sigma(y_{t}^{i})\subset \sigma(y_{t+1}^{i})$, where $(\zeta_{t}^{i})_{t}:=(w_{t}^{i},v_{t}^{i})_{t}$ denotes the disturbances of DM$^{i}$ (which is independent of disturbances of other DMs and independent of $\omega_{0}$). If $I_{t}^{i}:=\{y_{t}^{i}\}$ for all $i=1,\dots, N$ and $t=0,\dots,T-1$, then Assumption \ref{assump:ind1}(ii) holds.
\item[(ii)] For each $i=1,\dots, N$ and $t=0,\dots, T-1$, let
$x_{t+1}^{i}=f_{t}^{i}(\omega_{0},x_{0:t}^{1:N},u_{0:t}^{1:N})+w_{t}^{i}$, where $w_{t}^{i}$ admits zero-mean Gaussian density function $\theta_{t}^{i}$ with positive-definite covariance, and let $y_{t}^{i}={h}_{t}^{i}(x_{0:t}^{i}, y_{0:t-1}^{i}, v_{0:t}^{i})$ such that $\sigma(y_{t}^{i})\subset \sigma(y_{t+1}^{i})$, where $(v_{t}^{i})_{t}$ are independent of disturbances of other DMs and independent of $\omega_{0}$. If $I_{t}^{i}:=\{y_{t}^{i}\}$ for all $i=1,\dots, N$ and $t=0,\dots,T-1$, then Assumption \ref{assump:ind1}(ii) holds.
\end{itemize}
 Part(i) follows from the fact that the information structure is open-loop and nested for each DM, and hence, under this information structure the problem is static with the classical information structure through time for each DM.  Part(ii) is true since for all $t=0,\dots,T-1$ and $i=1,\dots, N$, 
\begin{flalign*}
&\hat{\phi}_{t}^{i}(x^{i}_{t}, \omega_{0}, x_{0:t-1}^{1:N}, u_{0:t-1}^{1:N}):=\frac{\theta_{t}^{i}(x_{t}^{i}-f_{t}^{i}(\omega_{0}, x_{0:t-1}^{1:N}, u_{0:t-1}^{1:N}))}{\theta_{t}^{i}(x_{t}^{i})},\\
& \hat{\eta}_{t}^{i}(dx_{t}^{i}):=\theta_{t}^{i}(x_{t}^{i})dx_{t}^{i},
\end{flalign*}
and since the information structure is nested through time for each DM. 
\end{example}
}

\subsection{Existence and structure of optimal policies for symmetric dynamic team problems with infinitely many decision makers}\label{sec:mftdy2}

In the following, we study the existence and structure of globally optimal policies for dynamic team problems with a symmetric information structure (that are not necessarily partially nested) and with a finite but large and also infinitely many decision makers. We note that a related result is given in \cite{sanjari2019optimal} where convex mean-field team problems have been considered under the assumption that the action space is convex for each DM and the cost function is convex in policies. We note that even if the cost function is convex in actions when there is a mean-field coupling in dynamics, convexity rarely holds since the information structure under decentralized setup is non-classical, and that may lead to the non-convexity of the team problem in policies (see for example \cite[Section 3.3]{YukselSaldiSICON17}). In the following, convexity is not imposed. Again, for our results in this subsection, we impose Assumption \ref{assump:ind1}.

\subsubsection{Exchangeability of optimal policies for symmetric dynamic team problems with a finite but large number of decision makers}

In this subsection, we focus on symmetric dynamic team problems with $N$-DMs, and we establish a structural result for optimal policies of this class of problems (which is more general than the pre-limit mean-field model ($\mathcal{P}_{T}^{N}$)). In the next subsection, we use this result to establish existence and structural properties of globally randomized optimal policies for mean-field dynamic team problems.

 Now, we recall the definition of the symmetric information structure from \cite{sanjari2019optimal} (note that symmetric information structures can be classical, partially nested, or non-classical). Several examples as well as a graph interpretation of dynamic teams with symmetric information structures have been presented in \cite[Section 4]{sanjari2019optimal}. In particular, pre-limit mean-field and mean-field dynamic team problems ($\mathcal{P}_{T}^{N}$) and ($\mathcal{P}_{T}^{\infty}$) introduced in Section \ref{sec:2.2} have a symmetric information structure.
\begin{definition}\label{def:s}\cite{sanjari2019optimal}
Let the information of DM$^{i}$ acting at time $t$ be described as $I_{t}^{i}:=\{y_{t}^{i}\}$. The information structure of a sequential $N$-DM team problem is \textit{symmetric} if  
\begin{itemize}
\item[(i)] $y^{i}_{t}=h_{t}(x_{0}^{i},x_{0}^{-i}, \zeta^{i}_{0:t}, {\zeta}^{-i}_{0:t}, u_{0:t-1}^{i}, u_{0:t-1}^{-i})$, where $h_{t}$ is identical for all $i=1,\dots,N$ (note that the arguments of the function depend on $i$) and $b^{-i}=(b^{1},\dots, b^{i-1},b^{i+1},\dots, b^{N})$ for $b=x_{0}, {\zeta}_{0:t},  u_{0:t-1}$.
\end{itemize}
\end{definition}

{We note that the above definition can be generalized to be applicable for teams with countably infinite number of DMs. Before, we present the result for dynamic mean-field teams,  we characterize team problems with symmetric information structures in which the search for an optimal policy can be restricted to policies in $L_{\text{EX}}^{N}$ without losing global optimality.
To this end, we focus on a more general setup of team problems within randomized policies $P_{\pi}\in L^{N}$ as
\begin{flalign}
\scalemath{0.95}{J_{N}^{\pi}(\underline{\gamma}^{1:N})}&\scalemath{0.95}{:=\int P_{\pi}(d\underline{\gamma})\mu^{N}(d\omega_{0},d\underline{\zeta})c^{N}(\underline{\zeta}, \underline{\gamma},\underline{y}, \omega_{0})\nu^{N}(d\underline{y}|\underline{\zeta}, \underline{\gamma},\omega_{0})}\nonumber\\
&\scalemath{0.95}{:=\int \bigg(\int c(\omega_{0},\underline{\zeta}^{1:N},\underline{u}^{1},\dots, \underline{u}^{N})\prod_{i=1}^{N}\underline{\gamma}^{i}(d\underline{u}^{i}|\underline{y}^{i})\bigg)P_{\pi}(d\underline\gamma^{1},\dots,d\underline\gamma^{N})\mu^{N}(d\omega_{0}, d\underline\zeta^{1:N})}\label{eq:finitecostd}\\
&\:\:\:\:\:\scalemath{0.95}{\times \prod_{t=0}^{T-1}{{\nu}_{t}^{N}}\left(d{y}^{1:N}_{t}\middle|\omega_{0}, x_{0}^{1:N},{\zeta}^{1:N}_{0:t-1},{y}_{0:t-1}^{1:N},{u}^{1:N}_{0:t-1}\right)}\nonumber,
\end{flalign}
where $c^{N}(\underline{\zeta},\underline{\gamma},\underline{y}, \omega_{0}):=\int c(\omega_{0},\underline{\zeta}^{1:N},\underline{u}^{1},\dots, \underline{u}^{N})\prod_{i=1}^{N}\underline{\gamma}^{i}(d\underline{u}^{i}|\underline{y}^{i})$ and the following assumptions hold.
\begin{assumption}\label{assump:3.1}
For any permutation $\sigma$ of the set $\{1,\dots,N\}$, we have for all $\omega_{0}$,
 \begin{flalign}\label{eq:3.3}
c(\omega_{0}, (\underline{\zeta}^{\sigma})^{1:N},(\underline{u}^{\sigma})^{1:N})=c(\omega_{0}, \underline{\zeta}^{1:N},\underline{u}^{1:N}),
 \end{flalign}
 where $(\underline{\zeta}^{\sigma})^{1:N}:=(\underline{\zeta}^{\sigma(1)},\dots,\underline{\zeta}^{\sigma(N)})$ and $(\underline{u}^{\sigma})^{1:N}:=(\underline{u}^{\sigma(1)},\dots,\underline{u}^{\sigma(N)})$.
\end{assumption}

\begin{assumption}\label{assump:sxcdy}
\hfill 
\begin{itemize}
\item [(a)] $(\underline\zeta^{1}, \dots,\underline\zeta^{N})$ are exchangeable conditioned on $\omega_{0}$; 
\item  [(b)] For all $t=0,\dots, T-1$, and all Borel sets $A^{i}$ on ${{\mathbb{Y}^{t}}}$ (with $A=A^{1}\times \dots \times A^{N}$) 
\end{itemize}
\begin{flalign}
&{\nu_{t}^{N}}\left(A\middle|\omega_{0}, x_{0}^{1:N},{\zeta}^{1:N}_{0:t-1},{y}_{0:t-1}^{1:N},{u}^{1:N}_{0:t-1}\right)=\prod_{i=1}^{N}{\nu_{t}^{i}}\left(A^{i}\middle|\omega_{0}, x_{0}^{i},{\zeta}^{i}_{0:t-1},{y}_{0:t-1}^{1:N},{u}^{1:N}_{0:t-1}\right)\nonumber,
\end{flalign}
where $\nu_{t}^{i}$ is a stochastic kernel of the observation DM$^{i}$ at time $t$, $y_{t}^{i}$, induced by $h_{t}$ (which is identical for each DM).
\end{assumption}

We note that dynamic mean-field team problems introduced in Section \ref{sec:dn} with the cost function \eqref{eq:mfcost}, dynamic \eqref{eq:mfdynamics2}, and observations \eqref{eq:mfobs2}, under Assumption \ref{assump:2}  satisfy Assumption \ref{assump:3.1} and Assumption \ref{assump:sxcdy}.

\begin{lemma}\label{lem:3.1dy}
Consider a dynamic team problem with a symmetric information structure. Let Assumption \ref{assump:3.1} and Assumption \ref{assump:sxcdy} hold. Assume $\bar{L}^{N}$ is an arbitrary convex subset of $L^{N}$. Then 
\begin{flalign*}
&\scalemath{0.95}{ \inf\limits_{P_{\pi} \in \bar{L}^{N}}\int P_{\pi}(d\underline{\gamma})\mu^{N}(d\omega_{0},d\underline{\zeta})c^{N}( \underline{\zeta}, \underline{\gamma},\underline{y}, \omega_{0})\nu^{N}(d\underline{y}|\underline{\zeta}, \underline{\gamma},\omega_{0})}\\
&{=\inf\limits_{P_{\pi} \in \bar{L}^{N} \cap L_{\text{EX}}^{N}}\int P_{\pi}(d\underline{\gamma})\mu^{N}(d\omega_{0},d\underline{\zeta})c^{N}( \underline{\zeta}, \underline{\gamma},\underline{y}, \omega_{0})\nu^{N}(d\underline{y}|\underline{\zeta}, \underline{\gamma},\omega_{0})}.
\end{flalign*}
\end{lemma}
}
\subsubsection{Existence and structure of optimal policies for mean-field dynamic team problems}

In the following, we establish the existence of a globally randomized optimal policy for dynamic mean-field teams with infinitely many decision makers.  Define state dynamics and observations as \eqref{eq:mfdynamics2} and \eqref{eq:mfobs2}. The information structure of DM$^{i}$ at time $t$ is $I_{t}^{i}=\{{y}^{i}_{t}\}$, and $\zeta_{t}^{i}:=(w_{t}^{i}, v_{t}^{i})$ (with $\zeta_{0}^{i}:=(x_{0}^{i}, w_{0}^{i}, v_{0}^{i})$) denotes the uncertainty corresponding to dynamics and observations at time $t$ for DM$^i$ which are exogenous random vectors in standard Borel spaces. First, we re-state the infinite decision maker mean-field team problem and its pre-limit within randomized policies.
\begin{itemize} [wide]
\item[\bf\text{Problem} \bf{($\mathcal{P}_{T}^{N}$)}:]
Consider an $N$-DM dynamic team with the expected cost of a randomized policy $P_{\pi}^{N}\in L^{N}$ as:
\begin{flalign}
&\scalemath{0.95}{\int P_{\pi}^{N}(d\underline{\gamma})\mu^{N}(d\omega_{0},d\underline{\zeta})c^{N}(\underline{\zeta}, \underline{\gamma},\underline{y}, \omega_{0})\nu^{N}(d\underline{y}|\underline{\zeta}, \underline{\gamma},\omega_{0})}\nonumber\\
&\scalemath{0.95}{:=\int \bigg(\int \frac{1}{N}\sum_{t=0}^{T-1}\sum_{i=1}^{N}c\bigg(\omega_{0}, x_{t}^{i},u_{t}^{i},\frac{1}{N}\sum_{p=1}^{N}u_{t}^{p}(y_{t}^{p}),\frac{1}{N}\sum_{p=1}^{N}x_{t}^{p}\bigg)\prod_{k=1}^{N}\underline{\gamma}^{k}(d\underline{u}^{k}|\underline{y}^{k})\bigg)}\label{eq:finiterandom}\\
&\:\scalemath{0.95}{\times P_{\pi}^{N}(d\underline\gamma^{1},\dots,d\underline\gamma^{N})\mu^{N}(d\omega_{0}, d\underline\zeta^{1:N})\prod_{t=0}^{T-1}{{\nu}^{N}_{t}}\left(d{y}^{1:N}_{t}\middle|\omega_{0}, x_{0}^{1:N},{\zeta}^{1:N}_{0:t-1},{y}_{0:t-1}^{1:N},{u}^{1:N}_{0:t-1}\right)}\nonumber,
\end{flalign}
\end{itemize}
where 
\begin{flalign*}
c^{N}(\underline{\zeta}, \underline{\gamma},\underline{y}, \omega_{0}):=\int \frac{1}{N}\sum_{t=0}^{T-1}\sum_{i=1}^{N}c\bigg(\omega_{0}, x_{t}^{i},u_{t}^{i},\frac{1}{N}\sum_{p=1}^{N}u_{t}^{p}(y_{t}^{p}),\frac{1}{N}\sum_{p=1}^{N}x_{t}^{p}\bigg)\prod_{k=1}^{N}\underline{\gamma}^{k}(d\underline{u}^{k}|\underline{y}^{k}).
\end{flalign*}

The above problem is considered as a pre-limit problem for our infinite-decision maker team problem.  {We note that $N$-DM teams of ($\mathcal{P}_{T}^{N}$) is a special case of \eqref{eq:finitecostd} since we have a special structure for the cost function $c^{N}$ and observations which satisfy Assumption \ref{assump:3.1} and Definition \ref{def:s}, respectively.} 
{
\begin{remark}
Our analysis below also allows a more general observations for each DM where the observations of each DM at time $t$ can be explicitly functions of average of previous states and actions as 
\begin{flalign*}
&y_{t}^{i}=h_{t}\bigg(x_{0:t}^{i},u_{0:t-1}^{i},\frac{1}{N}\sum_{p=1}^{N}x_{0:t-1}^{p}, \frac{1}{N}\sum_{p=1}^{N}u_{0:t-1}^{p}, v_{0:t}^{i}\bigg).
\end{flalign*}
However, to simplify the presentations of theorems and proofs and emphasize in the decentralization of optimal policy, for the rest of the paper, we consider \eqref{eq:mfobs2}. 
\end{remark}
}
\begin{itemize} [wide]
\item[\bf\text{Problem} \bf{($\mathcal{P}^{\infty}_{T}$)}:]
Consider infinite-DM static team with the following expected cost of a randomized policy $P_{\pi}\in L$ as:
\begin{flalign}
&\scalemath{0.95}{\limsup\limits_{N \to \infty}\int P_{\pi, N}(d\underline{\gamma})\mu^{N}(d\omega_{0},d\underline{\zeta})c^{N}(\underline{\zeta}, \underline{\gamma},\underline{y}, \omega_{0})\nu^{N}(d\underline{y}|\underline{\zeta}, \underline{\gamma},\omega_{0}),}\label{eq:infiniterandom}
\end{flalign}
where $P_{\pi, N}$ is the restriction of $P_{\pi}\in L$ to its first $N$ components and $\mu^{N}$ is the marginal of the fixed probability measure on $(\omega_{0}, \underline{\zeta}^{1},\underline{\zeta}^{2},\dots)$ to the first $N+1$ components.
\end{itemize}

\begin{assumption}\label{assump:ind2}
Assume Assumption \ref{assump:ind1} holds with functions $\psi_{t}^{i}$ and $\phi_{t}^{i}$ are of the following forms for every $i\in \mathcal{N}$ and $t=0,\dots, T-1$: 
\begin{flalign*}
&\scalemath{0.95}{\psi_{t}^i\bigg(y^{i}_{t}, \omega_{0}, x_{0}^{i},\zeta_{0:t-1}^{i},y_{0:t-1}^{i}, u_{0:t-1}^{i}, \frac{1}{N}\sum_{p=1}^{N}u_{0:t-1}^{p}, \frac{1}{N}\sum_{p=1}^{N}x_{0:t}^{p}\bigg)},\\
&\scalemath{0.95}{\phi_{t}^i\bigg(y^{i}_{t}, \omega_{0}, \frac{1}{N}\sum_{p=1}^{N}u_{0:t-1}^{p}, \frac{1}{N}\sum_{p=1}^{N}x_{0:t}^{p}\bigg)},
\end{flalign*}
where $\psi_{t}^{i}$ is continuous in the last three arguments (actions and the empirical mean of actions and states) and $\phi_{t}^{i}$ is continuous in the last two arguments (the empirical mean of actions and states).

\end{assumption}

Before presenting our main result for dynamic mean-field teams, we present sufficient conditions under which the expected cost function induced by randomized optimal policies in $L^{N}_{\text{EX}}$ and $L_{\text{EX}}$ are equal as $N$ goes to infinity, and hence, following from Lemma \ref{lem:3.1dy}, under symmetry, this shows that  without loss of global optimality, optimal policies of dynamic mean-field teams can be considered to be an infinitely-exchangeable type. 

\begin{lemma}\label{lem:findefdyn2}
Consider the team problem ($\mathcal{P}_{T}^{N}$) where Assumption \ref{assump:ind2}, Assumption \ref{assump:c}, and Assumption \ref{assump:2} hold. Assume further that the cost function bounded, then
\begin{flalign}
&\limsup\limits_{N \to \infty}\inf\limits_{P_{\pi}^{N} \in L^{N}_{\text{EX}}}\int P_{\pi}^{N}(d\underline{\gamma})\mu^{N}(d\omega_{0},d\underline{\zeta})c^{N}( \underline{\zeta}, \underline{\gamma},\underline{y}, \omega_{0})\nu^{N}(d\underline{y}|\underline{\zeta}, \underline{\gamma},\omega_{0})\nonumber\\
&=\limsup\limits_{N \to \infty}\inf\limits_{P_{\pi} \in L_{\text{EX}}}\int P_{\pi, N}(d\underline{\gamma})\mu^{N}(d\omega_{0},d\underline{\zeta})c^{N}( \underline{\zeta}, \underline{\gamma},\underline{y}, \omega_{0})\nu^{N}(d\underline{y}|\underline{\zeta}, \underline{\gamma},\omega_{0})\label{eq:ertdy2},
\end{flalign}
where $P_{\pi, N}$ is the restriction of $P_{\pi}\in L_{\text{EX}}$ to its first $N$ components and $\mu^{N}$ is the marginal of the fixed probability measure on $(\omega_{0}, \underline{\zeta}^{1},\underline{\zeta}^{2},\dots)$ to the first $N+1$ components.
\end{lemma}

In the following, we establish an existence and structural result for a randomized optimal policy of ($\mathcal{P}^{\infty}_{T}$). 
\begin{theorem}\label{the:exmftdy2}
Consider a mean-field team problem ($\mathcal{P}_{T}^{\infty}$) with ($\mathcal{P}_{T}^{N}$) having a symmetric information structure for every $N$. Let Assumption \ref{assump:ind2}, Assumption \ref{assump:c}, and Assumption \ref{assump:2} hold. Then, there exists a randomized optimal policy $P^{*}_{\pi}$ for {(${P}^{\infty}_{T}$)} which is  in $L_{\text{PR,SYM}}$,
\begin{flalign*}
&\scalemath{0.95}{\inf\limits_{P_{\pi} \in L_{\text{PR,SYM}}} \limsup\limits_{N \to \infty} \int P_{\pi, N}(d\underline{\gamma})\mu^{N}(d\omega_{0},d\underline{\zeta})c^{N}( \underline{\zeta}, \underline{\gamma},\underline{y}, \omega_{0})\nu^{N}(d\underline{y}|\underline{\zeta}, \underline{\gamma},\omega_{0})}\\
&\scalemath{0.95}{:=\limsup\limits_{N \to \infty} \int P_{\pi, N}^{*}(d\underline{\gamma})\mu^{N}(d\omega_{0},d\underline{\zeta})c^{N}( \underline{\zeta}, \underline{\gamma},\underline{y}, \omega_{0})\nu^{N}(d\underline{y}|\underline{\zeta}, \underline{\gamma},\omega_{0})}\\
&\scalemath{0.95}{=\inf\limits_{P_{\pi} \in L} \limsup\limits_{N \to \infty} \int P_{\pi, N}(d\underline{\gamma})\mu^{N}(d\omega_{0},d\underline{\zeta})c^{N}( \underline{\zeta}, \underline{\gamma},\underline{y}, \omega_{0})\nu^{N}(d\underline{y}|\underline{\zeta}, \underline{\gamma},\omega_{0})}.
\end{flalign*}
\end{theorem}
\section{Approximations of Optimal Policies for Symmetric $N$-DM Stochastic Team Problems}\label{sec:aproxi}
In this section, we present approximation results of optimal policies for $N$-DM team problems. We show that for large $N$, symmetric policies are nearly optimal and the restriction of the optimal infinite solution to the finite team problem is nearly optimal for large $N$. All the proofs regarding this section are presented in Appendix \ref{appC}. We first consider the static case. {To present ideas more effectively, we first introduce the following set of probability measures on policies as:
\begin{flalign*}
\scalemath{0.95}{L_{\text{D}}^{N}:=\bigg\{P_{\pi} \in L^{N}\bigg{|}\text{for all}~A_{i} \in \mathcal{B}(\Gamma^{i}): P_{\pi}(\gamma^{1} \in A_{1},\dots, \gamma^{N} \in A_{N})=\prod_{i=1}^{N}1_{\{\tilde\gamma^{i}\in A_{i}\}},~ \text{for}~\tilde{\gamma}^{i}\in \Gamma^{i}\bigg\}}\nonumber,
\end{flalign*}
where the above set corresponds to Dirac-delta measures  in $L_{\text{PR}}^{N}$.} 
\begin{theorem}\label{the:3}
Consider a static team problem ($\mathcal{P}_{N}$) (see \eqref{eq:pf}) where  Assumption \ref{assump:cont} and Assumption \ref{assump:oc} hold. Assume further the cost function is bounded. Then,
\begin{itemize}[wide]
\item [(i)]
\begin{flalign}
\scalemath{0.95}{\inf\limits_{P_{\pi}^{N} \in L_{\text{PR,SYM}}^{N}}\int P_{\pi}^{N}(d\underline{\gamma})\mu^{N}(d\omega_{0},d\underline{y}) c^{N}(\underline{\gamma},\underline{y}, \omega_{0})\leq \inf\limits_{P_{\pi}^{N} \in L_{\text{CO}}^{N}}\int P_{\pi}^{N}(d\underline{\gamma})\mu^{N}(d\omega_{0},d\underline{y}) c^{N}(\underline{\gamma},\underline{y}, \omega_{0})+ \epsilon_{N}}\label{eq:approx},
\end{flalign}
and
\begin{flalign}
\scalemath{0.95}{\inf\limits_{P_{\pi}^{N} \in L_{\text{PR,SYM}}^{N}}\int P_{\pi}^{N}(d\underline{\gamma})\mu^{N}(d\omega_{0},d\underline{y}) c^{N}(\underline{\gamma},\underline{y}, \omega_{0})\leq \inf\limits_{P_{\pi}^{N} \in L_{\text{D}}^{N}}\int P_{\pi}^{N}(d\underline{\gamma})\mu^{N}(d\omega_{0},d\underline{y}) c^{N}(\underline{\gamma},\underline{y}, \omega_{0})+\epsilon_{N}}\label{eq:approx1},
\end{flalign}
where $\epsilon_{N} \to 0$ as $N$ goes to infinity.
\item [(ii)]  If $P_{\pi}^{*}\in L_{\text{PR,SYM}}$ is a randomized optimal policy of  ($\mathcal{P}_{\infty}$), then there exist $\bar\epsilon_{N}\geq 0$ where $\bar\epsilon_{N} \to 0$ as $N$ goes to infinity and
\begin{flalign}
&\int P_{\pi,N}^{*}(d\underline{\gamma})\mu^{N}(d\omega_{0},d\underline{y})c^{N}(\underline{\gamma},\underline{y}, \omega_{0})\leq\inf\limits_{P_{\pi}^{N} \in L_{\text{D}}^{N}}\int P_{\pi}^{N}(d\underline{\gamma})\mu^{N}(d\omega_{0},d\underline{y}) c^{N}(\underline{\gamma},\underline{y}, \omega_{0})+\epsilon_{N}+\bar{\epsilon}_{N}\label{eq:approxl56},
\end{flalign}
where $P_{\pi,N}^{*}$ is the restriction of $P_{\pi}^{*}$ to the first $N$ components.
\end{itemize}
\end{theorem}
{The main idea for establishing Part(i) is to use Lemma \ref{lem:exc} and Lemma \ref{lem:findef} to provide an approximation of optimal expected cost by restricting the search for randomized policies to those that are restrictions  of randomized policies in $L_{\text{EX}}$ to the $N$ first components. We note that since the set of policies $L_{\text{D}}^{N}$ is not a convex subset of the set of randomized policies $L^{N}$, \eqref{eq:approx} does not immediately imply \eqref{eq:approx1} using Lemma \ref{lem:exc} but the result can be established using an extreme point argument and since policies in $L_{\text{D}}^{N}$ are optimal among all randomized policies $L^{N}_{\text{PR}}$ for $N$-DM teams thanks to Blackwell's irrelevant information theorem \cite{Blackwell2}. Part(ii) follows from Part(i) and Theorem \ref{the:2}, using the fact that a randomized optimal policy $P_{\pi}^{*}\in L_{\text{PR,SYM}}$ of  ($\mathcal{P}_{\infty}$) provides an approximation for the optimal expected cost when the search for randomized optimal policy for $N$-DM teams is restricted to those in $L_{\text{PR,SYM}}^{N}$.}

Similarly, we present approximation results of optimal policies for symmetric dynamic $N$-DM team problems.
\begin{theorem}\label{the:3d}
Consider a dynamic team problem ($\mathcal{P}^{N}_{T}$) (see \eqref{eq:mfcost}). Let Assumption \ref{assump:ind2}, Assumption \ref{assump:c}, and Assumption \ref{assump:2} hold. 
If the cost function is bounded, then
\begin{itemize}[wide]
\item [(i)]
\begin{flalign}
\inf\limits_{P_{\pi}^{N} \in L_{\text{PR, SYM}}^{N}}&\int P_{\pi}^{N}(d\underline{\gamma})\mu^{N}(d\omega_{0},d\underline{\zeta})c^{N}( \underline{\zeta}, \underline{\gamma},\underline{y}, \omega_{0})\nu^{N}(d\underline{y}|\underline{\zeta}, \underline{\gamma},\omega_{0})\nonumber\\
&\leq \inf\limits_{P_{\pi}^{N} \in L_{\text{CO}}^{N}}\int P_{\pi}^{N}(d\underline{\gamma})\mu^{N}(d\omega_{0},d\underline{\zeta})c^{N}( \underline{\zeta}, \underline{\gamma},\underline{y}, \omega_{0})\nu^{N}(d\underline{y}|\underline{\zeta}, \underline{\gamma},\omega_{0})+ \epsilon_{N}\label{eq:approxd},
\end{flalign}
and
\begin{flalign}
\inf\limits_{P_{\pi}^{N} \in L_{\text{PR, SYM}}^{N}}& \int P_{\pi}^{N}(d\underline{\gamma})\mu^{N}(d\omega_{0},d\underline{\zeta})c^{N}( \underline{\zeta}, \underline{\gamma},\underline{y}, \omega_{0})\nu^{N}(d\underline{y}|\underline{\zeta}, \underline{\gamma},\omega_{0})\nonumber\\
&\leq \inf\limits_{P_{\pi}^{N} \in L_{\text{D}}^{N}}\int P_{\pi}^{N}(d\underline{\gamma})\mu^{N}(d\omega_{0},d\underline{\zeta})c^{N}( \underline{\zeta}, \underline{\gamma},\underline{y}, \omega_{0})\nu^{N}(d\underline{y}|\underline{\zeta}, \underline{\gamma},\omega_{0})+\epsilon_{N}\label{eq:approx1d},
\end{flalign}
where $\epsilon_{N} \to 0$ as $N$ goes to infinity.
\item [(ii)] If $P_{\pi}^{*}\in L_{\text{PR,SYM}}$ is a randomized optimal policy of  ($\mathcal{P}^{\infty}_{T}$), then there exist $\bar\epsilon_{N}\geq 0$ where $\bar\epsilon_{N} \to 0$ as $N$ goes to infinity and
\begin{flalign}
&\int  P_{\pi, N}^{*}(d\underline{\gamma})\mu^{N}(d\omega_{0},d\underline{\zeta})c^{N}( \underline{\zeta}, \underline{\gamma},\underline{y}, \omega_{0})\nu^{N}(d\underline{y}|\underline{\zeta}, \underline{\gamma},\omega_{0})\nonumber\\
&\leq\inf\limits_{P_{\pi}^{N} \in L_{\text{D}}^{N}}\int P_{\pi}^{N}(d\underline{\gamma})\mu^{N}(d\omega_{0},d\underline{\zeta})c^{N}( \underline{\zeta}, \underline{\gamma},\underline{y}, \omega_{0})\nu^{N}(d\underline{y}|\underline{\zeta}, \underline{\gamma},\omega_{0})+\epsilon_{N}+\bar{\epsilon}_{N}\nonumber,
\end{flalign}
where $P_{\pi,N}^{*}$ is the restriction of $P_{\pi}^{*}$ to the first $N$ components.
\end{itemize}
\end{theorem}
\begin{proof}
Proof follows from a similar steps as the proof of Theorem \ref{the:3} using the results of Lemma \ref{lem:findefdyn2} and Theorem \ref{the:exmftdy2}.
\end{proof}
\begin{appendix}
\section{Connection between $L_{\text{CO}}^{N}$ and $L^{N}$ in Section \ref{sec:strategic}}
In following theorem, we show that sets of randomized policies $L_{\text{CO}}^{N}$ and $L^{N}$ are identical.

\begin{theorem}\label{the:lcoln}
The set of randomized policies $L^{N}$ is identical to the set of randomized policies  $L_{\text{CO}}^{N}$.
\end{theorem}
\begin{proof}
Clearly, we have $L_{\text{CO}}^{N}\subseteq L^{N}$. In the following, we show $L^{N}\subseteq L_{\text{CO}}^{N}$.  Following from \cite{BorkarRealization}, for each $i=1,\dots, N$, the set of marginals of randomized policies in $L^{N}$ to each coordinate $\Gamma^{i}$ is a convex combination of its extreme points which is in the set of Delta-dirac measures of elements in $\Gamma^{i}$. Hence, for the set of extreme points of the convex set $L^{N}$ (denoted by $\text{Extreme}(L^{N})$), we have 
\begin{flalign*}
&\text{Extreme}(L^{N})\\
&\subseteq\bigg\{P_{\pi} \in L^{N}\bigg{|}\text{for all}~A_{i} \in \mathcal{B}(\Gamma^{i}): P_{\pi}(\gamma^{1} \in A_{1},\dots, \gamma^{N} \in A_{N})=\prod_{i=1}^{N}1_{\{\tilde{\gamma}^{i}\in A_{i}\}},~ \text{for}~\tilde{\gamma}^{i}\in \Gamma^{i}\bigg\}
\end{flalign*}
 Hence,  $\text{Extreme}(L^{N})\subset L^{N}_{\text{CO}}$ and since both sets $L^{N}$ and $L^{N}_{\text{CO}}$ are convex, we have $L^{N}\subseteq L_{\text{CO}}^{N}$, and this completes the proof.
\end{proof}

\section{Proofs from Section \ref{sec:section4}}\label{appA}
\subsection{Proof of Lemma \ref{lem:exc}}
For any permutation $\sigma \in S_{N}$, we define a randomized policy $P_{\pi}^{\sigma} \in \bar{L}^{N}$ as a permutation, $\sigma$, of arguments of a randomized policy $P_{\pi} \in \bar{L}^{N}$, i.e., for $A^{i} \in \mathcal{B}(\Gamma^{i})$
\begin{flalign*}
&P_{\pi}^{\sigma}(\gamma^{1}\in A^{1},\dots, \gamma^{2}\in A^{N}):=P_{\pi}(\gamma^{\sigma(1)}\in A^{1},\dots, \gamma^{\sigma(N)}\in A^{N}).
\end{flalign*}
 We have
\begin{flalign}
&\int P^{\sigma}_{\pi}(d\underline{\gamma})\mu^{N}(d\omega_{0},d\underline{y})c^{N}(\underline{\gamma},\underline{y}, \omega_{0})\nonumber\\
&=\int c(\omega_{0}, u^{1},\dots, u^{N})\prod_{k=1}^{N}\gamma^{k}(du^k|y^k)\nonumber\\&\times\tilde{\mu}^{N}(dy^{1},\dots,dy^{N}|\omega_{0}) P^{\sigma}_{\pi}(d\gamma^{1},\dots,d\gamma^{N})\mathbb{P}_{0}(d\omega_{0})\nonumber\\
&=\int c(\omega_{0}, u^{1},\dots, u^{N}))\prod_{k=1}^{N}\gamma^{k}(du^k|y^k)\label{eq:4.4.4.5}\\
&\times\tilde{\mu}^{N}(dy^{1},\dots,dy^{N}|\omega_{0}) P_{\pi}(d\gamma^{\sigma(1)},\dots,d\gamma^{\sigma(N)})\mathbb{P}_{0}(d\omega_{0})\nonumber\\
&=\int c(\omega_{0}, u^{\sigma(1)},\dots, u^{\sigma(N)}))\prod_{k=1}^{N}\gamma^{\sigma(k)}(du^{\sigma(k)}|y^{\sigma(k)})\label{eq:4.4.4.4}\\
&\times\tilde{\mu}^{N}(dy^{\sigma(1)},\dots,dy^{\sigma(N)}|\omega_{0}) P_{\pi}(d\gamma^{1},\dots,d\gamma^{N})\mathbb{P}_{0}(d\omega_{0})\nonumber\\
&=\int c(\omega_{0}, u^{1},\dots, u^{N})\prod_{k=1}^{N}\gamma^{k}(du^k|y^k)\label{eq:4.4.4.6}\\&\times\tilde{\mu}^{N}(dy^{1},\dots,dy^{N}|\omega_{0})P_{\pi}(d\gamma^{1},\dots,d\gamma^{N})\mathbb{P}_{0}(d\omega_{0})\nonumber\\
&=\int P_{\pi}(d\underline{\gamma})\mu^{N}(d\omega_{0},d\underline{y})c^{N}(\underline{\gamma},\underline{y}, \omega_{0})\nonumber,
\end{flalign}
 where $\tilde{\mu}^{N}$ is the joint conditional distribution of observations $(y^{1},\dots,y^{N})$ given $\omega_{0}$, and \eqref{eq:4.4.4.5} follows from the definition of $P^{\sigma}_{\pi}$ and \eqref{eq:4.4.4.4} follows from relabeling $u^{\sigma(i)}, y^{\sigma(i)}$ with $u^{i}, y^{i}$ for all $i=1,\dots,N$. Equality \eqref{eq:4.4.4.6} follows from Assumption \ref{assump:exccost} and Assumption \ref{assump:exobs}. 
 
 Let $\epsilon\geq 0$, and consider a randomized policy $P^{*}_{\pi,\epsilon} \in \bar{L}^{N}$ such that  
 \begin{flalign*}
& \int P^{*}_{\pi,\epsilon}(d\underline{\gamma})\mu^{N}(d\omega_{0},d\underline{y})c^{N}(\underline{\gamma},\underline{y}, \omega_{0})\leq \inf\limits_{P_{\pi} \in \bar{L}^{N}}\int P_{\pi}(d\underline{\gamma})\mu^{N}(d\omega_{0},d\underline{y})c^{N}(\underline{\gamma},\underline{y}, \omega_{0})+\epsilon.
\end{flalign*}
Consider $\tilde{P}_{\pi,\epsilon}$ as a convex combination of all possible permutations of $P^{*}_{\pi,\epsilon}$ by averaging them. Since $\bar{L}^{N}$ is convex, we have $\tilde{P}_{\pi,\epsilon} \in \bar{L}^{N}$. Also, we have  $\tilde{P}_{\pi,\epsilon} \in L_{\text{EX}}^{N}$, and for any permutation $\sigma \in S_{N}$, we have
\begin{flalign}
\tilde{P}_{\pi,\epsilon}(d\gamma^{1},\dots,d\gamma^{N})&:=\sum_{\sigma\in S_{N}}\frac{1}{|S_{N}|}P^{*,\sigma}_{\pi,\epsilon}(d\gamma^{1},\dots,d\gamma^{N})\nonumber\\
&\:=\tilde{P}_{\pi,\epsilon}^{\sigma}(d\gamma^{1},\dots,d\gamma^{N})\nonumber,
\end{flalign}
where $|S_{N}|$ denotes the cardinality of the set $S_{N}$, and the second equality follows from the fact that the sum is over all permutation $\sigma$ by taking average of them. Therefore, a randomized policy $\tilde{P}_{\pi,\epsilon}$ is in  $\bar{L}^{N}\cap L_{\text{EX}}^{N}$.
 We have,
 \begin{flalign*}
\int \tilde{P}_{\pi,\epsilon}(d\underline{\gamma})\mu^{N}(d\omega_{0},d\underline{y})c^{N}(\underline{\gamma},\underline{y}, \omega_{0})&:=\int (\sum_{\sigma\in S_{N}}\alpha_{\sigma}P^{*,\sigma}_{\pi,\epsilon})(d\underline{\gamma})\mu^{N}(d\omega_{0},d\underline{y})c^{N}(\underline{\gamma},\underline{y}, \omega_{0})\\
&= \sum_{\sigma\in S_{N}}\alpha_{\sigma}\int {P}_{\pi,\epsilon}^{*,\sigma}(d\underline{\gamma})\mu^{N}(d\omega_{0},d\underline{y})c^{N}(\underline{\gamma},\underline{y}, \omega_{0})\\
&=\sum_{\sigma\in S_{N}}\alpha_{\sigma}\int {P}_{\pi,\epsilon}^{*}(d\underline{\gamma})\mu^{N}(d\omega_{0},d\underline{y})c^{N}(\underline{\gamma},\underline{y}, \omega_{0})\\
&\leq\inf\limits_{P_{\pi} \in \bar{L}^{N}}\int P_{\pi}(d\underline{\gamma})\mu^{N}(d\omega_{0},d\underline{y})c^{N}(\underline{\gamma},\underline{y}, \omega_{0})+\epsilon,
 \end{flalign*}
 where the second equality is true since the map $P_{\pi} \to \int P_{\pi}(d\underline{\gamma})\mu^{N}(d\omega_{0},d\underline{y})c^{N}(\underline{\gamma},\underline{y}, \omega_{0})$ is linear and the third equality follows from \eqref{eq:4.4.4.6}.
 Since $\tilde{P}_{N,\epsilon} \in \bar{L}^{N} \cap L_{\text{EX}}^{N}$, we have
\begin{equation*}
\int \tilde{P}_{\pi,\epsilon}(d\underline{\gamma})\mu^{N}(d\omega_{0},d\underline{y})c^{N}(\underline{\gamma},\underline{y}, \omega_{0}) \geq \inf\limits_{P_{\pi} \in \bar{L}^{N}\cap L_{\text{EX}}^{N}}\int P_{\pi}(d\underline{\gamma})\mu^{N}(d\omega_{0},d\underline{y})c^{N}(\underline{\gamma},\underline{y}, \omega_{0}).
\end{equation*}
 Hence, for any $\epsilon \geq 0$, we have 
 \begin{equation*}
\scalemath{0.95}{ \inf\limits_{P_{\pi} \in \bar{L}^{N}\cap L_{\text{EX}}^{N}}\int P_{\pi}(d\underline{\gamma})\mu^{N}(d\omega_{0},d\underline{y})c^{N}(\underline{\gamma},\underline{y}, \omega_{0}) \leq \inf\limits_{P_{\pi} \in \bar{L}^{N}}\int P_{\pi}(d\underline{\gamma})\mu^{N}(d\omega_{0},d\underline{y})c^{N}(\underline{\gamma},\underline{y}, \omega_{0})+\epsilon}.
 \end{equation*}
 Since $\epsilon$ is arbitrary, this completes the proof.

\subsection{Proof of Lemma  \ref{lem:findef}}
To prove Lemma  \ref{lem:findef}, we use two following results by Diaconis and Friedman \cite[Theorem 13]{diaconis1980finite} and Aldous \cite[Proposition 7.20]{aldous2006ecole} (see also \cite{kallenberg1973canonical} for more general results) which we recall for reader's convenience:
\begin{theorem}\cite[Theorem 13]{diaconis1980finite}\label{the:dia}
Let $Y=(Y_{1},\dots,Y_{n})$ be an $n$-exchangeable and $Z=(Z_{1},Z_{2},\dots)$ be an infinitely-exchangeable sequence of random variables with $\mathcal{L}(Z_{1},\dots,Z_{k})=\mathcal{L}(Y_{I_1},\dots,Y_{I_k})$ for all $k\geq 1$ where the indices $(I_{1}, I_{2}, \dots)$ are i.i.d. random variables with the uniform distribution on the set $\{1,\dots,n\}$. Then,  for all $m=1,\dots, n$,
\begin{flalign}
\bigg|\bigg|\mathcal{L}(Y_{1},\dots,Y_{m})-\mathcal{L}(Z_{1},\dots,Z_{m})\bigg|\bigg|_{TV}\leq \frac{m(m-1)}{2n},
\end{flalign}
where $\mathcal{L}(\cdot)$ denotes the law of random variables and $||\cdot||_{TV}$ is the total variation norm.
\end{theorem}
\begin{theorem}\cite[Proposition 7.20]{aldous2006ecole}\label{the:ald}
Let $X:=(X_{1},X_{2},\dots)$ be an infinitely-exchangeable sequence of random variables taking values in a Polish space $\mathbb{X}$ and directed by a random measure $\alpha$ (i.e., $\alpha$ is a $\mathcal{P}(\mathbb{X})$-valued random variable and $Pr(X\in A)=\int_{\mathcal{P}(\mathbb{X})} \prod_{i=1}^{\infty}\xi(A^{i})Pr(\alpha\in d\xi)$ where $A^{i}\in \mathcal{B}(\mathbb{X})$ and $(A=A^{1}\times A^{2}\times \dots)$, see \cite[Definition 2.6]{aldous2006ecole}). Suppose that either for each $n$
\begin{itemize}
\item[(1)] $X^{(n)}=(X_{1}^{(n)},X_{2}^{(n)},\dots)$  is infinitely-exchangeable directed by $\alpha_{n}$, or
\item[(2)]  $X^{(n)}=(X_{1}^{(n)},\dots, X_{n}^{(n)})$  is $n$-exchangeable with empirical measure $\alpha_{n}$.
\end{itemize}
Then, $X^{(n)}$ converges in distribution to $X$ $\big(X^{(n)}\xrightarrow[n \to \infty]{\text{d}}X\big)$  if and only if $\alpha_{n}\xrightarrow[n \to \infty]{\text{d}}\alpha$.
\end{theorem}
We note that by convergence in distribution to an infinite exchangeable sequence, we mean the following: $X^{(n)} \xrightarrow[n \to \infty]{\text{d}} X$ if and only if $(X_{1}^{(n)},\dots, X_{m}^{(n)})\xrightarrow[n \to \infty]{\text{d}}(X_{1},\dots, X_{m})$ for each $m\geq 1$ \cite[page 55]{aldous2006ecole}.

Using the above theorems, we now complete the Proof of Lemma  \ref{lem:findef}. {Since the action space $\mathbb{U}$ is compact and observations are i.i.d. with a fixed marginal (under Assumption \ref{assump:ind}, via a change of measure argument observations can be viewed to be independent of $\omega_{0}$), the set of probability measures $L^{N}$ is tight. Furthermore, by \cite[Theorem 5.1]{yuksel2018general}, $L^{N}$ is closed under the topology of weak convergence and hence $L^{N}$ is compact. Using the argument in \cite[Theorem 5.1]{yuksel2018general} under Assumption \ref{assump:cont}, the expected cost function is lower semicontinuous in policies $P_{N} \in L^{N}$. Hence, there exists an optimal policy for {($\mathcal{P}_{N}$)}, and by Lemma \ref{lem:exc}, this optimal policy can be assumed to be in $L_{\text{EX}}^{N}$.} Consider a sequence of $N$-exchangeable randomized policies $\{P_{\pi}^{*,N}\}_{N}$, where for every $N\geq 1$, $P_{\pi}^{*,N} \in L_{\text{EX}}^{N}$ and
\begin{equation}\label{eq:32.1t}
\scalemath{0.95}{\int P_{\pi}^{*, N}(d\underline{\gamma})\mu^{N}(d\omega_{0},d\underline{y})c^{N}(\underline{\gamma},\underline{y}, \omega_{0})=\inf\limits_{P^{N}_{\pi} \in L_{\text{EX}}^{N}}\int P^{N}_{\pi}(d\underline{\gamma})\mu^{N}(d\omega_{0},d\underline{y})c^{N}(\underline{\gamma},\underline{y}, \omega_{0})}.
\end{equation}
In the following, we show \eqref{eq:ert} in two steps. In the first step, for every $N$, we use the construction in Theorem \ref{the:dia} to construct an infinitely-exchangeable randomized policy $P^{*,\infty}_{\pi, N} \in L_{\text{EX}}$ using the $N$-exchangeable randomized policy $P_{\pi}^{*, N} \in  L^{N}_{\text{EX}}$ by considering the indices as a sequence of i.i.d. random variables with uniform distribution on the set $\{1,\dots,N\}$, and then, we show that there exists a weakly convergent subsequence  of joint measures on the first coordinate, observations, and the average of induced actions of randomized policies  $P^{*,\infty}_{\pi, N} \in L_{\text{EX}}$. Then, we show that the expected cost function induced by the $N$-exchangeable randomized policy $P_{\pi}^{*, N} \in  L^{N}_{\text{EX}}$  converges through a subsequence to a limit induced by an infinitely-exchangeable randomized policy $P^{*,\infty}_{\pi, N}$.
 
\begin{itemize}[wide]
\item[\textbf{(Step 1):}]
Let $(I_{1},I_{2},\dots)$ be i.i.d. random variables with uniform distribution on the set $\{1,\dots,N\}$. For a fixed $N$ and for any $N$-exchangeable randomized policy $P_{\pi}^{*, N} \in  L^{N}_{\text{EX}}$, we construct  an infinitely-exchangeable randomized policy $P^{*,\infty}_{\pi, N} \in L_{\text{EX}}$ as follows: for every $N$ and $m$ and for all $A^{i} \in \mathcal{B}(\Gamma^{i})$
\begin{flalign*}
&P^{*,\infty}_{\pi, N}(\gamma^{1}\in A^{1},\dots, \gamma^{m}\in A^{m}):=P_{\pi}^{*, N}(\gamma^{I_{1}}\in A^{1},\dots, \gamma^{I_{m}}\in A^{m}).
\end{flalign*}
where $P^{*,\infty}_{\pi, N}$ is the restriction of $P^{*,\infty}_{\pi, P_{\pi}^{N}}\in  L_{\text{EX}}$ to the first $N$ components.
We note that $P^{*,\infty}_{\pi, N} \in L_{\text{EX}}$ because we use i.i.d. sequence $(I_{1},I_{2},\dots)$ for indexing probability measures on the space of policies, hence, for every fixed $N$ and $N$-exchangeable randomized policy $P_{\pi}^{*, N}$, a randomized policy $P^{*,\infty}_{\pi, N}$ is i.i.d through DMs and hence it is infinitely-exchangeable. 

Let $u^{*, i}_{N}$ be the control action induced by $\gamma^{i}_{N}$ where random variables $(\gamma^{1}_{N},\dots,\gamma^{N}_{N})$ are determined by $N$-exchangeable randomized policy $P_{\pi}^{*, N}\in L_{\text{EX}}^{N}$. Let $u^{*, i}_{\infty, N}$ be the control action induced by $\gamma^{i}_{N, \infty}$ where random variables $(\gamma^{1}_{N, \infty},\dots,\gamma^{N}_{N,\infty})$ are determined by infinitely-exchangeable randomized policy $P_{\pi, N}^{*, \infty}\in L_{\text{EX}}$. Since under the reduction (Assumption \ref{assump:ind}), observations are i.i.d. and also independent of $\omega_{0}$, following from Theorem \ref{the:dia}, we have for every $m\geq 1$
\begin{flalign}
&\scalemath{0.94}{\bigg|\bigg|\mathcal{L}(\gamma^{1}_{N},\dots,\gamma^{m}_{N},y^{1},\dots,y^{m})-\mathcal{L}(\gamma^{1}_{N, \infty},\dots,\gamma^{m}_{N,\infty},y^{1},\dots,y^{m})\bigg|\bigg|}\nonumber\\
&\scalemath{0.94}{=\bigg|\bigg|\mathcal{L}(\gamma^{1}_{N},\dots,\gamma^{m}_{N})\prod_{i=1}^{m}\mathcal{L}(y^{i})-\mathcal{L}(\gamma^{1}_{N, \infty},\dots,\gamma^{m}_{N,\infty})\prod_{i=1}^{m}\mathcal{L}(y^{i})\bigg|\bigg|_{TV}\xrightarrow[N \to \infty]{}0}\label{eq:poA.7},
\end{flalign}
where \eqref{eq:poA.7} follows from the fact that $(\gamma^{1}_{N},\dots,\gamma^{N}_{N})$ and $(\gamma^{1}_{N, \infty},\dots,\gamma^{N}_{N, \infty})$ are random variables with joint probability measures $P_{\pi}^{*, N} \in L_{\text{EX}}^{N}$ and $P_{\pi, N}^{*, \infty} \in L_{\text{EX}}\big|_{N}$, respectively. 

Since $\mathbb{U}$ is compact, the marginal of probability measures on $\mathbb{U}$ is tight. Since the probability measure on $\mathbb{Y}$ is fixed, the marginal on $\mathbb{Y}$ is also tight. Since marginals are tight, then the collection of all measures on $(\mathbb{U} \times \mathbb{Y})$ with these tight marginals is also tight (see e.g., \cite[Proof of Theorem 2.4]{yukselSICON2017}), and hence, the set $\Gamma^{i}$ is tight for each $i\in \mathbb{N}$. Hence, $\{\mathcal{L}(\gamma^{i}_{\infty, N})\}_{N}$ is tight for each DM and by exchangeablity $\mathcal{L}(\gamma^{i}_{\infty, N})=\mathcal{L}(\gamma^{1}_{\infty, N})$. Hence, we can find a subsequence such that $\mathcal{L}(\gamma^{i}_{\infty, l})\xrightarrow[l\to \infty]{} \mathcal{L}(\gamma^{i}_{\infty})$ for all $i\in \mathbb{N}$. Since marginals of $\{\mathcal{L}(\gamma^{1}_{\infty, l}, \dots, \gamma^{m}_{\infty, l})\}_{l}$ are tight, for each $m\geq 1$, there exists a further subsequence 
\begin{flalign*}
\mathcal{L}(\gamma^{1}_{\infty, n}, \dots, \gamma^{m}_{\infty, n}) \xrightarrow[n \to \infty]{} \mathcal{L}(\gamma_{\infty}^{1}, \dots, \gamma_{\infty}^{m}),
\end{flalign*}
where $(\gamma_{\infty}^{1}, \gamma_{\infty}^{2}, \dots)$ is infinitely-exchangeable and induced by an infinitely-exchangeable randomized policy $P_{\pi}^{*, \infty} \in L_{\text{EX}}$ since the set of infinitely-exchangeable randomized policies is closed under the weak-convergence topology,  where by weak convergence for an infinite sequence, we mean weak convergence of finite restrictions. That is because, if $P_{\pi}^{\sigma, *, \infty}$ is the limit in the weak convergence topology of the sequence randomized policies $\{P_{\pi, n}^{\sigma, *, \infty}\}_{n}$ as $n \to \infty$, where for $A^{i} \in \mathcal{B}(\Gamma^{i})$ and for all $N\in\mathbb{N}$ and all finite permutations $\sigma \in S_{N}$
\begin{flalign*}
&P_{\pi, n}^{\sigma, *, \infty}(\gamma^{1}\in A^{1},\gamma^{2}\in A^{2},\dots):=P_{\pi, n}^{*, \infty}(\gamma^{\sigma(1)}\in A^{1},\gamma^{\sigma(2)}\in A^{2},\dots).
\end{flalign*}
Then, following from exchangeability, since sequences $\{P_{\pi, n}^{*, \infty}\}_{n}$ and $\{P_{\pi, n}^{\sigma, *, \infty}\}_{n}$ are identical, the limit in the weak convergence topology of both randomized policies $P_{\pi}^{*, \infty}$ and $P_{\pi}^{\sigma, *, \infty}$ are also identical, and hence, the limit $P_{\pi}^{*, \infty}$  is infinitely-exchangeable, $P_{\pi}^{*, \infty}\in L_{\text{EX}}$. Hence, following from \eqref{eq:poA.7}, for each $m\geq 1$
\begin{flalign*}
\mathcal{L}(\gamma^{1}_{n}, \dots, \gamma^{m}_{n}) \xrightarrow[n \to \infty]{} \mathcal{L}(\gamma_{\infty}^{1}, \dots, \gamma_{\infty}^{m}).
\end{flalign*}
By construction of random variables $u^{*, i}_{n}$ and $u^{*, i}_{\infty}$ and since random variables $\gamma^{i}_{n}$s are independent of $y^{i}$s, we have  for each $m\geq 1$
\begin{flalign*}
(u^{*, 1}_{n}, \dots, u^{*, m}_{n}) \xrightarrow[n \to \infty]{\text{d}} (u_{\infty}^{1}, \dots, u_{\infty}^{m}),
\end{flalign*}
 where $(u_{\infty}^{1}, u_{\infty}^{2}, \dots)$ is induced by an infinitely-exchangeable policy $P_{\pi}^{*, \infty} \in L_{\text{EX}}$.  
Following from Theorem \ref{the:ald}, we have for all $A \in \mathcal{U}$ and $\mathbb{P}$-almost surely
\begin{flalign}
F_{n}(A):=F_{n}^{\omega}(A):=\frac{1}{n}\sum_{i=1}^{n}\delta_{u^{*,i}_{n}(\omega)}(A)\xrightarrow[n \to \infty]{\text{d}}\alpha^{\omega}(A)\label{eq:emp},
\end{flalign}
where $\omega$ denotes the sample path dependency and $\alpha$ is the directing measure of an infinitely-exchangeable random variables $(u_{\infty}^{1},u_{\infty}^{2},\dots)$ (that is $\alpha(\omega, A)=Pr(u_{\infty}^{*, i}\in A|H)$ almost surely for all $A \in \mathcal{U}$ where $H$ is the $\sigma$-field generated by $\mathcal{P}(\mathbb{U})$-valued random variable $\alpha$ \cite{aldous2006ecole}). {Following from \eqref{eq:emp}, since the action space $\mathbb{U}$ is compact, we have $\mathbb{P}$-almost surely
\begin{flalign}
 \bar{F}_{n}:=\bar{F}_{n}^{\omega}:=\frac{1}{n}\sum_{i=1}^{n}u^{*,i}_{n}(\omega)=\int_{\mathbb{U}}uF_{n}(du) \xrightarrow[n \to \infty]{\text{d}} \bar{F}:=\int_{\mathbb{U}}u\alpha^{\omega}(du).
\end{flalign}}

Define $\tilde{P}^{*, n}$ as the joint probability measure of $(u^{*,1}_{n}, \bar{F}_n, \underline{y})$ where marginals on $\underline{y}:=(y^{1},y^{2},\dots)$ are fixed to be $\prod_{i=1}^{\infty}Q(dy^{i})$. Since marginals on $(u^{*,1}_{n}, \bar{F}_n)$ are tight and marginals on $\underline{y}$ are fixed, $\{\tilde{P}^{*, n}\}_{n}$ is tight. Hence, there exists a subsubsequence $\{\tilde{P}^{*, {k}}\}_{{k}}$ converges weakly to $\tilde{P}^{*}$ as $k$ goes to infinity. This implies that marginals $\{\tilde{P}^{*, {k}}\}_{{k}}$ on $(u^{*,1}_{{k}}, \bar{F}_{{k}})$ converges to the marginals of  $\tilde{P}^{*}$ on $(u^{*,1}, \bar{F})$, hence, $\tilde{P}^{*}$ is induced by $(u_{\infty}^{1},u_{\infty}^{2},\dots)$ which is infinitely-exchangeable and is induced by an infinitely-exchangeable randomized policy in $L_{\text{EX}}$.
 \item[\textbf{(Step 2):}]
We have
\begin{flalign}
&\scalemath{0.94}{ \limsup\limits_{N \to \infty}\int P_{\pi}^{*, N}(d\underline{\gamma})\mu^{N}(d\omega_{0},d\underline{y})c^{N}(\underline{\gamma},\underline{y}, \omega_{0})}\nonumber\\
&\scalemath{0.94}{ =\limsup\limits_{N \to \infty}\frac{1}{N}\sum_{i=1}^{N}\int c(\omega_{0}, u^{i}, \frac{1}{N}\sum_{p=1}^{N}u^{p})\prod_{k=1}^{N}\gamma^{k}(du^{k}|y^{k})P_{\pi}^{*, N}(d{\gamma}^{1},\dots,d\gamma^{N})}\nonumber\\
&\scalemath{0.94}{\times\prod_{i=1}^{N}\hat{\mu}(dy^{i}|\omega_{0})\mathbb{P}_{0}(d\omega_{0})}\nonumber\\
&\scalemath{0.94}{ =\limsup\limits_{N \to \infty} \frac{1}{N}\sum_{i=1}^{N}\int c(\omega_{0}, u^{i}, \frac{1}{N}\sum_{p=1}^{N}u^{p})\prod_{k=1}^{N}\gamma^{k}(du^{k}|y^{k})P_{\pi}^{*, N}(d{\gamma}^{1},\dots,d\gamma^{N})}\label{eq:A.9a}\\
&\scalemath{0.94}{\times\prod_{i=1}^{N}f(\omega_{0},y^{i})Q(dy^{i})\mathbb{P}_{0}(d\omega_{0})}\nonumber\\
&\scalemath{0.94}{ =\limsup\limits_{N \to \infty}\int \int_{\prod_{i=N+1}^{\infty}\mathbb{Y}}c(\omega_{0}, u^{1}, \bar{F}_{N})\tilde{P}^{*, N}(d{u}^{1}, d\bar{F}_{N}, d\underline{y})\prod_{i=1}^{\infty}f(\omega_{0},y^{i})\mathbb{P}_{0}(d\omega_{0})}\label{eq:A.9}\\
&\scalemath{0.94}{ \geq \lim\limits_{{k} \to \infty} \int \int_{\prod_{i={k}+1}^{\infty}\mathbb{Y}} c(\omega_{0}, u^{1}, \bar{F}_{{k}})\tilde{P}^{*, {k}}(d{u}^{1}, d\bar{F}_{k}, d\underline{y})\prod_{i=1}^{\infty}f(\omega_{0},y^{i})\mathbb{P}_{0}(d\omega_{0})}\label{eq:A.10}\\
&\scalemath{0.94}{ = \int c(\omega_{0}, u^{1}, \bar{F})\tilde{P}^{*}(d{u}^{1}, d\bar{F}, , d\underline{y})\prod_{i=1}^{\infty}f(\omega_{0},y^{i})\mathbb{P}_{0}(d\omega_{0})}\label{eq:A.11}\\
&\scalemath{0.94}{\geq \limsup\limits_{N \to \infty}\inf\limits_{P_{\pi} \in L_{\text{EX}}}\int P_{\pi, N}(d\underline{\gamma})\mu^{N}(d\omega_{0},d\underline{y})c^{N}(\underline{\gamma},\underline{y}, \omega_{0})}\label{eq:A.12}.
\end{flalign}
{where $\hat{\mu}$ is the conditional distribution of each observation $y^{i}$ given $\omega_{0}$}, and \eqref{eq:A.9a} follows from Assumption \ref{assump:oc}(i), hence,  under Assumption \ref{assump:ind}, in the new (equivalent) expected cost function, observations are i.i.d. and independent of $\omega_{0}$. \eqref{eq:A.9} follows from integrating over the set $\prod_{i=N+1}^{\infty}\mathbb{Y}$ and since $(u_{N}^{*, 1},\dots, u_{N}^{*, N})$ is $N$-exchangeable. Inequality \eqref{eq:A.10} follows from the assumption that the cost function is bounded and limsup is the greatest subsequence limit of a bounded sequence where $k$ is the index of the subsequence considered in (Step 1). Equality \eqref{eq:A.11} follows from the dominated convergence theorem and following from Assumption \ref{assump:cont} and since by (Step 1) $\{\tilde{P}^{*, {k}}\}_{{k}}$ converges weakly to $\tilde{P}^{*}$ as ${k}$ goes to infinity. Inequality \eqref{eq:A.12} follows from the fact that $\tilde{P}^{*}$ is the joint measure with the first coordinate $(u_{\infty}^{1},u_{\infty}^{2},\dots)$ which is infinitely-exchangeable and it is induced by an infinitely-exchangeable randomized policy in $L_{\text{EX}}$. The above inequalities become equalities since the opposite direction is true as well (that is because $L_{\text{EX}}\big|_{N}\subset L_{\text{EX}}^{N}$) and this completes the proof.
\end{itemize}

\subsection{Proof of Theorem \ref{the:2}}
We complete the proof in four steps.
\begin{itemize}[wide]
\item[\textbf{(Step 1):}] 
Similar to the proof of Lemma  \ref{lem:findef}, using  \cite[Theorem 5.1]{yuksel2018general}, we can show that there exists a randomized optimal policy for {($\mathcal{P}_{N}$)} which belongs to the set $L^{N}$, and by Lemma \ref{lem:exc}, this randomized optimal policy can be assumed to be in the set of $N$-exchangeable randomized policies $L_{\text{EX}}^{N}$. Consider a sequence of $N$-exchangeable randomized policies $\{P_{\pi}^{*,N}\}_{N}$, where for every $N\geq 1$, $P_{\pi}^{*,N} \in L_{\text{EX}}^{N}$ and
\begin{equation*}
\int P_{\pi}^{*, N}(d\underline{\gamma})\mu^{N}(d\omega_{0},d\underline{y})c^{N}(\underline{\gamma},\underline{y}, \omega_{0})=\inf\limits_{P^{N}_{\pi} \in L_{\text{EX}}^{N}}\int P^{N}_{\pi}(d\underline{\gamma})\mu^{N}(d\omega_{0},d\underline{y})c^{N}(\underline{\gamma},\underline{y}, \omega_{0}).
\end{equation*}
 
\item[\textbf{(Step 2):}]
In this step, we show that to establish an existence result, it is sufficient to show the convergence of the expected cost induced by a randomized optimal policy in $L_{\text{PR,SYM}}^{N}$ of $N$-DM teams to the expected cost induced by a randomized policy $L_{\text{PR,SYM}}$ of mean-field teams through a subsequence as $N$ goes to infinity. We first lift the space of randomized admissible policies, and we represent any admissible randomized policy as a probability measure in $L$ (which is convex) and $L_{\text{EX}}\subset L$. We have
\begin{flalign}
&\inf\limits_{P_{\pi} \in L}\limsup\limits_{N \to \infty} \int P_{\pi, N}(d\underline{\gamma})\mu^{N}(d\omega_{0},d\underline{y})c^{N}(\underline{\gamma},\underline{y}, \omega_{0})\nonumber\\
&\geq\limsup\limits_{N \to \infty}\inf\limits_{P_{\pi}^{N} \in L^{N}}\int P_{\pi}^{N}(d\underline{\gamma})\mu^{N}(d\omega_{0},d\underline{y})c^{N}(\underline{\gamma},\underline{y}, \omega_{0})\label{eq:q4.1}\\
&=\limsup\limits_{N \to \infty}\inf\limits_{P_{\pi}^{N} \in L^{N}_{\text{EX}}}\int P_{\pi}^{N}(d\underline{\gamma})\mu^{N}(d\omega_{0},d\underline{y})c^{N}(\underline{\gamma},\underline{y}, \omega_{0})\label{eq:q4.2}\\
&\geq  \lim\limits_{M \to \infty}\limsup\limits_{N \to \infty}\inf\limits_{P_{\pi}^{N} \in L^{N}_{\text{EX}}}\int P_{\pi}^{N}(d\underline{\gamma})\mu^{N}(d\omega_{0},d\underline{y})\min{\{M, c^{N}(\underline{\gamma},\underline{y}, \omega_{0})\}}\label{eq:q4.s2}\\
&=  \lim\limits_{M \to \infty}\limsup\limits_{N \to \infty}\inf\limits_{P_{\pi} \in L_{\text{EX}}}\int P_{\pi, N}(d\underline{\gamma})\mu^{N}(d\omega_{0},d\underline{y})\min{\{M, c^{N}(\underline{\gamma},\underline{y}, \omega_{0})\}}\label{eq:q4.3}\\
&= \lim\limits_{M \to \infty}\limsup\limits_{N \to \infty}\inf\limits_{P_{\pi}^{N} \in L_{\text{CO,SYM}}^{N}}\int P_{\pi}^{N}(d\underline{\gamma})\mu^{N}(d\omega_{0},d\underline{y})\min{\{M, c^{N}(\underline{\gamma},\underline{y}, \omega_{0})\}}\label{eq:4.1}\\
&= \lim\limits_{M \to \infty}\limsup\limits_{N \to \infty}\inf\limits_{P_{\pi}^{N} \in L_{\text{PR,SYM}}^{N}}\int P_{\pi}^{N}(d\underline{\gamma})\mu^{N}(d\omega_{0},d\underline{y})\min{\{M, c^{N}(\underline{\gamma},\underline{y}, \omega_{0})\}}\label{eq:4.3}\\
&\geq \inf\limits_{P_{\pi} \in L_{\text{PR,SYM}}}\limsup\limits_{N \to \infty} \int P_{\pi, N}(d\underline{\gamma})\mu^{N}(d\omega_{0},d\underline{y})c^{N}(\underline{\gamma},\underline{y}, \omega_{0})\label{eq:4.4}\\
&\geq \inf\limits_{P_{\pi} \in L_{\text{CO,SYM}}}\limsup\limits_{N \to \infty} \int P_{\pi, N}(d\underline{\gamma})\mu^{N}(d\omega_{0},d\underline{y})c^{N}(\underline{\gamma},\underline{y}, \omega_{0})\label{eq:4.4i}\\
&\geq\inf\limits_{P_{\pi} \in L} \limsup\limits_{N \to \infty} \int P_{\pi, N}(d\underline{\gamma})\mu
^{N}(d\omega_{0},d\underline{y})c^{N}(\underline{\gamma},\underline{y}, \omega_{0})\label{eq:4.6},
\end{flalign}
where \eqref{eq:q4.1} follows from exchanging limsup with inf and the fact that the restriction to $N$-first coordinate $P_{\pi, N} \in L^{N}$ for any randomized policy $P_{\pi}\in L$, and \eqref{eq:q4.2} follows from Lemma \ref{lem:exc}. Inequality \eqref{eq:q4.s2} follows from $\min{\{M, c^{N}(\underline{\gamma},\underline{y}, \omega_{0})\}}\leq c^{N}(\underline{\gamma},\underline{y}, \omega_{0})$. Equality \eqref{eq:q4.3} follows from Lemma \ref{lem:findef} and \eqref{eq:4.1} follows from Theorem \ref{the:defin}. The set of extreme points of the convex set $L_{\text{CO,SYM}}^{N}$ is in $L_{\text{PR,SYM}}^{N}$ (that is because, $L_{\text{CO,SYM}}^{N}$ corresponds to the randomized policies with common and individual independent randomness where each DM choose an identical randomized policy), hence, \eqref{eq:4.3} is true since $L_{\text{CO,SYM}}^{N}$ is convex, and the map $\int P_{\pi}^{N}(d\underline{\gamma})\mu^{N}(d\omega_{0},d\underline{y})c^{N}(\underline{\gamma},\underline{y}, \omega_{0}):L_{\text{CO,SYM}}^{N} \to \mathbb{R}$ is linear. Inequalities \eqref{eq:4.4i} and \eqref{eq:4.6} follow from the fact that $L_{\text{PR,SYM}} \subset L_{\text{CO,SYM}}\subset L$. Hence, by \eqref{eq:4.6}, this chain of inequalities must be chain of equalities.  

In the next two steps, we justify \eqref{eq:4.4} through showing that there exists a subsequence of policies induced by symmetric/identical private randomization whose weak-limit achieves \eqref{eq:4.4}. First, we establish compactness of the set of randomized policies $L_{\text{PR,SYM}}^{N}$, and then, we show a lower semicontinuity of the induced expected cost function justifying \eqref{eq:4.4}.
\item[\textbf{(Step 3):}] Consider the set of randomized policies $L_{\text{PR,SYM}}^{N}$. For each DM, we can equivalently represent any randomized policy as a probability measure on $(\mathbb{U}\times \mathbb{Y})$, where the marginal on observations is fixed. Since the team is static, this decouples the policy spaces from the policies of the previous decision makers. Following from symmetry, we can represent each DM's privately randomized policy space as  $\{P\in \mathcal{P}(\mathbb{U}\times \mathbb{Y})| P(B)=\int_{B} \Pi(du^{i}|y^{i})Q(dy^{i})\}$ where $B \in \mathcal{B}(\mathbb{U}\times \mathbb{Y})$ and $\Pi$ is an identical randomized policy from the set of stochastic kernels from space of observations to space of actions for each DM. Since $\mathbb{U}$ is compact, the marginals on $\mathbb{U}$ are relatively compact.  Since the marginals are relatively compact, the collection of all measures with these relatively compact marginals are also relatively compact (see e.g., \cite[Proof of Theorem 2.4]{yukselSICON2017}), and hence, the randomized policy space is relatively compact. Following from symmetry, the set of individual randomized policies for each DM is closed under product topology where each coordinate converges in the weak convergence topology.
{Hence, (Step 3) implies that there exists a subsequence of (symmetric) individually randomized policies for each DM that converges weakly to the limit which is identical for each DM. In (Step 4), we show that the limit randomized policy is optimal by showing a lower semicontinuity of the induced expected cost function. 
\item[\textbf{(Step 4):}]
Define the empirical measure on actions and observations as follows:
\begin{equation*}
\Lambda_{N}(B):=\frac{1}{N}\sum_{i=1}^{N}\delta_{\beta^{i}_{N}}(B),
\end{equation*}
where for each $N$, $\beta_{N}^{i}:=(u^{i, *}_{N},y^{i})$, $B \in \mathcal{Z}:=(\mathbb{U} \times \mathbb{Y})$, $u^{i,*}_{N}$ is the action induced by the randomized policy $\Pi^{*}_{N}$ in (Step 3).

Now, we have
\begin{flalign}
&\scalemath{0.95}{\lim\limits_{M \to \infty}\limsup\limits_{N \to \infty}\inf\limits_{P_{\pi}^{N} \in L_{\text{PR,SYM}}^{N}}\int P_{\pi}^{N}(d\underline{\gamma})\mu^{N}(d\omega_{0},d\underline{y})\min{\{M, c^{N}(\underline{\gamma},\underline{y}, \omega_{0})\}}}\nonumber\\
&\scalemath{0.95}{ =\lim\limits_{M\rightarrow \infty} \limsup\limits_{N\rightarrow \infty}\int\bigg(\int\min{\left\{M, {c}\left(\omega_{0}, u,\int_{\mathbb U}u\Lambda_{N}(du \times \mathbb{Y})\right)\right\}}\Lambda_{N}(du,dy)\bigg)}\label{eq:4.32r2}\\
&\scalemath{0.95}{\times \prod_{i=1}^{\infty}P_{N}^{*, \omega_{0}}(du^{i},dy^{i})\mathbb{P}_{0}(d\omega_{0})}\nonumber\\
&\scalemath{0.95}{ \geq \lim\limits_{M\rightarrow \infty} \lim\limits_{n\rightarrow \infty}\int\bigg(\int\min{\left\{M, {c}\left(\omega_{0}, u,\int_{\mathbb U}u\Lambda_{n}(du \times \mathbb{Y})\right)\right\}}\Lambda_{n}(du,dy)\bigg)}\label{eq:4.71}\\
&\scalemath{0.95}{\times\prod_{i=1}^{\infty}P_{n}^{*, \omega_{0}}(du^{i},dy^{i})\mathbb{P}_{0}(d\omega_{0})}\nonumber\\
& \scalemath{0.95}{=\lim\limits_{M\rightarrow \infty}\int\lim\limits_{n\rightarrow \infty}\int\bigg(\int\min{\left\{M, {c}\left(\omega_{0}, u,\int_{\mathbb U}u\Lambda_{n}(du \times \mathbb{Y})\right)\right\}}\Lambda_{n}(du,dy)\bigg)} \label{eq:4.91}\\
&\scalemath{0.95}{\times \prod_{i=1}^{\infty}P_{n}^{*, \omega_{0}}(du^{i},dy^{i})\mathbb{P}_{0}(d\omega_{0})}\nonumber\\
& \scalemath{0.95}{\geq\lim\limits_{M\rightarrow \infty}\int\int\bigg(\int\min{\left\{M, {c}\left(\omega_{0}, u,\int_{\mathbb U}u\Lambda(du \times \mathbb{Y})\right)\right\}}\Lambda(du,dy)\bigg)}\label{eq:4.101}\\
&\scalemath{0.95}{\times \prod_{i=1}^{\infty}P^{*, \omega_{0}}(du^{i},dy^{i})\mathbb{P}_{0}(d\omega_{0})}\nonumber\\
& \scalemath{0.95}{= \int\bigg(\int {c}\left(\omega_{0}, u,\int_{\mathbb U}u\Lambda(du \times \mathbb{Y})\right)\Lambda(du,dy)\bigg)\prod_{i=1}^{\infty}P^{*, \omega_{0}}(du^{i},dy^{i})\mathbb{P}_{0}(d\omega_{0})}\label{eq:4.111}\\
&\scalemath{0.95}{ =\limsup\limits_{N \to \infty}\int\frac{1}{N}\sum_{i=1}^{N}{c}\bigg(\omega_{0}, u^{i},\frac{1}{N}\sum_{p=1}^{N}u^{p}\bigg) \prod_{i=1}^{N} P^{*, \omega_{0}}(du^{i},dy^{i})\mathbb{P}_{0}(d\omega_{0})}\label{eq:4.32c2}\\
&\scalemath{0.95}{\geq \inf\limits_{P_{\pi} \in L_{\text{PR,SYM}}}\limsup\limits_{N \to N} \int P_{\pi, N}(d\underline{\gamma})\mu^{N}(d\omega_{0},d\underline{y})c^{N}(\underline{\gamma},\underline{y}, \omega_{0})}\label{eq:4.35r1}
\end{flalign}
where $P_{N}^{*, \omega_{0}}(du^{i},dy^{i}):=\Pi_{N}^{*}(du^{i}|dy^{i})\hat{\mu}(dy^{i}|\omega_{0})=\Pi_{N}^{*}(du^{i}|y^{i})f(\omega_{0},y^{i})Q(dy^{i})$. Equality \eqref{eq:4.32r2} follows from the definition of the empirical measure, and Assumption \ref{assump:oc}(i), and follows from symmetry of optimal policies since every DM apply an identical policy, the set of policies can be extended to infinite product space and then we can consider the expected cost by integrating over $\prod_{i=N}^{\infty}(\mathbb{U}\times \mathbb{Y})$. Inequality \eqref{eq:4.71} follows from the fact that limsup is the greatest convergent subsequence limit for a bounded sequence, where we denoted the convergent subsequence of coordinates of policies in (Step 3) with $n \in \mathbb{I} \subset \mathbb{N}$.  Equality \eqref{eq:4.91} follows from the law of total expectation, and the dominated convergence theorem. 

Fix the convergent subsequence $n$, following from symmetry and Assumption \ref{assump:ind} and Assumption \ref{assump:oc}(i), we have $\beta_{n}^{i}=(u^{*, i}_{n},y^{i})$ are i.i.d. Now, using a similar argument as the proof of \cite[Theorem 8]{sanjari2018optimal}, through choosing a suitable subsubsequence and using the strong law of large numbers, we can show that for a continuous bounded function $g \in C_{b}(\mathcal{Z})$
\begin{flalign}\label{eq:091}
\mathbb{P}\bigg(\bigg\{\omega \in \Omega: \lim\limits_{n\rightarrow \infty}\bigg|\frac{1}{n}\sum_{i=1}^{n}g(\beta^{i}_{n})-\mathbb{E}(g(\beta^{i}_{\infty}))\bigg|=0\bigg\}\bigg)=1.
\end{flalign}

By considering a countable family of measure determining functions $\mathcal{T}\subset C_{b}(\mathcal{Z})$, \eqref{eq:091} implies that the empirical measures $\{\Lambda_{n}\}_{n}$ converges weakly to $\Lambda={\mathcal{L}}(\beta^{i}_{\infty})$ $\mathbb{P}$-almost surely, and $\Lambda$ is induced by the limit randomized policy $P^{*, \omega_{0}}$. We define the $w$-$s$ topology on the above set of probability measures on $(\Omega_{0}\times \mathbb{U}\times \mathbb{Y})$. That is, the coarsest topology on $\mathcal{P}(\Omega_{0} \times \mathbb{U} \times \mathbb{Y})$ under which $\int \hat{f}(\omega_{0}, u,y)\kappa(d\omega_{0}, du, dy): \mathcal{P}(\Omega_{0} \times \mathbb{U} \times \mathbb{Y}) \to \mathbb{R}$ is continuous for every measurable and bounded $\hat{f}$ which is continuous in $u$ but need not to be continuous in $y$ and $\omega_{0}$ (see e.g., \cite{Schal} and \cite[Theorem 5.6]{yuksel2018general}). Following from Assumption \ref{assump:c} and Assumption \ref{assump:oc}(ii), and since actions induced by identical policies are i.i.d. (thanks to symmetry), we have $\mathbb{P}$-almost surely
\begin{flalign*}
&\scalemath{0.93}{f_{n}:=\min{\bigg\{M, {c}\left(\omega_{0}, \cdot,\int_{\mathbb U}u \Lambda_{n}(du \times \mathbb{Y})\right)\bigg\}} \xrightarrow{\text{cont}} f:=\min{\bigg\{M, {c}\left(\omega_{0}, \cdot,\int_{\mathbb U}u \Lambda(du \times \mathbb{Y})\right)\bigg\}}}, 
\end{flalign*}
where we recall that the sequence $f_{n}$ converges continuously to $f$ ($f_{n}\xrightarrow{\text{cont}} f$) if and only if $f_{n}(a_{n})\to f(a)$ whenever $a_{n}\to a$ as $n \to \infty$. Now, \eqref{eq:4.101} follows from the generalized dominated convergence theorem for varying measures in \cite[Theorem 3.5]{serfozo1982convergence}. Equality \eqref{eq:4.111} follows from the monotone convergence theorem, and \eqref{eq:4.32c2} follows from the fact that the limit randomized policy, $P^{*, \omega_{0}}$, does not depend on $N$ and symmetry, hence, \eqref{eq:4.32c2} is true using a similar analysis as \eqref{eq:4.101}.  Inequality \eqref{eq:4.35r1} follows from the fact that the limit policy, $P^{*, \omega_{0}}(du^{i},dy^{i}):=\Pi^{*}(du^{i}|y^{i})f(\omega_{0}, y^{i})Q(dy^{i})$, achieving \eqref{eq:4.32c2} belongs to $L_{\text{PR,SYM}}$. That is because, following from (Step 3), for each DM, the set of randomized policies is closed under the product topology where each coordinate converges weakly, and hence, the limit policy is also a randomized policy induced by a subsequence of $N$-DM optimal policies (which are symmetric through DMs). This implies \eqref{eq:4.35r1} and completes the proof.}  
\end{itemize}

\section{Proofs from Section \ref{sec:4}}\label{appB}
\subsection{Independent measurement reduction under Assumption \ref{assump:ind1}}\label{sec:statreduc}
Under Assumption \ref{assump:ind1}(i), we can represent the expected cost  as
\begin{flalign}
\scalemath{0.95}{J_{N}(\underline{\gamma}^{1:N})}&\scalemath{0.95}{:=\int c(\omega_{0},u_{0:T-1}^{1},\dots,u_{0:T-1}^{N})\mu^{N}(d\omega_{0}, d\underline\zeta^{1:N})}\nonumber\\
&\:\scalemath{0.95}{\times \prod_{i=1}^{N}\prod_{t=0}^{T-1}1_{\{\gamma_{t}^{i}(y_{t}^{i})\in du_{t}^{i}\}}{{\nu}_{t}^{i}}\left(d{y}^{i}_{t}\middle|\omega_{0}, x_{0}^{1:N},{\zeta}^{1:N}_{0:t-1},{y}_{0:t-1}^{1:N},{u}_{0:t-1}^{1:N}\right)}\nonumber\\
&\scalemath{0.95}{=\int c(\omega_{0},u_{0:T-1}^{1},\dots,u_{0:T-1}^{N})\mu^{N}(d\omega_{0}, d\underline\zeta^{1:N})}\label{eq:stred}\\
&\:\scalemath{0.95}{\times \prod_{i=1}^{N}\prod_{t=0}^{T-1}1_{\{\gamma_{t}^{i}(y_{t}^{i})\in du_{t}^{i}\}}{{\psi}_{t}^{i}}\bigg({y}^{i}_{t}, \omega_{0}, x_{0}^{1:N},{\zeta}^{1:N}_{0:t-1},{y}_{0:t-1}^{1:N},{u}_{0:t-1}^{1:N}\bigg)\tau_{t}^{i}(dy_{t}^{i})}\nonumber\\
&\scalemath{0.95}{=\int c_{s}(\omega_{0},\underline{\zeta}^{1:N}, u_{0:T-1}^{1:N},y_{0:T-1}^{1:N})\mu^{N}(d\omega_{0}, d\underline\zeta^{1:N})\prod_{i=1}^{N}\prod_{t=0}^{T-1}1_{\{\gamma_{t}^{i}(y_{t}^{i})\in du_{t}^{i}\}}\tau_{t}^{i}(dy_{t}^{i})}\nonumber,
\end{flalign}
where the new (equivalent) cost function is defined as
\begin{flalign*}
\scalemath{0.95}{c_{s}(\omega_{0},\underline{\zeta}^{1:N}, u_{0:T-1}^{1:N},y_{0:T-1}^{1:N}):=c(\omega_{0},u_{0:T-1}^{1:N})\prod_{i=1}^{N}\prod_{t=0}^{T-1}{{\psi}_{t}^{i}}\bigg({y}^{i}_{t}, \omega_{0}, x_{0}^{1:N},{\zeta}^{1:N}_{0:t-1},{y}_{0:t-1}^{1:N},{u}_{0:t-1}^{1:N}\bigg)},
\end{flalign*}
and \eqref{eq:stred} follows from Assumption \ref{assump:ind1}(i). Similar derivation holds when randomized policies are considered. Similarly, we can define the new (equivalent) cost function under Assumption \ref{assump:ind1}(ii). We note that in the above, we considered control actions induced by deterministic policies; however, the above analysis can be extended to randomized policies by just replacing $\prod_{i=1}^{N}\prod_{t=0}^{T-1}1_{\{\gamma_{t}^{i}(y_{t}^{i})\in du_{t}^{i}\}}$ with $\prod_{i=1}^{N}\prod_{t=0}^{T-1}\gamma_{t}^{i}(du_{t}^{i}|y_{t}^{i})$.

\subsection{Proof of Lemma \ref{lem:3.1dy}}
We follow the steps of the proof of Lemma \ref{lem:exc}. For any permutation $\sigma \in S_{N}$, we define a randomized policy $P_{\pi}^{\sigma} \in \bar{L}^{N}$ as a permutation $\sigma$ of arguments of a randomized policy $P_{\pi} \in \bar{L}^{N}$, i.e., for $A^{i} \in \mathcal{B}(\Gamma^{i})$
\begin{flalign*}
&P_{\pi}^{\sigma}(\underline\gamma^{1}\in A^{1},\dots, \underline\gamma^{2}\in A^{N}):=P_{\pi}(\underline\gamma^{\sigma(1)}\in A^{1},\dots, \underline\gamma^{\sigma(N)}\in A^{N}).
\end{flalign*}
 We have
\begin{flalign}
&\scalemath{0.94}{\int P^{\sigma}_{\pi}(d\underline{\gamma})\mu^{N}(d\omega_{0},d\underline{\zeta})c^{N}( \underline{\zeta}, \underline{\gamma},\underline{y}, \omega_{0})\nu^{N}(d\underline{y}|\underline{\zeta}, \underline{\gamma},\omega_{0})}\nonumber\\
&\scalemath{0.94}{{=\int c(\omega_{0},\underline{u}^{1},\dots,\underline{u}^{N})\prod_{k=1}^{N}\underline{\gamma}^{k}(d\underline{u}^{k}|\underline{y}^{k})P_{\pi}(d\underline\gamma^{\sigma(1)},\dots,d\underline\gamma^{\sigma(N)})}}\label{eq:4.4.4.5dy}\\
&\:\scalemath{0.94}{{\times \tilde{\mu}^{N}(d\underline\zeta^{1:N}|\omega_{0})\mathbb{P}_{0}(d\omega_{0})\prod_{t=0}^{T-1}\prod_{i=1}^{N}{\nu_{t}^{i}}\left(d{y}^{i}_{t}\middle|\omega_{0}, x_{0}^{i},{\zeta}^{i}_{0:t-1},{y}_{0:t-1}^{1:N},{u}^{1:N}_{t-1}\right)}}\nonumber\\
&\scalemath{0.94}{=\int c(\omega_{0},\underline{u}^{\sigma(1)},\dots,\underline{u}^{\sigma(N)})\prod_{k=1}^{N}\underline{\gamma}^{\sigma(k)}(d\underline{u}^{\sigma(k)}|\underline{y}^{\sigma(k)})P_{\pi}(d\underline\gamma^{1},\dots,d\underline\gamma^{N})}\label{eq:4.4.4.4dy}\\
&\:\scalemath{0.94}{\times \tilde{\mu}^{N}(d(\underline{\zeta}^{\sigma})^{1:N}|\omega_{0})\mathbb{P}_{0}(d\omega_{0})\prod_{t=0}^{T-1}\prod_{i=1}^{N}{\nu_{t}^{i}}\left(d{y}^{\sigma(i)}_{t}\middle|\omega_{0}, x_{0}^{\sigma(i)},{\zeta}^{\sigma(i)}_{0:t-1},({y}_{0:t-1}^{\sigma})^{1:N},({u}^{\sigma}_{0:t-1})^{1:N}\right)}\nonumber\\
&\scalemath{0.94}{=\int c(\omega_{0},\underline{u}^{1},\dots,\underline{u}^{N})\prod_{k=1}^{N}\underline{\gamma}^{k}(d\underline{u}^{k}|\underline{y}^{k})P_{\pi}(d\underline\gamma^{1},\dots,d\underline\gamma^{N})}\label{eq:4.4.4.6dy}\\
&\:\scalemath{0.94}{\times \tilde{\mu}^{N}(d\underline\zeta^{1:N}|\omega_{0})\mathbb{P}_{0}(d\omega_{0})\prod_{t=0}^{T-1}\prod_{i=1}^{N}{\nu_{t}^{i}}\left(d{y}^{i}_{t}\middle|\omega_{0}, x_{0}^{i},{\zeta}^{i}_{0:t-1},{y}_{0:t-1}^{1:N},{u}^{1:N}_{0:t-1}\right)}\nonumber\\
&\scalemath{0.94}{=\int P_{\pi}(d\underline{\gamma})\mu^{N}(d\omega_{0},d\underline{\zeta})c^{N}( \underline{\zeta}, \underline{\gamma},\underline{y}, \omega_{0})\nu^{N}(d\underline{y}|\underline{\zeta}, \underline{\gamma},\omega_{0})}\nonumber
\end{flalign}
where $\tilde{\mu}^{N}$ is the conditional distribution of uncertainties $\underline{\zeta}^{1:N}$ given $\omega_{0}$, and \eqref{eq:4.4.4.5dy} follows from Assumption \ref{assump:sxcdy}(b) and the definition of randomized policy $P^{\sigma}_{\pi}$ and \eqref{eq:4.4.4.4dy} follows from relabeling $\underline{u}^{\sigma(i)}, \underline{y}^{\sigma(i)}, \underline{\zeta}^{\sigma(i)}$ with $\underline{u}^{i}, \underline{y}^{i}, \underline{\zeta}^{i}$ for all $i=1,\dots,N$ and the fact that  $y^{i}_{t}=h_{t}(x_{0}^{i},x_{0}^{-i}, \zeta^{i}_{0:t}, {\zeta}^{-i}_{0:t}, u_{0:t-1}^{i}, u_{0:t-1}^{-i})$. Equality \eqref{eq:4.4.4.6dy} follows from Assumption \ref{assump:sxcdy}(a), Assumption \ref{assump:3.1} and the hypothesis that the information structure is symmetric. The rest of the proof follows from similar steps in that of Lemma \ref{lem:exc}.

\subsection{Proof of Lemma \ref{lem:findefdyn2}}\label{section:prlem5.5}
We follow steps of the proof of Lemma \ref{lem:findef}. Under Assumption \ref{assump:ind2} and Assumption \ref{assump:c},  for every finite $N$, there exists an optimal policy in   $L_{\text{EX}}^{N}$. Consider a sequence $\{P_{\pi}^{*,N}\}_{N}$, where for every $N\geq 1$, $P_{\pi}^{*,N} \in L_{\text{EX}}^{N}$ and
\begin{flalign}
&{\int P_{\pi}^{*, N}(d\underline{\gamma})\mu^{N}(d\omega_{0},d\underline{\zeta})c^{N}( \underline{\zeta}, \underline{\gamma},\underline{y}, \omega_{0})\nu^{N}(d\underline{y}|\underline{\zeta}, \underline{\gamma},\omega_{0})}\nonumber\\
&{=\inf\limits_{P^{N}_{\pi} \in L_{\text{EX}}^{N}}\int P^{N}_{\pi}(d\underline{\gamma})\mu^{N}(d\omega_{0},d\underline{\zeta})c^{N}( \underline{\zeta}, \underline{\gamma},\underline{y}, \omega_{0})\nu^{N}(d\underline{y}|\underline{\zeta}, \underline{\gamma},\omega_{0})}\label{eq:32.1tdy}.
\end{flalign}

\begin{itemize}[wide]
\item[\textbf{(Step 1):}]
Let $(I_{1},I_{2},\dots)$ be i.i.d. random variables with the uniform distribution on the set $\{1,\dots,N\}$. For a fixed $N$ and for any $P_{\pi}^{*, N} \in  L^{N}_{\text{EX}}$, we construct  $P^{*,\infty}_{\pi, N} \in L_{\text{EX}}$ as follows:  for every fixed $N$ and for all $A^{i} \in \mathcal{B}(\Gamma^{i})$
\begin{flalign*}
&P^{*,\infty}_{\pi, N}(\underline\gamma^{1}\in A^{1},\dots, \underline\gamma^{2}\in A^{N}):=P_{\pi}^{*, N}(\underline\gamma^{I_{1}}\in A^{1},\dots, \underline\gamma^{I_{N}}\in A^{N}),
\end{flalign*}
where $P^{*,\infty}_{\pi, N}$ is the restriction of $P^{*,\infty}_{\pi, P_{\pi}^{N}}\in  L_{\text{EX}}$ to the first $N$ components. 

Let $u^{*, i}_{t, N}$ be the control action induced by $\gamma^{i}_{N, t}$ where random variables $(\gamma^{1}_{N, t},\dots,\gamma^{N}_{N, t})$ for all $t=0,\dots,T-1$ are determined by $P_{\pi}^{*, N}\in L^{N}_{\text{EX}}$. Let $u^{*, i}_{t, \infty, N}$ be the control action induced by $\gamma^{i}_{\infty, t}$ where random variables $(\gamma^{1}_{t,  \infty, N},\dots,\gamma^{N}_{t,  \infty, N})$ are determined by $P_{\pi, N}^{*, \infty}\in  L_{\text{EX}}$.  Let $\underline\gamma^{i}_{N}:=(\gamma^{i}_{N, 0},\dots,\gamma^{i}_{N, T-1})$, $\underline\gamma^{i}_{N, \infty}:=(\gamma^{i}_{0, \infty, N},\dots,\gamma^{i}_{T-1, \infty, N})$, $\underline u^{i}_{N}:=(u^{i}_{N, 0},\dots,u^{i}_{N, T-1})$ and $\underline u^{i}_{N, \infty}:=(u^{i}_{0, \infty, N},\dots,u^{i}_{T-1, \infty, N})$ for each DM. Since under the reduction (Assumption \ref{assump:ind1}), observations are i.i.d. through DMs and also independent of $\omega_{0}$, following from Theorem \ref{the:dia}, we have for every $m\geq 1$
\begin{flalign}
&\scalemath{0.94}{\bigg|\bigg|\mathcal{L}(\underline\gamma^{1}_{N},\dots,\underline\gamma^{m}_{N},\underline{y}^{1},\dots,\underline{y}^{m})-\mathcal{L}(\underline\gamma^{1}_{N, \infty},\dots,\underline\gamma^{m}_{N,\infty},\underline{y}^{1},\dots,\underline{y}^{m})\bigg|\bigg|}\nonumber\\
&\scalemath{0.94}{=\bigg|\bigg|\mathcal{L}(\underline\gamma^{1}_{N},\dots,\underline\gamma^{m}_{N})\prod_{i=1}^{m}\mathcal{L}(\underline{y}^{i})-\mathcal{L}(\underline\gamma^{1}_{N, \infty},\dots,\underline\gamma^{m}_{N,\infty})\prod_{i=1}^{m}\mathcal{L}(\underline{y}^{i})\bigg|\bigg|_{TV}\xrightarrow[N \to \infty]{}0}\label{eq:podA.7}.
\end{flalign}
where \eqref{eq:poA.7} follows from the fact that $(\underline\gamma^{1}_{N},\dots,\underline\gamma^{N}_{N})$ and $(\underline\gamma^{1}_{N, \infty},\dots,\underline\gamma^{N}_{N, \infty})$ are random variables with joint probability measures $P_{\pi}^{*, N} \in L_{\text{EX}}^{N}$ and $P_{\pi, N}^{*, \infty} \in L_{\text{EX}}\big|_{N}$, respectively. Since $\mathbb{U}$ is compact, and under the reduction the probability measure on observation is fixed, any joint probability measures on acttions and observations is tight, hence, $\{\mathcal{L}(\underline\gamma^{i}_{\infty, N})\}_{N}$ is tight for each DM and by exchangeablity $\mathcal{L}(\underline\gamma^{i}_{\infty, N})=\mathcal{L}(\underline\gamma^{1}_{\infty, N})$. Hence, we can find a subsequence such that $\mathcal{L}(\underline\gamma^{i}_{\infty, l})\xrightarrow[l\to \infty]{} \mathcal{L}(\underline\gamma^{i}_{\infty})$ for all $i\in \mathbb{N}$. Since marginals of $\{\mathcal{L}(\underline\gamma^{1}_{\infty, l}, \dots,\underline\gamma^{m}_{\infty, l})\}_{l}$ are tight, for each $m\geq 1$, there exists a further subsequence 
\begin{flalign*}
\mathcal{L}(\underline\gamma^{1}_{\infty, n}, \dots, \underline\gamma^{m}_{\infty, n}) \xrightarrow[n \to \infty]{} \mathcal{L}(\underline\gamma_{\infty}^{1}, \dots, \underline\gamma_{\infty}^{m}),
\end{flalign*}
where $(\underline\gamma_{\infty}^{1}, \underline\gamma_{\infty}^{2}, \dots)$ is infinitely-exchangeable and induced by $P_{\pi}^{*, \infty} \in L_{\text{EX}}$ since the set of infinitely-exchangeable random variables is closed under the weak-convergence topology. Hence, following from \eqref{eq:poA.7}, for each $m\geq 1$
\begin{flalign*}
\mathcal{L}(\underline\gamma^{1}_{n}, \dots, \underline\gamma^{m}_{n}) \xrightarrow[n \to \infty]{} \mathcal{L}(\underline\gamma_{\infty}^{1}, \dots, \underline\gamma_{\infty}^{m}).
\end{flalign*}
By construction $\underline u^{*, i}_{n}$ and $\underline u^{*, i}_{\infty}$ and since random variables $\underline\gamma^{i}_{n}$s are independent of $\underline{y}^{i}$s, we have  for each $m\geq 1$
\begin{flalign*}
(\underline u^{*, 1}_{n}, \dots,\underline u^{*, m}_{n}) \xrightarrow[n \to \infty]{\text{d}} (\underline u_{\infty}^{1}, \dots, \underline u_{\infty}^{m}),
\end{flalign*}
 where $(\underline u_{\infty}^{1},\underline u_{\infty}^{2}, \dots)$ is induced by an infinitely-exchangeable policies $P_{\pi}^{*, \infty} \in L_{\text{EX}}$.  
Following from Theorem \ref{the:ald}, $\mathbb{P}$-almost surely 
\begin{flalign}
F_{n, t}(A):=F_{n, t}^{\omega}(A):=\frac{1}{n}\sum_{i=1}^{n}\delta_{u^{*,i}_{n, t}(\omega)}(A)\xrightarrow[n \to \infty]{\text{d}}\alpha_{t}^{u, \omega}(A)\label{eq:emp1},
\end{flalign}
 where $A \in \mathcal{U}$ and $\omega$ denotes the sample path dependence and $\alpha_{t}^{u}$ is the directing random measure of an infinitely-exchangeable random variables $(\underline u_{\infty, t}^{1},\underline u_{\infty, t}^{2},\dots)$. By \eqref{eq:emp1}, since the action space is compact, for all $t=0,\dots,T-1$, we have $\mathbb{P}$-almost surely
\begin{flalign}
 \mu_{n, t}^{u}:=\mu_{n, t}^{u, \omega}:=\frac{1}{n}\sum_{i=1}^{n}u^{*,i}_{n, t}=\int_{\mathbb{U}}uF_{n, t}(du) \xrightarrow[n \to \infty]{\text{d}} \mu_{t}^{u}:=\int_{\mathbb{U}}u\alpha_{t}^{u, \omega}(du).
\end{flalign}

\item[\textbf{(Step 2):}]
 Let ${x}^{*, i}_{t, n}$ be the state of DM$^i$ at time $t$  under ${u}_{0:t-1,n}^{*, i}:=({u}_{0,n}^{*, i},\dots,{u}_{t-1,n}^{*, i})$: 
\begin{flalign}
&{x}^{*, i}_{t, n}=f_{t-1}\bigg({x}^{*, i}_{t-1, n},{u}^{*, i}_{n, t-1},\frac{1}{n}\sum_{p=1}^{n}{x}^{*, p}_{t-1, n}, \frac{1}{n}\sum_{p=1}^{n}{u}^{*, p}_{n, t-1}, w_{t-1}^{i}\bigg).\label{eq:Amfdynamics2}
\end{flalign}
Let $t=1$. We have
\begin{flalign}
&{x}^{*, i}_{1, n}=f_{0}\bigg({x}^{i}_{0},{u}^{*, i}_{n, 0},\frac{1}{n}\sum_{p=1}^{n}{x}^{p}_{0}, \frac{1}{n}\sum_{p=1}^{n}{u}^{*, p}_{n, 0}, w_{0}^{i}\bigg).
\end{flalign}
Since initial states are i.i.d. conditioned on $\omega_{0}$, by continuity of the function $f_{0}$ in actions and states, we have  ${x}^{*, i}_{1, n} \xrightarrow[n \to \infty]{\text{d}} {x}^{*, i}_{1, \infty}$ for all DMs. Hence, $\{\mathcal{L}(x^{*, 1}_{1, n}, \dots, x^{*, n}_{1, n})\}_{n}$ is tight, and hence, for each $m\geq 1$, there exists a subsubsequence $k$ such that $(x^{*, 1}_{1, k}, \dots, x^{*, m}_{1, k}) \xrightarrow[k \to \infty]{\text{d}} (x^{*, 1}_{1, \infty}, \dots, x^{*, m}_{1, \infty})$. Following from Theorem \ref{the:ald}, since $f_{0}$ is bounded, we have $\mathbb{P}$-almost surely
\begin{flalign}
\mu_{{k}, 1}^{x}:=\frac{1}{k}\sum_{i=1}^{n}x^{*,i}_{1, k}= \mu_{{k}, 1}^{x, \omega}=\int_{\mathbb{X}} x \frac{1}{{k}}\sum_{i=1}^{{k}}\delta_{x^{*,i}_{1, {k}}}(dx) \xrightarrow[k \to \infty]{\text{d}} \mu_{1}^{x}:=\int_{\mathbb{X}}x\alpha_{1}^{x, \omega}(dx),
\end{flalign}
where $\alpha_{1}^{x}$ is the directing measure for $(x^{*, 1}_{\infty, 1}, x^{*, 2}_{\infty, 1}, \dots)$.
Similarly, we can show that for $t=2$, 
\begin{flalign}
&{x}^{*, i}_{2, {k}}=f_{1}\bigg({x}^{*, i}_{1, {k}},{u}^{*, i}_{{k}, 1},\mu_{{k},1}^{x}, \mu_{{k},1}^{u}, w_{1}^{i}\bigg).
\end{flalign}
By continuity of the function $f_{1}$ and the analysis for $t=1$, we have  ${x}^{*, i}_{2, {k}} \xrightarrow[k \to \infty]{\text{d}} {x}^{*, i}_{2, \infty}$ for all DMs. Hence, $\{\mathcal{L}(x^{*, 1}_{2, {k}}, \dots, x^{*, {k}}_{2, {k}})\}_l$ is tight and for each $m\geq 1$, there exists a further subsubsequence $k_{l}$ such that $(x^{*, 1}_{2, {k_{l}}}, \dots, x^{*, m}_{2, {k_{l}}}) \xrightarrow[{k_{l}} \to \infty]{\text{d}} (x^{*, 1}_{2, \infty}, \dots, x^{*, m}_{2, \infty})$ . Following from Theorem \ref{the:ald}, since $f_{1}$ is bounded, , we have $\mathbb{P}$-almost surely
\begin{flalign}
 \mu_{{k_{l}}, 2}^{x}:= \mu_{{k_{l}}, 2}^{x, \omega}=\int_{\mathbb{X}} x \frac{1}{{k_{l}}}\sum_{i=1}^{{k_{l}}}\delta_{x^{*,i}_{2, {k_{l}}}}(dx) \xrightarrow[{k_{l}} \to \infty]{\text{d}} \mu_{2}^{x}:=\int_{\mathbb{X}}x\alpha_{2}^{x, \omega}(dx),
\end{flalign}
where $\alpha_{2}^{x}$ is the directing measure for $(x^{*, 1}_{\infty, 2}, x^{*, 2}_{\infty, 2}, \dots)$.
 By induction, for each $m\geq 1$, there exists a further subsubsequence $n$ (which we indicate by $n$ to omit further sub-subscript) such that $(\underline x^{*, 1}_{n}, \dots, \underline{x}^{*, m}_{n}) \xrightarrow[n_{{l}} \to \infty]{\text{d}} (\underline x^{*, 1}_{\infty}, \dots,\underline x^{*, m}_{\infty})$  and $\mu_{n_{l}, t}^{x}\xrightarrow[n \to \infty]{\text{d}}\mu_{t}^{x}$  for all $t=0,\dots,T-1$.
 
 Now, we follow the steps of Lemma  \ref{lem:findef}, however, in addition to actions and observations, we consider states and disturbances in our analysis and we use the result of (Step 2). Define $\tilde{P}^{*, n}$ as the joint probability measures of $(\underline u^{*,1}_{n},\underline x^{*,1}_{n}, \mu_{n, 0:T-1}^{u}, \mu_{n, 0:T-1}^{x}, \underline{y}, \underline{\zeta})$. Since marginals on $(\underline u^{*,1}_{n},\underline x^{*,1}_{n}, \mu_{n, 0:T-1}^{u}, \mu_{n, 0:T-1}^{x})$ are tight and under the reduction marginals on  $(\underline{y}, \underline{\zeta})$ are fixed, $\{\tilde{P}^{*, n}\}_{n}$ is tight. Hence, there exists a further subsubsequence $\{\tilde{P}^{*, n_{k}}\}_{n_{k}}$ converges weakly to $\tilde{P}^{*}$ as $n_{k}$ goes to infinity. This implies that marginals $\{\tilde{P}^{*, n_{k}}\}_{n_{k}}$ converge to the marginals of  $\tilde{P}^{*}$, hence, $\tilde{P}^{*}$ is induced by $(\underline u_{\infty}^{1},\underline u_{\infty}^{2},\dots)$ which is infinitely-exchangeable and is  induced by a policy in $L_{\text{EX}}$.
\item[\textbf{(Step 3):}] Since the cost function is continuous in states and actions, under the reduction (Assumption \ref{assump:ind2}), we have $\mathbb{P}$-almost surely
\begin{flalign}
&\scalemath{0.9}{\frac{1}{N}\sum_{i=1}^{N}\sum_{t=0}^{T-1}c\bigg(\omega_{0}, x_{t}^{i},u_{t}^{i},\frac{1}{N}\sum_{p=1}^{N}u_{t}^{p},\frac{1}{N}\sum_{p=1}^{N}x_{t}^{p}\bigg)}\nonumber\\
&\scalemath{0.9}{\times\prod_{i=1}^{N}\prod_{t=0}^{T-1}{{\phi}_{t}^{i}}\bigg({y}^{i}_{t}, \omega_{0}, x_{0}^{i},{\zeta}^{i}_{0:t-1},{y}_{0:t-1}^{i},{u}_{0:t-1}^{i}, \frac{1}{N}\sum_{p=1}^{N}u^{p}_{0:t-1}, \frac{1}{N}\sum_{p=1}^{N}x^{p}_{0:t}\bigg)}\nonumber\\
&\scalemath{0.9}{=\frac{1}{N}\sum_{i=1}^{N}\bar{c}\bigg(\omega_{0},\underline{\zeta}^{i}, \underline{x}^{i}, \underline{u}^{i},\frac{1}{N}\sum_{p=1}^{N}\underline{u}^{p},\frac{1}{N}\sum_{p=1}^{N}\underline{x}^{p}\bigg)}\label{eq:star12di}\\
&\scalemath{0.9}{\prod_{i=1}^{N}{{\underline\phi}^{i}}\bigg(\underline{y}^{i}, \omega_{0},\underline{\zeta}^{i},\underline{u}^{i},\frac{1}{N}\sum_{p=1}^{N}\underline{u}^{p}, \frac{1}{N}\sum_{p=1}^{N}\underline{x}^{p}\bigg)}\nonumber,
\end{flalign}
where \eqref{eq:star12di} is true following from \eqref{eq:mfdynamics2} and Assumption \ref{assump:c} for some function $\hat{c}:\Omega_{0}\times {\bf S} \times {\bf X} \times {\bf U} \times {\bf U} \times {\bf X} \to \mathbb{R}_{+}$ which is continuous in states and actions and 
\begin{flalign*}
&\scalemath{0.9}{{{\underline\phi}^{i}}\bigg(\underline{y}^{i}, \omega_{0},\underline{\zeta}^{i},\underline{u}^{i},\frac{1}{N}\sum_{p=1}^{N}\underline{u}^{p}, \frac{1}{N}\sum_{p=1}^{N}\underline{x}^{p}\bigg)}\\
&\scalemath{0.9}{:=\prod_{t=0}^{T-1}{{\phi}_{t}^{i}}\bigg({y}^{i}_{t}, \omega_{0}, x_{0}^{i},{\zeta}^{i}_{0:t-1},{y}_{0:t-1}^{i},{u}_{0:t-1}^{i}, \frac{1}{N}\sum_{p=1}^{N}u^{p}_{0:t-1}, \frac{1}{N}\sum_{p=1}^{N}x^{p}_{0:t}\bigg)}
\end{flalign*}
 We have,
\begin{flalign}
&\scalemath{0.92}{\limsup\limits_{N \to \infty}\inf\limits_{P_{\pi}^{N} \in L^{N}_{\text{EX}}}\int P_{\pi}^{N}(d\underline{\gamma})\mu^{N}(d\omega_{0},d\underline{\zeta})c^{N}( \underline{\zeta}, \underline{\gamma},\underline{y}, \omega_{0})\nu(d\underline{y}|\underline{\zeta}, \underline{\gamma},\omega_{0})}\nonumber\\
&\scalemath{0.92}{=\limsup\limits_{N \to \infty}\int \int_{\prod_{i=n_{l}+1}^{\infty}{\bf{Y}\times {\bf{S}}}}\bar{c}\bigg(\omega_{0},\underline{\zeta}^{i}, \underline{x}^{i}, \underline{u}^{i},\mu_{N, 0:T-1}^{u},\mu_{N, 0:T-1}^{x}\bigg)}\label{eq:B.17}\\
&\:\scalemath{0.92}{\times \tilde{P}^{*, N}(d\underline u^{*,i},d\underline x^{*,i}, d\mu_{N, 0:T-1}^{u}, d\mu_{N, 0:T-1}^{x},\underline{y}, \underline{\zeta})}\nonumber\\
&\scalemath{0.92}{\times\prod_{i=1}^{\infty}{{\underline\phi}^{i}}\bigg(\underline{y}^{i}, \omega_{0},\underline{\zeta}^{i},\underline{u}^{*, i},\mu_{N, 0:T-1}^{u}, \mu_{N, 0:T-1}^{x}\bigg)\mathbb{P}_{0}(d\omega_{0})}\nonumber\\
&\scalemath{0.92}{\geq\lim\limits_{n_{k} \to \infty}\int \int_{\prod_{i=n_{k}+1}^{\infty}{\bf{Y}\times {\bf{S}}}} \bar{c}\bigg(\omega_{0},\underline{\zeta}^{1}, \underline{x}^{1}, \underline{u}^{1},\mu_{n_{k}, 0:T-1}^{u},\mu_{n_{k}, 0:T-1}^{x}\bigg)}\label{eq:B.18}\\
&\:\scalemath{0.92}{\times \tilde{P}^{*, n_{k}}(d\underline u^{*,i},d\underline x^{*,i}, d\mu_{n_{k}, 0:T-1}^{u}, d\mu_{n_{k}, 0:T-1}^{x}, \underline{y}, \underline{\zeta})}\nonumber\\
&{\times\prod_{i=1}^{\infty}{{\underline\phi}^{i}}\bigg(\underline{y}^{i}, \omega_{0},\underline{\zeta}^{i},\underline{u}^{*, i},\mu_{n_{k}, 0:T-1}^{u}, \mu_{n_{k}, 0:T-1}^{x}\bigg)\mathbb{P}_{0}(d\omega_{0})}\nonumber\\
&\scalemath{0.92}{=\int \bar{c}\bigg(\omega_{0},\underline{\zeta}^{1}, \underline{x}^{1}, \underline{u}^{1},\mu_{0:T-1}^{u},\mu_{0:T-1}^{x}\bigg)}\scalemath{0.92}{\tilde{P}^{*}(d\underline u^{*,i},d\underline x^{*,i}, d\mu_{0:T-1}^{u}, d\mu_{0:T-1}^{x}, \underline{y}, \underline{\zeta})}\label{eq:B.19}\\
&\scalemath{0.94}{\times \prod_{i=1}^{\infty}{{\underline\phi}^{i}}\bigg(\underline{y}^{i}, \omega_{0},\underline{\zeta}^{i},\underline{u}^{*, i},\mu_{0:T-1}^{u}, \mu_{0:T-1}^{x}\bigg)\mathbb{P}_{0}(d\omega_{0})}\nonumber\\
&\scalemath{0.92}{\geq\limsup\limits_{N \to \infty}\inf\limits_{P_{\pi} \in L_{\text{EX}}}\int P_{\pi, N}(d\underline{\gamma})\mu^{N}(d\omega_{0},d\underline{\zeta})c^{N}( \underline{\zeta}, \underline{\gamma},\underline{y}, \omega_{0})\nu(d\underline{y}|\underline{\zeta}, \underline{\gamma},\omega_{0})}\label{eq:B.20}.
\end{flalign}
where \eqref{eq:B.17} follows from integrating over the set $(\prod_{i=n_{l}+1}^{\infty}{\bf{Y}\times {\bf{S}}})$ and the fact that under the reduction, observations and disturbances, initial states are i.i.d. and $(u_{N}^{*, 1},\dots, u_{N}^{*, N})$ is $N$-exchangeable. Inequality \eqref{eq:B.18} follows from the assumption that the cost function is bounded and limsup is the greatest subsequence limit of a bounded sequence. Equality \eqref{eq:B.19} follows from the dominated convergence theorem and following from Assumption \ref{assump:c} and Assumption \ref{assump:ind2} and since probability measures on observations disturbances are fixed and since by (Step 2) $\{\tilde{P}^{*, n_{k}}\}_{n_{k}}$ converges weakly to $\tilde{P}^{*}$ and $\prod_{i=1}^{\infty}{{\underline\phi}^{i}}(\underline{y}^{i}, \omega_{0},\underline{\zeta}^{i},\underline{u}^{*, i},\mu_{n_{k}, 0:T-1}^{u}, \mu_{n_{k}, 0:T-1}^{x})$ converges weakly to the limit in the product topology as $n_{k}$ goes to infinity. Inequality \eqref{eq:B.20} follows from the fact that $\tilde{P}^{*}$ is the joint measure with the first coordinate $(u_{\infty}^{1},u_{\infty}^{2},\dots)$ which is infinitely-exchangeable and  is induced by a policy in $L_{\text{EX}}$. The above inequalities are equalities since the opposite direction is true (that is because $L_{\text{EX}}\big|_{N}\subset L_{\text{EX}}^{N}$) and this completes the proof.
\end{itemize}
  \subsection{Proof of Theorem \ref{the:exmftdy2}}
We complete the proof in five steps where the steps are similar to the proof of Theorem \ref{the:2}.
\begin{itemize}[wide]
\item[\textbf{(Step 1):}] Under Assumption \ref{assump:c} and Assumption \ref{assump:ind2}, by Lemma \ref{lem:3.1dy}, for every finite $N$, there exists a randomized optimal policy in   $L_{\text{EX}}^{N}$. Consider a sequence $\{P_{\pi}^{*,N}\}_{N}$, where for every $N\geq 1$, $P_{\pi}^{*,N} \in L_{\text{EX}}^{N}$ and
\begin{flalign}
&\int P_{\pi}^{*, N}(d\underline{\gamma})\mu^{N}(d\omega_{0},d\underline{\zeta})c^{N}( \underline{\zeta}, \underline{\gamma},\underline{y}, \omega_{0})\nu^{N}(d\underline{y}|\underline{\zeta}, \underline{\gamma},\omega_{0})\label{eq:32.1dy}\\
&=\inf\limits_{P^{N}_{\pi} \in L_{\text{EX}}^{N}}\int P^{N}_{\pi}(d\underline{\gamma})\mu^{N}(d\omega_{0},d\underline{\zeta})c^{N}( \underline{\zeta}, \underline{\gamma},\underline{y}, \omega_{0})\nu^{N}(d\underline{y}|\underline{\zeta}, \underline{\gamma},\omega_{0})\nonumber.
\end{flalign}
\item[\textbf{(Step 2):}]
 Similar to (Step 2) of the proof of Theorem \ref{the:2} using Lemma \ref{lem:findefdyn2} and Theorem \ref{the:defin}, we can show that to complete the proof, it is sufficient to show 
\begin{flalign}
&\lim\limits_{M \to \infty}\limsup\limits_{N \to \infty}\inf\limits_{P_{\pi}^{N} \in L_{\text{PR,SYM}}^{N}}\int P_{\pi}^{N}(d\underline{\gamma})\mu^{N}(d\omega_{0},d\underline{\zeta})\min{\{M, c^{N}(\underline{\zeta}, \underline{\gamma},\underline{y}, \omega_{0})\}\nu^{N}(d\underline{y}|\underline{\zeta}, \underline{\gamma},\omega_{0})}\nonumber\\
&\geq \inf\limits_{P_{\pi} \in L_{\text{PR,SYM}}}\limsup\limits_{N \to \infty} \int P_{\pi, N}(d\underline{\gamma})\mu^{N}(d\omega_{0},d\underline{\zeta})c^{N}( \underline{\zeta}, \underline{\gamma},\underline{y}, \omega_{0})\nu^{N}(d\underline{y}|\underline{\zeta}, \underline{\gamma},\omega_{0})\label{eq:4.4dyf}.
\end{flalign}

In the next two steps, we justify \eqref{eq:4.4dyf} through showing that there exists a subsequence of randomized policies induced by symmetric/identical private randomization whose weak subsequent limit achieves the right hand side of \eqref{eq:4.4dyf}.
\item[\textbf{(Step 3):}] Consider the set of randomized policies  $L_{\text{PR,SYM}}^{N}$. We note that under a symmetric information structure and since each DM applies an identical policy, under Assumption \ref{assump:ind1}, $\underline{y}^{i}$ are i.i.d. through DMs and also independent of $\omega_{0}$. Hence, following from the information structure, the randomized policy spaces of each DM is separated from the policies of the other decision makers. Hence, we can equivalently represent any privately randomized policy for each DM acting through time separately as a probability measures  induced by symmetric (identical randomized policies), i.e.,  probability measures $q$ on $({\bf{U}}\times {\bf{Y}})$ where randomized policies of each DM for every $t=0,\dots,T-1$ satisfy
\begin{flalign}
&\scalemath{0.95}{\int g(\omega_{0}, x_{0}^{i}, \zeta_{0:t-1}^{i}, y_{0:t}^{i},u_{0:t}^{i})q(dy^{i}_{0:t},du^{i}_{0:t}|\omega_{0}, x_{0}^{i}, \zeta_{0:t-1}^{i})}\nonumber\\
  &\scalemath{0.95}{=\int g(\omega_{0}, x_{0}^{i}, \zeta_{0:t-1}^{i},y_{0:t}^{i},u_{0:t}^{i})}{\prod_{k=0}^{t} \Pi_{k}^{N}(du^{i}_{k}|y^{i}_{k})\eta_{k}(dy^{i}_{k}|\omega_{0}, x_{0}^{i}, \zeta_{0:k-1}^{i}, y^{i}_{0:k-1}, u_{0:k-1}^{i})}\nonumber,
\end{flalign}
for all bounded functions $g$ which is continuous  in actions and observations and measurable in other arguments and for some stochastic kernel $\Pi_{k}^{N}$ representing a randomized policy of DMs at time $k$ (which is identical through DMs).

Since $\mathbb{U}$ is compact, the marginals on $\bf{U}$ is relatively compact under the weak convergence topology. Hence, the collection of all probability measures with these relatively compact marginals are also relatively compact (see e.g., \cite[Proof of Theorem 2.4]{yukselSICON2017}). Since every DM applies an identical policy and since observations are i.i.d., the randomized policy space is relatively compact (where each coordinate is relatively compact in the weak convergence topology), and hence, there exists a subsequence of randomized policies $\tilde{q}_{n}\in \mathcal{P}(\prod_{i}({\bf{Y}} \times {\bf{U}}))$  converges weakly (each coordinate converges weakly) to a limit $\tilde{q}$ (as an infinite product of policies of DMs), where $n$ is the index of the subsequence and $n$ goes to infinity. Now, we show that randomized policy space is closed under the weak convergence topology. Assume $\hat{q}_{n}\in \mathcal{P}({\bf{Y}} \times {\bf{U}})$ (induced by identical randomized policies $\Pi_{t}^{n}$ for each DM at time $t=0,\dots, T-1$) converges weakly to $\hat{q}$. If Assumption \ref{assump:ind1}(i) (under the structure Assumption \ref{assump:ind2}) holds, then there exists an independent static reduction for each DM through time, and hence,  following from the discussion in the proof of \cite[Theorem 5.2]{yuksel2018general}, each coordinate of policy spaces corresponds to DM$^{i}$ at time $t$ is closed under the weak convergence topology. Also, if  Assumption \ref{assump:ind1}(ii)  (under the structure \ref{assump:ind2}) holds, then \cite[Theorem 5.6]{yuksel2018general} leads to the same conclusion. Hence, this implies that for $\tilde{q}^{*}_{N}\in \mathcal{P}(\prod_{i=1}^{N}({\bf{Y}} \times {\bf{U}}))$ induced by optimal randomized policies $\Pi^{*,N}_{t}$ for each DM at time $t$, there exists a subsequence $\tilde{q}_{n}^{*}\in \mathcal{P}(\prod_{i=1}^{\infty}({\bf{Y}} \times {\bf{U}}))$ (as an infinite product of policies of DMs  $\Pi^{*,n}_{t}$) converges weakly (each coordinate converges weakly) to a limit $\tilde{q}^{*}$ which is in $L_{\text{PR,SYM}}$ and it is induced by a randomized policy $\Pi_{t}^{*,\infty}$ for each DM at time $t$.

 \item[\textbf{(Step 4):}]
 Let $\{\hat{q}^{*}_{N}\}_{N}$ be a policy for each DM induced by optimal randomized policies $\Pi^{*,N}_{t}$ for $N$-DM team problems, and let $\underline{u}^{i,*}_{N}:=({u}^{i,*}_{N, 0},\dots,{u}^{i,*}_{N, T-1})$ be the action of DM$^{i}$ through time induced by $\Pi^{*,N}_{t}$. Following from (Step 3), there exists a weak subsequential limit $\hat{q}^{*}$ of $\{\hat{q}^{*}_{n}\}_{n}$ as $n \to \infty$ for each DM, which is induced by  $\Pi^{*,\infty}_{t}$. Let $\underline{u}^{i,*}_{\infty}:=({u}^{i,*}_{\infty, 0},\dots,{u}^{i,*}_{\infty, T-1})$ be the action of DM$^{i}$ induced by the identically randomized policy $\Pi^{*,\infty}_{t}$. Define
\begin{equation}\label{eq:empdystate}
\Upsilon_{N}(B):=\frac{1}{N}\sum_{i=1}^{N}\delta_{(\underline{x}^{i}_{N}, \tilde\alpha_{N}^{i})}(B),
\end{equation}
where $\tilde\alpha_{N}^{i}:=(\underline{u}^{i, *}_{N},\underline{y}^{i}, \underline{\zeta}^{i})$, $B \in {\bf{X}} \times \mathcal{Z}:={\bf{X}}\times {\bf{U}}\times {\bf{Y}} \times {\bf{S}}$, ${\bf{U}}:=(\prod_{t=0}^{T-1}\mathbb{U})$, ${\bf{Y}}:=(\prod_{t=0}^{T-1}\mathbb{Y})$, ${\bf{S}}:=(\prod_{t=0}^{T-1}{\mathbb{S}})=\mathbb{X}\times (\prod_{t=0}^{T-1}{\mathbb{W}\times \mathbb{V}})$, ${\bf X}:=(\prod_{t=0}^{T-1}\mathbb{X})$, $\underline{y}^{i}:=(y_{0}^{i},\dots,y_{T-1}^{i})$, $\underline{\zeta}^{i}:=(\zeta^{i}_{0},\dots,\zeta_{T-1}^{i})$, and $\underline{x}^{i}_{N}:=(x_{0}^{i},\dots,x_{T-1}^{i})$ with states are driven by a sequence of $N$-DM randomized optimal policies of $\Pi^{*,N}_{t}$. In the following, we show that, conditioned on $\omega_{0}$,   the subsequence of empirical measures $\{\Upsilon_{n}\}_{n}$ converges to $\Upsilon:=\mathcal{L}((\underline{x}^{1}_{\infty}, \tilde\alpha_{\infty}^{1})|\omega_{0})$ in $w$-$s$ topology, where $\tilde\alpha^{i}_{\infty}=(\underline{u}^{*, i}_{\infty},\underline{y}^{i}, \underline{\zeta}^{i})$ and $\underline{x}^{i}_{\infty}$ denotes that states of DM$^{i}$ driven by $\underline{u}^{*, i}_{\infty}$  (we note that the convergence is weakly, but since $\underline{\zeta}^{i}$s are exogenous random variables with a fixed marginal, the convergence is also in the $w$-$s$ topology). Define 
\begin{equation}\label{eq:empdy}
\bar{Q}_{N}(B):=\frac{1}{N}\sum_{i=1}^{N}\delta_{\tilde\alpha^{i}_{N}}(B),
\end{equation}
 where $B \in \mathcal{Z}$. Under the reduction in Assumption \ref{assump:ind1}, observations of each DM are independent of actions and observations of other DMs through time, $t=0,\dots,T-1$, and hence, a similar argument, which is used to show \eqref{eq:091}, implies that the subsequence of empirical measures $\{\bar{Q}_{n}\}_{n \in \mathbb{I}}$  converges $\mathbb{P}$-almost surely to $\bar{Q}=\mathcal{L}(\tilde\alpha^{i}_{\infty}|\omega_{0})$ in $w$-$s$ topology. Define 
\begin{equation}\label{eq:empdystate2}
\Upsilon_{n}^{t}(A):=\frac{1}{n}\sum_{i=1}^{n}\delta_{(\underline{x}^{i}_{t, n}, \tilde\alpha_{t, n}^{i})}(A),
\end{equation}
where  $\tilde\alpha_{t, n}^{i}:=({u}^{i, *}_{n, t},{y}^{i}_{t}, {\zeta}^{i}_{t})$, $A \in \mathbb{X} \times \mathbb{U} \times \mathbb{Y} \times \mathbb{S}$. Since conditioned on $\omega_{0}$, initial states are i.i.d, the empirical measure of initial states converges weakly to $\mathcal{L}(x^{1}_{0}|\omega_{0})$ $\mathbb{P}$-almost surely. Since $\{\bar{Q}_{n}\}_{n}$ converges $\mathbb{P}$-almost surely to $\bar{Q}$ in $w$-$s$ topology, we can conclude that $\Upsilon_{N}^{0}$ converges $\Upsilon^{0}:=\mathcal{L}(({x}^{i}_{0}, \tilde\alpha_{0, \infty}^{i})|\omega_{0})$ in $w$-$s$ topology $\mathbb{P}$-almost surely. Following from \eqref{eq:mfdynamics2}, for $t=0$, we have for every continuous and bounded function  $g\in C_{b}(\mathbb{X})$, conditioned on $\omega_{0}$, $\mathbb{P}$-almost surely
\begin{flalign}
\scalemath{0.9}{\lim\limits_{n \to \infty}\frac{1}{n}\sum_{i=1}^{n}g(x_{1,n}^{i})}&\scalemath{0.9}{=\lim\limits_{n \to \infty}\frac{1}{n}\sum_{i=1}^{n}g\bigg(f_{0}\bigg(x_{0}^{i},u^{*, i}_{n, 0},\frac{1}{n}\sum_{p=1}^{n}x_{0}^{p}, \frac{1}{n}\sum_{p=1}^{n}u_{n,0}^{*,p}, w_{0}^{i}\bigg)\bigg)}\nonumber\\
&\scalemath{0.9}{=\lim\limits_{n \to \infty}\int g\bigg(f_{0}\bigg(x,u,\int x \Upsilon_{n}^{0}(dx \times \mathbb{U} \times \mathbb{Y} \times \mathbb{S}), \int u \Upsilon_{n}^{0}(\mathbb{X}\times  du \times \mathbb{Y} \times \mathbb{S}) ,\zeta\bigg)\bigg)}\label{eq:5.67}\\
&\scalemath{0.9}{\times \Upsilon_{n}^{0}(dx,du,dy,d\zeta)}\nonumber\\
&\scalemath{0.9}{=\int g\bigg(f_{0}\bigg(x,u,\int x \Upsilon^{0}(dx \times \mathbb{U} \times \mathbb{Y} \times \mathbb{S}), \int u \Upsilon^{0}(\mathbb{X}\times  du \times \mathbb{Y} \times \mathbb{S}) ,\zeta\bigg)\bigg)}\label{eq:5.68}\\
&\scalemath{0.9}{\times\Upsilon^{0}(dx,du,dy,d\zeta)}\nonumber
\end{flalign}
where \eqref{eq:5.67} follows from \eqref{eq:empdystate2}, and \eqref{eq:5.68} follows from the generalized dominated convergence theorem for varying measures. That is because, function $g$ is continuous and bounded, $f_{0}$ is a bounded function which is continuous in actions and observations and measurable in uncertainties, and the fact that under the reduction, conditioned on $\omega_{0}$,  $\Upsilon_{N}^{0}$ converges $\Upsilon^{0}:=\mathcal{L}(({x}^{i}_{0}, \tilde\alpha_{0,\infty}^{i})|\omega_{0})$ in $w$-$s$ topology $\mathbb{P}$-almost surely. Hence, since $\{\bar{Q}_{n}\}_{n}$ converges $\mathbb{P}$-almost surely to $\bar{Q}$ in $w$-$s$ topology conditioned on $\omega_{0}$,  $\Upsilon_{N}^{1}$ converges $\Upsilon^{1}:=\mathcal{L}(({x}^{i}_{1,\infty},\tilde\alpha_{1,\infty}^{i})|\omega_{0})$ in $w$-$s$ topology $\mathbb{P}$-almost surely. By induction, one can show that  conditioned on $\omega_{0}$,  $\Upsilon_{N}^{t}$ converges $\Upsilon^{t}:=\mathcal{L}(({x}^{i}_{t,\infty}, \tilde\alpha_{t,\infty}^{i})|\omega_{0})$ in $w$-$s$ topology $\mathbb{P}$-almost surely for $t=0,\dots, T-1$. Hence, conditioned on $\omega_{0}$, $\{\Upsilon_{n}\}_{n \in \mathbb{I}}$ converges to $\Upsilon:=\mathcal{L}((\underline{x}^{i}_{\infty}, \tilde\alpha_{\infty}^{i})|\omega_{0})$ in $w$-$s$ topology.

 \item[\textbf{(Step 5):}]
 By Assumption \ref{assump:c}, similar to the proof of Theorem \ref{lem:findefdyn2}, we have \eqref{eq:star12di}. Under the reduction, we can consider policy spaces for each DM individually. Let for every $t=0,\dots,T-1$, $P^{*, \omega_{0}}_{n}$ be a probability measure on actions, observations and uncertainties induced by optimal randomized policies for each DM (which is identical because of symmetry) for $N$-DM teams conditioned on $\omega_{0}$, i.e., a probability measure that satisfies
\begin{flalign}
&\scalemath{0.95}{\int g(\omega_{0}, x_{0}^{i}, \zeta_{0:t-1}^{i}, y_{0:t}^{i},u_{n, 0:t}^{i, *})P^{*, \omega_{0}}_{n}(dx_{0}^{i}, d\zeta_{0:t-1}^{i}, dy^{i}_{0:t},du^{i, *}_{n, 0:t}|\omega_{0})}\nonumber\\
  &\scalemath{0.95}{=\int g(\omega_{0}, x_{0}^{i}, \zeta_{0:t-1}^{i},y_{0:t}^{i},u_{n, 0:t}^{i, *})\mu^{i}(dx_{0}^{i}, d\zeta_{0:t-1}^{i}|\omega_{0})}\label{eq:eq18dyn}\\
  &\:\:\:\:\:\:\scalemath{0.95}{ \times \prod_{k=0}^{t} \Pi_{k}^{*, n}(du^{*, i}_{n, k}|y^{i}_{k})\eta_{k}(dy^{i}_{k}|\omega_{0}, x_{0}^{i}, \zeta_{0:k-1}^{i}, y^{i}_{0:k-1}, u_{n, 0:k-1}^{i, *})}\nonumber,
\end{flalign}
for all bounded functions $g$ which is continuous in actions and observations and measurable in other arguments. Similarly, we denote $P^{*, \omega_{0}}$ as a probability measure induced by the limit policy, i.e., a probability measure satisfying \eqref{eq:eq18dyn} induced by $ \Pi_{k}^{*, \infty}$. Hence, following from a similar argument as in the (Step 4) of the proof of Theorem \ref{the:exmftdy2}, we have
\begin{flalign}
&\scalemath{0.9}{\lim\limits_{M \to \infty}\limsup\limits_{N \to \infty}\inf\limits_{P_{\pi}^{N} \in L_{\text{PR,SYM}}^{N}}\int P_{\pi}^{N}(d\underline{\gamma})\mu^{N}(d\omega_{0},d\underline{\zeta})\nu(d\underline{y}|\underline{\zeta}, \underline{\gamma},\omega_{0})\min{\{M, c^{N}(\underline{\zeta},\underline{\gamma},\underline{y}, \omega_{0})\}}}\nonumber\\
&\scalemath{0.9}{\geq  \lim\limits_{M\rightarrow \infty}\lim\limits_{n\rightarrow \infty}\int\int \min\bigg\{M,}\scalemath{0.9}{\bar{c}\left(\omega_{0},\underline\zeta, \underline x, \underline u,\int \underline u\Upsilon_{n}({\bf{X}} \times d \underline u \times {\bf{Y}}\times {\bf{S}}),\int \underline x {\Upsilon}_{n}(d\underline x \times {\bf{U}} \times {\bf{Y}}\times {\bf{S}})\right)\bigg\}}\label{eq:35.5.12}\\
&\scalemath{0.9}{\times {\Upsilon}_{n}(d \underline x, d\underline u,d\underline y,d\underline\zeta)\prod_{i=1}^{\infty}P_{n}^{*,\omega_{0}}(d\underline{u}^{i, *}_{n},d\underline{y}^{i},d\underline{\zeta}^{i})\prod_{i=1}^{\infty}{{\underline\phi}^{i}}\bigg(\underline{y}^{i}, \omega_{0},\underline{\zeta}^{i},\underline{u}^{i, *}_{n},\frac{1}{n}\sum_{p=1}^{n}\underline{u}^{p, *}_{n}, \frac{1}{n}\sum_{p=1}^{n}\underline{x}^{p}_{n}\bigg)\mathbb{P}_{0}(d\omega_{0})}\nonumber\\
&\scalemath{0.9}{=  \lim\limits_{M\rightarrow \infty}\int \lim\limits_{n\rightarrow \infty}\int \min\bigg\{M,}\scalemath{0.9}{\bar{c}\left(\omega_{0},\underline\zeta,\underline x,\underline u,\int \underline u\Upsilon_{n}({\bf{X}} \times d \underline u \times {\bf{Y}}\times {\bf{S}}),\int \underline x {\Upsilon}_{n}(d\underline x \times {\bf{U}} \times {\bf{Y}}\times {\bf{S}})\right)\bigg\}}\label{eq:34.5.12}\\
&\scalemath{0.9}{\times {\Upsilon}_{n}(d \underline x, d\underline u,d\underline y,d\underline\zeta)\prod_{i=1}^{\infty}P_{n}^{*,\omega_{0}}(d\underline{u}^{i, *}_{n},d\underline{y}^{i},d\underline{\zeta}^{i})\prod_{i=1}^{\infty}{{\underline\phi}^{i}}\bigg(\underline{y}^{i}, \omega_{0},\underline{\zeta}^{i},\underline{u}^{i, *}_{n},\frac{1}{n}\sum_{p=1}^{n}\underline{u}^{p, *}_{n}, \frac{1}{n}\sum_{p=1}^{n}\underline{x}^{p}_{n}\bigg)\mathbb{P}_{0}(d\omega_{0})}\nonumber\\
&\scalemath{0.9}{=\lim\limits_{M\rightarrow \infty}\int\int \min\bigg\{M,}\scalemath{0.9}{\bar{c}\left(\omega_{0},\underline\zeta,\underline x,\underline u,\int \underline u\Upsilon({\bf{X}} \times d\underline u \times {\bf{Y}}\times {\bf{S}}),\int \underline x {\Upsilon}(d\underline x \times {\bf{U}} \times {\bf{Y}}\times {\bf{S}})\right)\bigg\}}\label{eq:33.5.12}\\
&\scalemath{0.9}{\times {\Upsilon}(d \underline x, d\underline u,d\underline y,d\underline\zeta)\prod_{i=1}^{\infty}P^{*,\omega_{0}}(d\underline{u}^{i, *}_{\infty},d\underline{y}^{i},d\underline{\zeta}^{i})\prod_{i=1}^{\infty}{{\underline\phi}^{i}}\bigg(\underline{y}^{i}, \omega_{0},\underline{\zeta}^{i},\underline{u}^{i, *}_{\infty},{E}[\underline{u}^{1, *}_{\infty}|\omega_{0}], {E}[\underline{x}^{1}_{\infty}|\omega_{0}]\bigg)\mathbb{P}_{0}(d\omega_{0})}\nonumber\\
&\scalemath{0.9}{=\int\int \bar{c}\left(\omega_{0},\underline\zeta,\underline x,\underline u,\int \underline u\Upsilon({\bf{X}} \times d\underline u \times {\bf{Y}}\times {\bf{S}}),\int \underline x {\Upsilon}(d\underline x \times {\bf{U}} \times {\bf{Y}}\times {\bf{S}})\right)}\label{eq:37.5.12}\\
&\scalemath{0.9}{\times {\Upsilon}(d \underline x, d\underline u,d\underline y,d\underline\zeta)\prod_{i=1}^{\infty}P^{*,\omega_{0}}(d\underline{u}^{i, *}_{\infty},d\underline{y}^{i},d\underline{\zeta}^{i})\prod_{i=1}^{\infty}{{\underline\phi}^{i}}\bigg(\underline{y}^{i}, \omega_{0},\underline{\zeta}^{i},\underline{u}^{i, *}_{\infty},{E}[\underline{u}^{1, *}_{\infty}|\omega_{0}], {E}[\underline{x}^{1}_{\infty}|\omega_{0}]\bigg)\mathbb{P}_{0}(d\omega_{0})}\nonumber\\
&\scalemath{0.9}{\geq \inf\limits_{P_{\pi} \in L_{\text{PR,SYM}}}\limsup\limits_{N \to N} \int P_{\pi, N}(d\underline{\gamma})\mu^{N}(d\omega_{0},d\underline{\zeta})c^{N}( \underline{\zeta}, \underline{\gamma},\underline{y}, \omega_{0})\nu(d\underline{y}|\underline{\zeta}, \underline{\gamma},\omega_{0})}\nonumber,
\end{flalign}
 where \eqref{eq:35.5.12} follows from \eqref{eq:empdystate},  \eqref{eq:star12di}, and since limsup is the greatest convergent subsequence limit for a bounded sequence, and \eqref{eq:34.5.12} follows from the dominated convergence theorem. Following from a similar argument as the analysis in (Step 4) of the proof of Theorem \ref{the:2}, since $\{{\Upsilon}_{n}\}_{n \in \mathbb{I}}$ converges weakly to ${\Upsilon}$ $\mathbb{P}$-almost surely, an argument based on the generalized dominated convergence theorem for varying measures in \cite[Theorem 3.5]{serfozo1982convergence} implies \eqref{eq:33.5.12}, and \eqref{eq:37.5.12} follows from  the monotone convergence theorem. Hence, \eqref{eq:4.4dyf} holds and this completes the proof. 
 \end{itemize}
\section{Proofs from Section \ref{sec:aproxi}}\label{appC}
 \subsection{Proof of Theorem \ref{the:3}}
\begin{itemize}[wide]
\item [(i)]
We first show \eqref{eq:approx}. We have
\begin{flalign}
\scalemath{0.95}{\inf\limits_{P_{\pi}^{N} \in L_{\text{CO}}^{N}}\int}& \scalemath{0.95}{P_{\pi}^{N}(d\underline{\gamma})\mu^{N}(d\omega_{0},d\underline{y})c^{N}(\underline{\gamma},\underline{y}, \omega_{0})}\nonumber\\ &\scalemath{0.95}{\geq \inf\limits_{P_{\pi}^{N} \in L_{\text{CO}}^{N}\cap L_{\text{EX}}\big{|}_{N}}\int P_{\pi}^{N}(d\underline{\gamma})\mu^{N}(d\omega_{0},d\underline{y})c^{N}(\underline{\gamma},\underline{y}, \omega_{0})-\epsilon_{N}}\label{eq:ap4.2}\\
&\scalemath{0.95}{=\inf\limits_{P_{\pi}^{N} \in L_{\text{PR,SYM}}^{N}}\int P_{\pi}^{N}(d\underline{\gamma})\mu^{N}(d\omega_{0},d\underline{y})c^{N}(\underline{\gamma},\underline{y}, \omega_{0})-\epsilon_{N}}\label{eq:ap4.4},
\end{flalign}
where $L_{\text{EX}}\big{|}_{N}$ denotes the set of $N$-DM randomized policies which are the restrictions of policies in $L_{\text{EX}}$ to the $N$ first components. By Lemma \ref{lem:exc} since $ L_{\text{CO}}^{N}$ is convex, without losing global optimality, we can optimize over $L_{\text{CO}}^{N}\cap L_{\text{EX}}^{N}$. Let $\epsilon >0$, and consider $P_{\pi, \epsilon}^{*,N} \in L_{\text{CO}}^{N}\cap L^{N}_{\text{EX}}$ such that
\begin{flalign}
\inf\limits_{P^{N}_{\pi}\in L_{\text{CO}}^{N}\cap L_{\text{EX}}^{N}}&\int P^{N}_{\pi}(d\underline{\gamma})\mu^{N}(d\omega_{0},d\underline{y})c^{N}(\underline{\gamma},\underline{y}, \omega_{0})\geq \int P_{\pi,\epsilon}^{*, N}(d\underline{\gamma})\mu^{N}(d\omega_{0},d\underline{y})c^{N}(\underline{\gamma},\underline{y}, \omega_{0})-\epsilon\label{eq:asy32.1}.
\end{flalign}
Following from the proof of Lemma \ref{lem:findef}, using $P_{\pi,\epsilon}^{*, N} \in L_{\text{CO}}^{N}\cap L^{N}_{\text{EX}}$ and by considering the indexes as a sequence of i.i.d. random variables with uniform distribution on the set $\{1,\dots,N\}$, we can construct an infinitely-exchangeable policy $P_{\pi, \epsilon}^{*,\infty}$ where the restriction of an infinitely-exchangeable policy to $N$ first components $P^{*,\infty}_{\pi, N, \epsilon} \in L_{\text{CO}}^{N}\cap L_{\text{EX}}\big{|}_{N}$, satisfies 
\begin{flalign}
&\int P_{\pi, N, \epsilon}^{*,\infty}(d\underline{\gamma})\mu^{N}(d\omega_{0},d\underline{y})c^{N}(\underline{\gamma},\underline{y}, \omega_{0}) \int P_{\pi,\epsilon}^{*, N}(d\underline{\gamma})\mu^{N}(d\omega_{0},d\underline{y}) c^{N}(\underline{\gamma},\underline{y}, \omega_{0})+\epsilon_{N}\label{eq:asymendif}.
\end{flalign}
  Hence, \eqref{eq:asy32.1} and \eqref{eq:asymendif} imply that
\begin{flalign*}
\scalemath{0.95}{\inf\limits_{P^{N}_{\pi}\in L_{\text{CO}}^{N}\cap L_{\text{EX}}^{N}}}&\scalemath{0.95}{\int P^{N}_{\pi}(d\underline{\gamma})\mu^{N}(d\omega_{0},d\underline{y})c^{N}(\underline{\gamma},\underline{y}, \omega_{0})}\\
&\scalemath{0.95}{\geq \inf\limits_{P^{N}_{\pi}\in L_{\text{CO}}^{N}\cap  L_{\text{EX}}\big{|}_{N}}\int P^{N}_{\pi}(d\underline{\gamma})\mu^{N}(d\omega_{0},d\underline{y})c^{N}(\underline{\gamma},\underline{y}, \omega_{0})-\epsilon-\epsilon_{N}}.
\end{flalign*}
 Since $\epsilon$ is arbitrary, this implies \eqref{eq:ap4.2}. By Theorem \ref{the:defin}, without losing optimality, we can optimize over $L_{\text{CO,SYM}}^{N}$. Equality \eqref{eq:ap4.4} is true since $L_{\text{CO,SYM}}^{N}$ is convex with extreme points in  $L_{\text{PR,SYM}}^{N}$, and the map $\int P_{\pi}^{N}(d\underline{\gamma})\mu^{N}(d\omega_{0},d\underline{y})c^{N}(\underline{\gamma},\underline{y}, \omega_{0}):L_{\text{CO,SYM}}^{N} \to \mathbb{R}$ is linear. 
 
  Now, we show \eqref{eq:approx1} holds. We have
\begin{flalign}
\scalemath{0.95}{\inf\limits_{P_{\pi}^{N} \in L_{\text{D}}^{N}}\int P_{\pi}^{N}(d\underline{\gamma})\mu^{N}(d\omega_{0},d\underline{y}) c^{N}(\underline{\gamma},\underline{y}, \omega_{0})}&\scalemath{0.95}{=\inf\limits_{P_{\pi}^{N} \in L_{\text{PR}}^{N}}\int P_{\pi}^{N}(d\underline{\gamma})\mu^{N}(d\omega_{0},d\underline{y})c^{N}(\underline{\gamma},\underline{y}, \omega_{0})}\label{eq:ap14.23}\\
&\scalemath{0.9}{\geq \inf\limits_{P_{\pi}^{N} \in L_{\text{PR,SYM}}^{N}}\int P_{\pi}^{N}(d\underline{\gamma})\mu^{N}(d\omega_{0},d\underline{y})c^{N}(\underline{\gamma},\underline{y}, \omega_{0})-\epsilon_{N}}\label{eq:ap14.4},
\end{flalign}
where \eqref{eq:ap14.23} follows from Blackwell's irrelevant information theorem \cite{Blackwell2} and since $L_{\text{CO}}^{N}$ is convex with extreme  points in $L_{\text{PR}}^{N}$ and the map $\int P_{\pi}^{N}(d\underline{\gamma})\mu^{N}(d\omega_{0},d\underline{y})c^{N}(\underline{\gamma},\underline{y}, \omega_{0}):L_{\text{CO}}^{N} \to \mathbb{R}$ is linear, hence, without losing optimality, we can optimaize over $L_{\text{CO}}^{N}$. Inequality  \eqref{eq:ap14.4} follows from \eqref{eq:approx} and this completes the proof of (i).
\item [(ii)] Let $P_{\pi}^{*}\in L_{\text{PR,SYM}}$ be an optimal policy of  ($\mathcal{P}_{\infty}$) and $P_{\pi,N}^{*}$ is the restriction of $P_{\pi}^{*}$ to the first $N$ components. Define for all $N\in \mathbb{N}$
\begin{flalign*}
&a_{N}:= \int P_{\pi,N}^{*}(d\underline{\gamma})\mu^{N}(d\omega_{0},d\underline{y})c^{N}(\underline{\gamma},\underline{y}, \omega_{0})\\
&b_{N}:=\inf\limits_{P_{\pi}^{N} \in L_{\text{PR,SYM}}^{N}}\int P_{\pi}^{N}(d\underline{\gamma})\mu^{N}(d\omega_{0},d\underline{y})c^{N}(\underline{\gamma},\underline{y}, \omega_{0}).
\end{flalign*}
Following from (Step 4) of the proof of Theorem \ref{the:2}, since the cost function is bounded, 
\begin{flalign}
\scalemath{0.93}{\limsup\limits_{N \to \infty}\int P_{\pi,N}^{*}(d\underline{\gamma})}&\scalemath{0.93}{\mu^{N}(d\omega_{0},d\underline{y})c^{N}(\underline{\gamma},\underline{y}, \omega_{0})}\nonumber\\&\scalemath{0.93}{= \limsup\limits_{N \to \infty}\inf\limits_{P_{\pi}^{N} \in L_{\text{PR,SYM}}^{N}}\int P_{\pi}^{N}(d\underline{\gamma})\mu^{N}(d\omega_{0},d\underline{y})c^{N}(\underline{\gamma},\underline{y}, \omega_{0})}\label{eq:L23}.
\end{flalign}
Hence, $\limsup\limits_{N \to \infty}a_{N}=\limsup\limits_{N \to \infty}b_{N}$. Following from (Step 4) of the proof of Theorem \ref{the:2}, and symmetry, $\lim\limits_{N \to \infty}a_{N}=a<\infty$ and also there exists a subsequence such that $\lim\limits_{k \to \infty}b_{N_{k}}=a<\infty$. On the other hand, since $a_{N}\geq b_{N}$ for all $N \in \mathbb{N}$, we can find $\tilde{\epsilon}_{N} \geq 0$ such that $a_{N}=b_{N}+ \tilde{\epsilon}_{N}$. Taking limit as $k$ goes to infinity from both sides, we have $a=\lim\limits_{{k} \to \infty}(b_{N_{k}}+\epsilon_{N_{k}})=a+\lim\limits_{{k} \to \infty}\epsilon_{N_{k}}$.
Hence, $\lim\limits_{{k} \to \infty}\epsilon_{N_{k}}=0$ since $\tilde\epsilon_{N}\geq 0$.
Hence, there exists $\bar\epsilon_{N}\geq 0$ where $\bar\epsilon_{N} \to 0$ as $N$ goes to infinity such that 
\begin{flalign}
\int P_{\pi,N}^{*}(d\underline{\gamma})&\mu^{N}(d\omega_{0},d\underline{y})c^{N}(\underline{\gamma},\underline{y}, \omega_{0})\nonumber\\
&\leq\inf\limits_{P_{\pi}^{N} \in L_{\text{D}}^{N}}\int P_{\pi}^{N}(d\underline{\gamma})\mu^{N}(d\omega_{0},d\underline{y}) c^{N}(\underline{\gamma},\underline{y}, \omega_{0})+\epsilon_{N}+\bar{\epsilon}_{N}\label{eq:approxl}
\end{flalign}
where \eqref{eq:approxl} follows from \eqref{eq:approx1}, and this completes the proof of (ii).   
\end{itemize}
\end{appendix}
\bibliographystyle{plain}

\end{document}